\documentclass[preprint,11pt]{article}

\usepackage{graphicx}
\usepackage{amsmath,amssymb}
\usepackage{enumerate}
\usepackage[numbers]{natbib}
\usepackage{lscape}

\newtheorem{theo}{Theorem}[section]
\newtheorem{corol}[theo]{Corollary}
\newtheorem{lem}[theo]{Lemma}
\newtheorem{prop}[theo]{Proposition}
\newtheorem{exam}[theo]{Example}

\begin{document}
\begin{center}
{\Large \bf A computational algebraic geometry approach to classify partial Latin rectangles}
\end{center}

\begin{center}
{\large \em R. M. Falc\'on}
\vspace{0.5cm}

{\small School of Building Engineering, University of Seville, Spain.\\
{\em rafalgan@us.es}}
\end{center}

\vspace{0.5cm}

\noindent {\large \bf Abstract.} This paper provides an in-depth analysis of how computational algebraic geometry can be used to deal with the problem of counting and classifying $r\times s$ partial Latin rectangles based on $n$ symbols of a given size, shape, type or structure. The computation of Hilbert functions and triangular systems of radical ideals enables us to solve this problem for all $r,s,n\leq 6$. As a by-product, explicit formulas are determined for the number of partial Latin rectangles of size up to six. We focus then on the study of non-compressible regular partial Latin squares and their equivalent incidence structure called seminet, whose distribution into main classes is explicitly determined for point rank up to eight. We prove in particular the existence of two new configurations of point rank eight.

\vspace{0.5cm}

\noindent{\bf Keywords:} Partial Latin rectangle, seminet, polynomial ring, ideal.\\
\noindent{\bf 2000 MSC:} 05B15, 05B25, 13F20.

\section{Introduction}\label{intro}

An {\em $r\times s$ partial Latin rectangle based on $[n]=\{1,\ldots,n\}$} is an $r \times s$ array $P$ in which each cell is either empty or contains one symbol chosen from $[n]$, such that each symbol occurs at most once in each row and in each column. Its {\em size} is the number of non-empty cells. If there are not empty cells, then $P$ is a {\em Latin rectangle}. If $r=s=n$, then $P$ is a {\em partial Latin square} of order $n$ (a {\em Latin square} if there are not empty cells). Hereafter, $\mathcal{R}_{r,s,n}$ and $\mathcal{R}_{r,s,n:m}$ denote, respectively, the set of $r\times s$ partial Latin rectangles based on $[n]$ and its subset of elements of size $m$.

\vspace{0.1cm}

The problem of counting $r\times s$ Latin rectangles based on $n$ symbols is a classical problem in combinatorial design theory that is currently solved for $r,s,n\leq 11$ (see \cite{Stones2010} and the references therein). Their distribution into isotopism, isomorphism and main classes is only known for Latin squares of order $n\leq 11$ \cite{Hulpke2011, McKay2005, McKay2007}. Nevertheless, the more general problem of counting and classifying partial Latin rectangles in the sets $\mathcal{R}_{r,s,n}$ and $\mathcal{R}_{r,s,n:m}$ has not yet been dealt with in depth. It is only known the cardinality and distribution into isotopism and isomorphism classes of $\mathcal{R}_{r,s,n}$ for $r,s,n\leq 6$ \cite{FalconStones2015}, and the cardinality of $\mathcal{R}_{r,s,n:m}$ for $r,s,n\leq 4$ \cite{Falcon2013, Falcon2015}. This paper contributes to this line of research and provides an in-depth analysis of how computational algebraic geometry can be used to enumerate and classify partial Latin rectangles according not only to their size, but also to their shape, type and structure. The implementation of this algebraic method in the study of non-compressible regular partial Latin squares also enable us to deal with the equivalent problem of classifying seminets, a type of incident structure introduced by U{\v{s}}an \cite{Uvsan1977} as a natural generalization of nets.

\vspace{0.1cm}

The remainder of the paper is organized as follows. Section 2 deals with some preliminary concepts and results on partial Latin rectangles, seminets and computational algebraic geometry. These results are implemented in Section 3 to determine the cardinality of $\mathcal{R}_{r,s,n:m}$ for all $r,s,n\leq 6$. In Section 4, the distribution of non-empty cells per row and column and the number of occurrences of each symbol enable us to use computational algebraic geometry in order to identify the set of partial Latin rectangles of a given shape, type or structure. The distribution of $\mathcal{R}_{r,s,n}$ into isotopism and main classes is then determined for all $r,s,n\leq 6$. As a by-product, we establish explicit formulas for the number of partial Latin rectangles of any order and size up to six. Finally, Section 5 deals with the distribution into main classes of seminets of point rank up to eight. We also prove the existence of two new configurations of seminets with point rank eight that complete the classification given by Lyakh \cite{Lyakh1988}.

\section{Preliminaries}\label{sec:preliminaries}

We review in this section some basic results on partial Latin rectangles, seminets and computational algebraic geometry that are used throughout the paper. We refer to the monographs of D\'enes and Keedwell \cite{Denes1991} and Cox et al. \cite{Cox2007} and to the original paper of U{\v{s}}an \cite{Uvsan1977} for more details about these topics.

\subsection{Classification of partial Latin rectangles}

An {\em entry} of a partial Latin rectangle $P\in\mathcal{R}_{r,s,n}$ is a triple $(i,j,k)\in [r]\times [s]\times [n]$ that is uniquely related to a non-empty cell of $P$ which is situated in the $i^{th}$ row and $j^{th}$ column and contains the symbol $k$. The set of entries of $P$ is denoted as $E(P)$. Let $S_m$ denote the symmetric group on $m$ elements. An {\em isotopism} of $\mathcal{R}_{r,s,n}$ is any triple $\Theta=(\alpha,\beta,\gamma)\in S_r\times S_s\times S_n$, where $\alpha$, $\beta$ and $\gamma$ constitute, respectively, a permutation of the rows, columns and symbols of any partial Latin rectangle $P\in \mathcal{R}_{r,s,n}$. This gives rise to the {\em isotopic} partial Latin rectangle $P^{\Theta}\in\mathcal{R}_{r,s,n}$, whose set of entries is $E(P^{\Theta})=\{(\alpha(i),\beta(j),\gamma(k))\colon\, (i,j,k)\in E(P)\}$. Permutations among the components of the entries of $P$ also give rise to new partial Latin rectangles. The {\em parastrophic} partial Latin rectangle of $P$ according to a permutation $\pi\in S_3$ is denoted by $P^{\pi}$ and has as set of entries the set $E(P^{\pi})=\{(p_{\pi(1)},p_{\pi(2)},p_{\pi(3)})\colon\, (p_1,$ $p_2,p_3)\in E(P)\}$. If the permutation $\pi$ preserves the set $\mathcal{R}_{r,s,n}$, then $\pi$ is said to be a {\em parastrophism}. The set of parastrophisms of $\mathcal{R}_{r,s,n}$ is, therefore,
\begin{itemize}
\item $\{\mathrm{Id}\}$ if $r$, $s$ and $n$ are pairwise distinct.
\item $\{\mathrm{Id},(12)\}$ if $r=s\neq n$.
\item $\{\mathrm{Id},(13)\}$ if $r=n\neq s$.
\item $\{\mathrm{Id},(23)\}$ if $s=n\neq r$.
\item $S_3$ if $r=s=n$.
\end{itemize}
Two partial Latin rectangles are {\em paratopic} if one of them is isotopic to a parastrophic partial Latin rectangle of the other. To be isotopic, parastrophic or paratopic are equivalence relations among partial Latin rectangles. They make possible the respective distribution of partial Latin rectangles into {\em isotopism}, {\em parastrophism} and {\em paratopism} or {\em main} classes.

\subsection{Compressibility and regularity of partial Latin squares}

Let $P$ be a partial Latin square of order $n$. It is said to be {\em non-compressible} if this does not contain empty rows or empty columns, or if
all the $n$ symbols appear as entries in $E(P)$. The partial Latin square $P$ is said to be {\em regular} if the next two conditions hold.
\begin{itemize}
\item It does not contain a cell that is, simultaneously, the only non-empty cell in its row and the only non-empty cell in its column.
\item If there exists a row or a column with exactly one non-empty cell, then the symbol contained in this cell appears at least twice in $E(P)$.
\end{itemize}

\subsection{Seminets}

Bates \cite{Bates1947} defined a {\em halfnet} as an incidence structure of points and lines such that there exist three distinct {\em parallel classes} of lines, every point is on at most one line of each class and any two lines belonging to distinct classes meet in at most one point. The number of points constitutes the {\em point rank} of a halfnet. Two halfnets are in the same {\em isomorphism class} if there exists a permutation among the points that preserves collinearity in each parallel class. If this happens after relabeling their parallel classes, then they are in the same {\em main class}. Currently, the distribution of halfnets into isomorphism and main classes is only partially known for nets and, to a much lesser extent, seminets.

\vspace{0.1cm}

Bruck \cite{Bruck1951} defined a {\em net} of order $n$ as a halfnet of $n^2$ points and $3n$ lines in which every point is on exactly one line of each parallel class, any two lines meet in exactly one point and there exists at least one line with exactly $n$ distinct points. Hence, every line contains $n$ points and every parallel class is formed by $n$ lines. More recently and motivated by its application in coding theory, U{\v{s}}an \cite{Uvsan1977} introduced the concept of {\em seminet} as a halfnet in which every point is on exactly one line of each parallel class and any two lines meet in at most one point. Unlike nets, the lines of a seminet can contain different numbers of points and its parallel classes can have different numbers of lines. The {\em $L$-order} of a seminet is the maximum number of lines in a parallel class. If all the lines have the same number $n$ of points, then all the parallel classes have the same number $m$ of lines. In this case, the seminet is said to be {\em $n$-regular}. If, furthermore, $m=n$, then it is a net of order $n$.

\vspace{0.1cm}

\begin{figure*}[ht]
  \begin{center}
  \begin{tabular}{ccc}
   \begin{tabular}{c} \includegraphics[width=.3\textwidth]{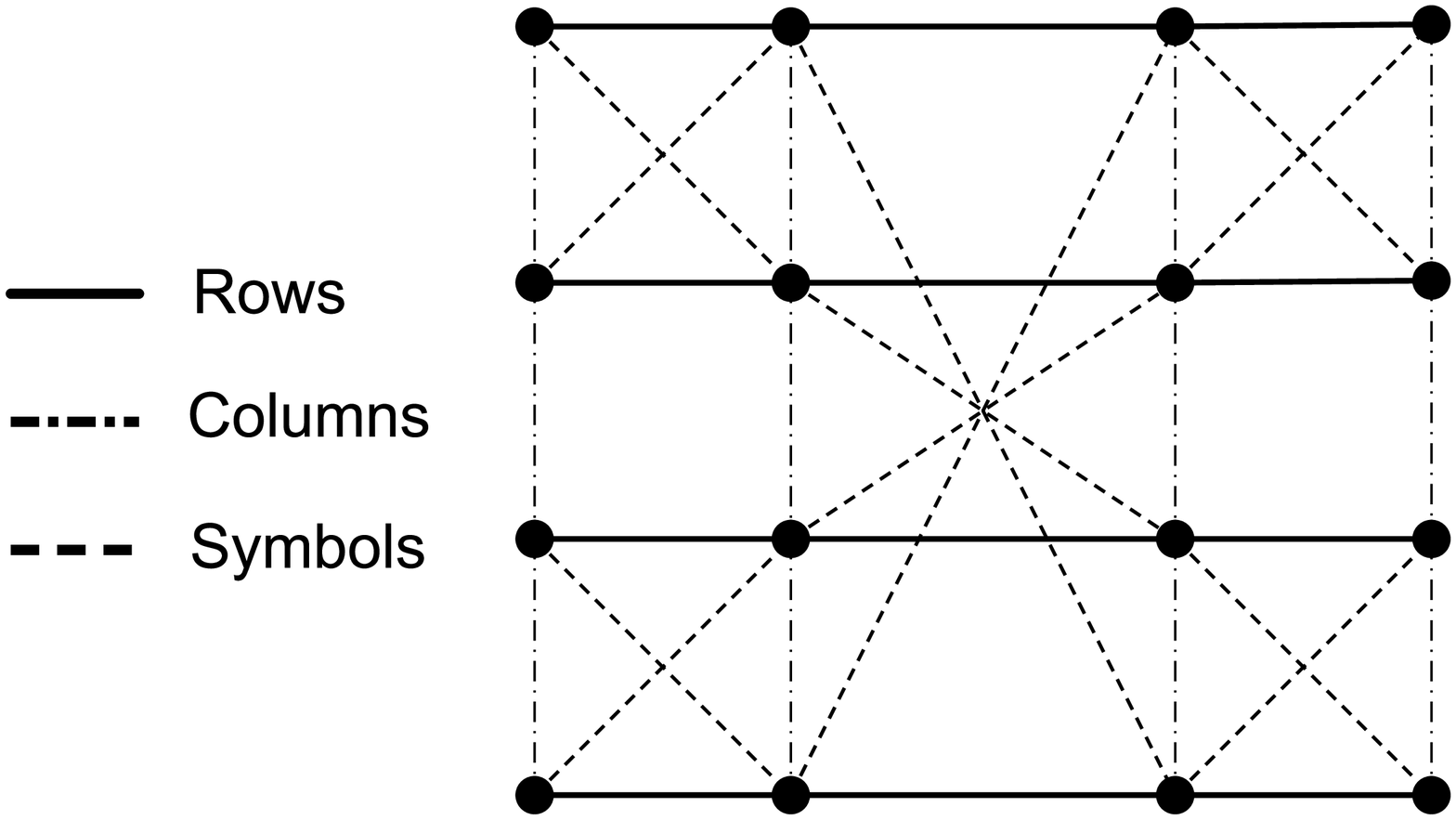}
    \end{tabular} & $\equiv$ & {\small \begin{tabular}{|c|c|c|c|}\hline
   1 & 2 & 3 & 4\\ \hline
   2 & 1 & 4 & 3\\ \hline
   3 & 4 & 1 & 2\\ \hline
   4 & 3 & 2 & 1\\ \hline
   \end{tabular}}
  \end{tabular}
  \end{center}
  \caption{Net identified with a Latin square of order $4$.}
  \label{Fig1}
\end{figure*}

\begin{figure}[ht]
  \begin{center}
  \begin{tabular}{ccc}
   \begin{tabular}{c} \includegraphics[width=.4\textwidth]{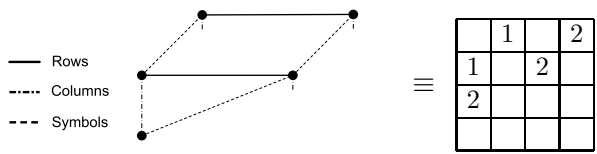}
    \end{tabular} & $\equiv$ & {\small \begin{tabular}{|c|c|c|c|}\hline
   \  & 1 & \  & 2\\ \hline
   1 & \  & 2 & \ \\ \hline
   2 &  \ &  \ & \ \\ \hline
   \  & \  & \  & \ \\ \hline
   \end{tabular}}
  \end{tabular}
  \end{center}
  \caption{Seminet identified with a partial Latin square of order $4$ and size $5$.}
  \label{Fig2}
\end{figure}

Every net of order $n$ can be identified with a Latin square of the same order. The points and parallel classes of the net are respectively identified with the cells of the Latin square and its sets of cells sharing the same row, column or symbol (see Figure \ref{Fig1}). In addition, Stojakovi{\'c} and U{\v{s}}an \cite{Stojakovic1979} proved that every seminet of $L$-order $n$ can be identified with a non-compressible regular partial Latin square of order $n$ in a similar way that nets do with Latin squares. In this case, the points of the seminet are identified with the non-empty cells of the partial Latin square (see Figure \ref{Fig2}). As a consequence, the distribution of nets and seminets into isomorphism and main classes results, respectively, from the equivalent distribution of Latin squares and non-compressible regular partial Latin squares into isotopism and main classes.

\renewcommand{\tabcolsep}{4pt}

Havel \cite{Havel1985} defined a {\em configuration} as a seminet containing at least four points such that every line contains at least two points and any two points $P$ and $Q$ of the seminet are {\em connected}, that is to say, there exists a sequence of points and lines, $P_0,l_0,P_1,l_1,\ldots,P_m$, such that $P_0=P$, $P_m=Q$ and each pair of points $P_{i-1}$ and $P_i$ are on the line $l_{i-1}$, for all $i\leq m$. Havel determined the main classes of those configurations with point rank up to seven and, shortly after, Lyakh \cite{Lyakh1988} gave a classification of those configurations with point rank eight.

\subsection{Computational algebraic geometry}

Let $X$ and $\mathbb{K}[X]$ respectively be the ordered set of $n$ variables $\{x_1,\ldots,x_n\}$ and the related multivariate polynomial ring $\mathbb{K}[x_1,\ldots,x_n]$ over a base field $\mathbb{K}$. The {\em class} of a polynomial $p\in \mathbb{K}[X]$ is the minimum $i\leq n$ such that $p\in \mathbb{K}[x_1,\ldots,x_i]$. A {\em triangular system} in $\mathbb{K}[X]$ is a finite ordered set of polynomials $\{p_1,\ldots,p_m\}\subset \mathbb{K}[X]$ such that the class of $p_i$ is less than the class of $p_{i+1}$, for all $i<m$.

\vspace{0.1cm}

An {\em ideal} of polynomials in $\mathbb{K}[X]$ is any subset $I\subseteq \mathbb{K}[X]$ such that $0\in I$; $p+q\in I$, for all $p,q\in I$; and $pq\in I$ for all $p\in I$ and $q\in \mathbb{K}[X]$. A {\em subideal} of $I$ is any subset $J\subseteq I$ that is also an ideal in $\mathbb{K}[X]$. The ideal {\em generated by} a finite set of polynomials $\{p_1,\ldots,p_m\}\subset \mathbb{K}[X]$ is defined as the set $\{q_1p_1+\ldots+q_np_n\colon\, q_1,\ldots,q_n\in\mathbb{K}[X]\}$. The {\em affine variety} $\mathcal{V}(I)$ is the set of points in $\mathbb{K}^n$ that are zeros of all the polynomials in $I$. If this is finite, then the ideal $I$ is {\em zero-dimensional}. It is {\em radical} if it contains all the polynomials $p\in \mathbb{K}[X]$ so that $p^m\in I$ for some natural $m$.

\vspace{0.1cm}

A {\em term order} on the set of monomials of $\mathbb{K}[X]$ is a multiplicative well-ordering whose smallest element is the constant monomial $1$. Thus, for instance, the {\em lexicographic} term order $<_{\mathrm{lex}}$ is defined so that, given two monomials $X^a=x_1^{a_1}\ldots x_n^{a_n}$ and $X^b=x_1^{b_1}\ldots x_n^{b_n}$, one has that $X^a<_{\mathrm{lex}}X^b$ if there exists a natural $m\leq n$ such that $a_i=b_i$ for all $i\leq m$ and $a_m<b_m$. The largest monomial of a polynomial with respect to a term order is its {\em leading monomial}. The {\em initial ideal} of an ideal $I\subseteq\mathbb{K}[X]$ is the ideal generated by the leading monomials of the non-zero polynomials of $I$. Any subset $G\subseteq I$ whose leading monomials generate this initial ideal is called a {\em Gr\"obner basis} of $I$ with respect to the underlying term order. Any monomial of $I$ that is not contained in its initial ideal is called {\em standard}. Regardless of the monomial term ordering, if the ideal $I$ is zero-dimensional and radical, then the number of standard monomials in $I$ coincides with the Krull dimension of the quotient ring $\mathbb{K}[X]/I$ and with the cardinality of $\mathcal{V}(I)$. This is obtained by means of the {\em Hilbert function}, which maps each non-negative integer $m$ onto $\mathrm{HF}_{\mathbb{K}[X]/I}(m)=\mathrm{dim}_{\mathbb{K}}(\mathbb{K}[X]_m/(\mathbb{K}[X]_m\cap I))$. Here, $\mathbb{K}[X]_m$ denotes the set of homogeneous polynomials in $\mathbb{K}[X]$ of degree $m$ and $\mathrm{HF}_{\mathbb{K}[X]/I}(m)$ coincides with the number of standard monomials in $I$ of degree $m$. The problem of computing Hilbert functions is NP-complete \cite{Bayer1992}. Its computation is based on that of a Gr\"obner basis of the ideal, whose complexity in case of dealing with a zero-dimensional ideal is $d^{O(n)}$ \cite{Lakshman1990a}, where $d$ is the maximal degree of the polynomials and $n$ is the number of variables.

\vspace{0.05cm}

The next result indicates how computational algebraic geometry can be used to enumerate and count the partial Latin rectangles in the set $\mathcal{R}_{r,s,n}$. Hereafter, the set of variables and the base field of the polynomial ring to be considered are, respectively, $X=\{x_{111},\ldots,x_{rsn}\}$ and the finite field $\mathbb{F}_2$.

\begin{theo}[\cite{Falcon2015}] \label{thm0} The set $\mathcal{R}_{r,s,n}$ is identified with the set of zeros of the zero-dimensional radical ideal in $\mathbb{F}_2[X]$
$$I_{r,s,n}=\langle\, x_{ijk}x_{i'jk},\, x_{ijk}x_{ij'k},\,x_{ijk}x_{ijk'}\colon\, i, i'\leq r; j,j'\leq s; k,k'\leq n\,\rangle.$$
Besides,
$$|\mathcal{R}_{r,s,n:m}|=\mathrm{HF}_{\mathbb{F}_2[X]/I_{r,s,n}}(m), \text{ for all } m\geq 0,$$
and
$$|\mathcal{R}_{r,s,n}|=\mathrm{dim}_{\mathbb{F}_2}(\mathbb{F}_2[X]/I_{r,s,n}).$$ \hfill $\Box$
\end{theo}

\vspace{0.3cm}

The proof of Theorem \ref{thm0} is based on the fact that every standard monomial $x_{111}^{a_{111}}\ldots x_{rsn}^{a_{rsn}}$ of the ideal $I_{r,s,n}$ can be identified with a partial Latin rectangle in $\mathcal{R}_{r,s,n}$ with set of entries $\{(i,j,k)\in [r]\times [s]\times [n] \colon\, a_{ijk}=1\}$. Particularly, the presence of the monomial $x_{ijk}x_{i'jk}$ as generator of the ideal $I_{r,s,n}$ involves the non-existence of the symbol $k$ twice in the $j^{th}$ column; that of $x_{ijk}x_{ij'k}$ involves the non-existence of the symbol $k$ twice in the $i^{th}$ row; and that of $x_{ijk}x_{ijk'}$ involves the non-existence of two distinct symbols in the cell $(i,j)$. Based on this result, the specialized algorithm described by Dickenstein and Tobis \cite{Dickenstein2012} was implemented in \cite{Falcon2015} for computing the cardinality of $\mathcal{R}_{r,s,n:m}$, for all $r, s, n \leq 4$. For higher orders, however, the required computational cost turned out to be excessive due to large memory storage requirements. This cost is only due to the computation of the corresponding Hilbert function, because the set of generators of $I_{r,s,n}$ constitutes itself a lexicographic Gr\"obner basis of the ideal. To reduce it, an alternative procedure is introduced in the next section. This is based on the similarity that exists among those generators in $I_{r,s,n}$ that correspond to distinct rows in a partial Latin rectangle. A preliminary version of this procedure was exposed in \cite{FalconStones2015}, where the cardinality of $\mathcal{R}_{r,s,n}$ was computed for all $r, s, n \leq 6$. For a better understanding of this procedure, the corresponding computation of $|\mathcal{R}_{3,3,3:2}|$ is illustrated in Example 1.

\section{An alternative procedure to compute $|\mathcal{R}_{r,s,n}|$}

For each positive integer $i\leq r$ we define the zero-dimensional subideal
$$I^{(i)}_{r,s,n}=\langle\,x_{ijk}x_{ij'k},\, x_{ijk}x_{ijk'}\colon\, j,j'\leq s; k,k'\leq n\,\rangle\subset I_{r,s,n}.$$

\vspace{0.1cm}

There exist distinct algorithms \cite{Hillebrand1999, Lazard1992, Moller1993} that enable us to decompose the zero-dimensional ideal $I^{(1)}_{r,s,n}$ into a finite set $\{J_{1,1},\ldots,J_{1,t}\}$ of subideals generated by triangular systems and whose affine varieties constitute a partition of $\mathcal{V}(I^{(1)}_{r,s,n})$. The complexity of this computation in the mentioned algorithms is polynomial once a lexicographic Gr\"obner basis of the ideal is known. This is our case, because the set of generators of $I^{(1)}_{r,s,n}$ constitutes itself one such a basis. Now, for each $i>1$ and $l\leq t$, let $J_{i,l}$ be the subideal of $I^{(i)}_{r,s,n}$ whose generators coincide with those of $J_{1,l}$ after replacing each variable $x_{1jk}$ by $x_{ijk}$. For each tuple $(t_1,\ldots,t_r)\in [t]^r$ we define the ideal
\begin{equation}\label{eqK}
K_{t_1,\ldots,t_r}=J_{1,t_1}+\ldots+ J_{r,t_r} + \langle\, x_{ijk}x_{i'jk}\colon\, i,i'\leq r; j\leq s; k\leq n\,\rangle.
\end{equation}

\vspace{0.1cm}

The triangularity of the underlying systems involves each subideal $J_{i,t_j}$ to have at least one generator of the form $x_{ij'k}$ or $x_{ij'k}-1$. The number of generators of the second form in the ideal $K_{t_1,\ldots,t_r}$ constitutes the minimum number of entries in a partial Latin rectangle that is identified with a point in $\mathcal{V}(K_{t_1,\ldots,t_r})$. We denote this number by $m_{t_1,\ldots,t_r}$.

\vspace{0.1cm}

\begin{prop}\label{prop1} Let $m$ be a non-negative integer. Then
$$\mathrm{HF}_{\mathbb{F}_2[X]/I_{r,s,n}}(m)=\sum_{\substack{(t_1,\ldots,t_r)\in [t]^r\\ m_{t_1,\ldots,t_r}\leq m}}\mathrm{HF}_{\mathbb{F}_2[X]/K_{t_1,\ldots,t_r}}(m-m_{t_1,\ldots,t_r}).$$
\end{prop}

{\bf Proof.} Let $X^a=x_{111}^{a_{111}}\ldots x_{rsn}^{a_{rsn}}$ be a standard monomial of degree $m$ in $I_{r,s,n}$. Since the ideals described in (\ref{eqK}) constitute a partition of the affine variety $\mathcal{V}(I_{r,s,n})$, there exists exactly one ideal $K_{t_1,\ldots,t_r}$ that contains the point $(a_{111},\ldots,a_{rsn})\in \mathcal{V}(I_{r,s,n})$. The result follows then from the fact that the monomial $X^a$ is uniquely related to the standard monomial $x_{111}^{a'_{111}}\ldots x_{rsn}^{a'_{rsn}}$ of degree $m-m_{t_1,\ldots,t_r}$ in $K_{t_1,\ldots,t_r}$, where $a'_{ijk}=0$ if $x_{ijk}-1$ is a generator of $K_{t_1,\ldots,t_r}$ and $a'_{ijk}=a_{ijk}$, otherwise. \hfill $\Box$

\vspace{0.5cm}

The smaller number of variables that are required to compute each addend in Proposition \ref{prop1}, together with the triangularity of the involved system and the possible parallel computation to determine distinct addends at the same time, reduce the running time and cost of computation of $\mathrm{HF}_{\mathbb{F}_2[X]/I_{r,s,n}}(m)$ in comparison with Theorem \ref{thm0}. Moreover, we do not need to compute all these addends, because $\mathrm{HF}_{\mathbb{F}_2[X]/K_{t_1,\ldots,t_r}}(m)=\mathrm{HF}_{\mathbb{F}_2[X]/K_{t_{\pi(1)}, \ldots,t_{\pi(r)}}}(m)$, for all $(t_1,\ldots,t_r)\in [t]^r$, $m\geq 0$ and $\pi\in S_r$.

\begin{exam}\label{example1} The ideal $I^{(1)}_{3,3,3}$ related to the first row of a partial Latin square of order $3$ can be decomposed into the next six disjoint subideals
\begin{enumerate}[i)]
\item $J_{1,1}=I^{(1)}_{3,3,3} + \langle\, x_{111},x_{121},x_{131}\,\rangle$.
\item $J_{1,2}=I^{(1)}_{3,3,3} + \langle\, x_{111},x_{121},x_{131}-1,x_{132},x_{133}\,\rangle$.
\item $J_{1,3}=I^{(1)}_{3,3,3} + \langle\, x_{111},x_{121}-1,x_{122}, x_{123},x_{131}\,\rangle$.
\item $J_{1,4}=I^{(1)}_{3,3,3} + \langle\, x_{111}-1,x_{112},x_{113},x_{121},x_{122},x_{131},x_{132}\,\rangle$.
\item $J_{1,5}=I^{(1)}_{3,3,3} + \langle\, x_{111}-1,x_{112},x_{113},x_{121},x_{122},x_{131},x_{132}-1,x_{133} \,\rangle$.
\item $J_{1,6}=I^{(1)}_{3,3,3} + \langle\, x_{111}-1,x_{112},x_{113},
x_{121},x_{122}-1,x_{123},x_{131},x_{132}\,\rangle$.
\end{enumerate}

\vspace{0.1cm}

Partial Latin squares of order $3$ are then distributed as points of
\begin{enumerate}
\item $\mathcal{V}(J_{1,1})$ if they do not contain the symbol $1$ in their first row.
\item $\mathcal{V}(J_{1,2})$ if they contain the symbol $1$ in the cell $(1,3)$.
\item $\mathcal{V}(J_{1,3})$ if they contain the symbol $1$ in the cell $(1,2)$.
\item $\mathcal{V}(J_{1,4})$ if they contain the symbol $1$ in the cell $(1,1)$ but do not contain the symbol $2$ in their first row.
\item $\mathcal{V}(J_{1,5})$ if they contain the symbol $1$ in the cell $(1,1)$ and the symbol $2$ in the cell $(1,3)$.
\item $\mathcal{V}(J_{1,6})$ if they contain the symbol $1$ in the cell $(1,1)$ and the symbol $2$ in the cell $(1,2)$.
\end{enumerate}

\vspace{0.1cm}

For each triple $(t_1,t_2,t_3)\in [6]^3$, we consider the ideal
$$K_{t_1,t_2,t_3}=J_{1,t_1}+J_{2,t_2}+ J_{3,t_3} + \langle\, x_{ijk}x_{i'jk}\colon\, i,i',j,k\leq 3\,\rangle.$$

The values of $\mathrm{HF}_{\mathbb{F}_2[X]/K_{t_1,t_2,t_3}}$ are exposed in Table \ref{tableHF}.

\begin{table}[h]
{
\caption{Hilbert functions related to the set $\mathcal{R}_{3,3,3}$.}
\label{tableHF}
\resizebox{\textwidth}{!}{
\begin{tabular}{rrrrrrrrrrrrrrrrr} \hline
\ & \multicolumn{13}{l}{$\mathrm{HF}_{\mathbb{F}_2[X]/K_{t_1,t_2,t_3}}(m)$}\\ \cline{2-17}
\ & \multicolumn{13}{l}{$t_1.t_2.t_3$}\\ \cline{2-17}
$m$ & 1.1.1 & 1.1.2 & 1.1.3 & 1.1.4 & 1.1.5 & 1.1.6 & 1.2.3 & 1.2.4 & 1.2.5 & 1.2.6 & 1.3.4 & 1.3.5 & 1.3.6 & 2.3.4 & 2.3.5 & 2.3.6\\ \hline
0 & 1 & 1 & 1 & 1 & 1 & 1 & 1 & 1 & 1 & 1 & 1 & 1 & 1 & 1 & 1 & 1\\
1 & 18 & 16 & 16 & 14 & 11 & 11 & 14 & 12 & 10 & 9 & 12 & 9 & 10 & 10 & 8 & 8\\
2 & 108 & 84 & 84 & 62 & 36 & 36 & 64 & 45 & 29 & 24 & 45 & 24 & 29 & 32 & 19 & 19\\
3 & 264 & 176 & 176 & 104 & 42 & 42 & 116 & 63 & 29 & 23 & 63 & 23 & 29 & 38 & 16 & 16\\
4 & 270 & 150 & 150 & 66 & 18 & 18 & 84 & 32 & 11 & 8 & 32 & 8 & 11 & 16 & 5 & 5 \\
5 & 108 & 48 & 48 & 12 & 2 & 2 & 24 & 5 & 1 & 1 & 5 & 1 & 1 & 2 & 1 & 1\\
6 & 12 & 4 & 4 & 0 & 0 & 0 & 2 & 0 & 0 & 0 & 0 & 0 & 0 & 0 & 0 & 0\\ \hline
\end{tabular}
}}
\end{table}

Let $m_{t_1,t_2,t_3}$ be the number of generators of the form $x_{ijk}-1$ in the ideal $K_{t_1,t_2,t_3}$. Thus, for instance, every point of the affine variety $\mathcal{V}(K_{6,3,2})$ is uniquely related to a partial Latin square of order $3$ and size at least $m_{6,3,2}=4$. This last value holds from the fact that the set of entries of any such a partial Latin square always contains the subset $\{(1,1,1),(1,2,2),(2,2,1),(3,3,1)\}$. From Proposition \ref{prop1}, we have, for example, that

{\small
$$\begin{array}{c}|\mathcal{R}_{3,3,3:2}|=
\mathrm{HF}_{\mathbb{F}_2[X]/K_{1,1,1}}(2)+ 3\ \mathrm{HF}_{\mathbb{F}_2[X]/K_{1,1,2}}(1)+ 3\ \mathrm{HF}_{\mathbb{F}_2[X]/K_{1,1,3}}(1)+\\
3\ \mathrm{HF}_{\mathbb{F}_2[X]/K_{1,1,4}}(1)+
3\ \mathrm{HF}_{\mathbb{F}_2[X]/K_{1,1,5}}(0)+
3\ \mathrm{HF}_{\mathbb{F}_2[X]/K_{1,1,6}}(0)+\\
6\ \mathrm{HF}_{\mathbb{F}_2[X]/K_{1,2,3}}(0)+
6\ \mathrm{HF}_{\mathbb{F}_2[X]/K_{1,2,4}}(0)+
6\ \mathrm{HF}_{\mathbb{F}_2[X]/K_{1,3,4}}(0)=
270.\end{array}$$}
\hfill $\lhd$
\end{exam}

\vspace{0.5cm}

This computational algebraic method has been implemented in the procedure {\em PLR} of the library {\em pls.lib}, available online on \texttt{http://personales.us.es/} \texttt{raufalgan/LS/pls.lib}, for the open computer algebra system for polynomial computations {\sc Singular} \cite{Decker2016}. The correctness and termination of this procedure are based on those of the algorithms described in \cite{Dickenstein2012, Hillebrand1999, Moller1993} for computing Hilbert functions. In order to prove its efficiency, we have firstly checked the known cardinality of $\mathcal{R}_{r,s,n:m}$, for all $r,s,n\leq 4$ (see Table \ref{table1}), which was already computed in \cite{Falcon2015}. In the same computer system, an {\em Intel Core i7-2600  CPU (8 cores), with a 3.4 GHz processor and 16 GB of RAM}, the maximum running time decreases from 50 seconds in \cite{Falcon2015} to $0$ seconds. This corresponds to the computation of the series $|\mathcal{R}_{4,4,4;m}|$. The procedure has then been applied for computing in Tables \ref{table2}--\ref{table4} the rest of cases so that $r\leq s\leq n\leq 6$. The running time ranges here from 0 seconds to 32 hours. This maximum running time corresponds to the computation of the series $|\mathcal{R}_{6,6,6;m}|$, for which 2,3 GB of RAM is required. For higher orders, the first series whose computation turned out to be excessive for our computer system due to large memory storage requirements was $|\mathcal{R}_{6,7,7;m}|$. In order to improve the efficiency of this computational algebraic method, we propose in the next section to impose some extra algebraic conditions to our base ideal. They are referred to the distribution of non-empty cells per row and column in a partial Latin rectangle and to the number of occurrences of each symbol.

\begin{table}[h]
{
\caption{Distribution of $\mathcal{R}_{r,s,n}$ according to the size, for $r\leq s\leq n\leq 4$.}
\label{table1}
\resizebox{\textwidth}{!}{
\begin{tabular}{rrrrrrrrrrrrrrrrrrrrrrr}
 \hline
\ & \multicolumn{5}{l}{$|\mathcal{R}_{r,s,n:m}|$}\\ \cline{2-21}
\ & \multicolumn{5}{l}{$r.s.n$}\\ \cline{2-21}
$m$ & 1.1.1 & 1.1.2 & 1.1.3 & 1.1.4 & 1.2.2 & 1.2.3 & 1.2.4 & 1.3.3 & 1.3.4 & 1.4.4 & 2.2.2 & 2.2.3 & 2.2.4 & 2.3.3 & 2.3.4 & 2.4.4 & 3.3.3 & 3.3.4 & 3.4.4 & 4.4.4\\ \hline
0 & 1 & 1 & 1 & 1 & 1 & 1 & 1 & 1 & 1 & 1 & 1 & 1 & 1 & 1 & 1 & 1 & 1 & 1 & 1 & 1\\
1 & 1 & 2 & 3 & 4 & 4 & 6  & 8  & 9 & 12 & 16 & 8 & 12 & 16 & 18 & 24 & 32 & 27 & 36 & 48 & 64\\
2 & &  &  &   & 2 & 6  & 12  & 18 & 36 & 72 & 16 & 42 & 80 & 108 & 204 & 384 & 270 & 504 & 936 & 1,728\\
3 & &  &  &   &  &   &   &  6 & 24 & 96 & 8 & 48 & 144 & 264 & 768 & 2,208 & 1,278 & 3,552 & 9,696 & 25,920\\
4 & &  &  &   &  &   &   &  & & 24 & 2 & 18 & 84 & 270 & 1,332 & 6,504 & 3,078 & 13,716 & 58,752 & 239,760\\
5 & &  &  &   &  &   &   &  & &   &   &   &   &  108 & 1,008 & 9,792 & 3,834 & 29,808 & 216,864 & 1,437,696\\
6 & &  &  &   &  &   &   &  & &   &   &   &   &  12 & 264 & 7,104 & 2,412 & 36,216 & 494,064 & 5,728,896\\
7 & &  &  &   &  &   &   &  & &   &   &   &   &    &   &  2,112 & 756 & 23,760 & 691,200 & 15,326,208\\
8 & &  &  &   &  &   &   &  & &   &   &   &   &    &   &  216 & 108 & 7,776 & 581,688 & 27,534,816\\
9 & &  &  &   &  &   &   &  & &   &   &   &   &    &   &    & 12 & 1,056 & 283,584 & 32,971,008\\
10 & &  &  &   &  &   &   &  & &   &   &   &   &    &   &    &   &   &  75,744 & 25,941,504\\
11 & &  &  &   &  &   &   &  & &   &   &   &   &    &   &    &   &   & 10,368 & 13,153,536\\
12 & &  &  &   &  &   &   &  & &   &   &   &   &    &   &    &   &   & 576 & 4,215,744\\
13 & &  &  &   &  &   &   &  & &   &   &   &   &    &   &    &   &   &   & 847,872\\
14 & &  &  &   &  &   &   &  & &   &   &   &   &    &   &    &   &   &   & 110,592\\
15 & &  &  &   &  &   &   &  & &   &   &   &   &    &   &    &   &   &   & 9,216\\
16 & &  &  &   &  &   &   &  & &   &   &   &   &    &   &    &   &   &   & 576\\ \hline
Total & 2 & 3 & 4
& 5 & 7 & 13 & 21 & 34 & 73 & 209 & 35 & 121 & 325 & 781 & 3,601 & 28,353 & 11,776 & 116,425 & 2,423,521 & 127,545,137\\ \hline
\end{tabular}
}}
\end{table}

\begin{table}[h]
{
\caption{Distribution of $\mathcal{R}_{r,s,5}$ according to the size, for $r\leq s\leq 5$.}
\label{table2}
\resizebox{\textwidth}{!}{
\begin{tabular}{rrrrrrrrrrrrrrrr} \hline
\ & \multicolumn{5}{l}{$|\mathcal{R}_{r,s,5:m}|$}\\ \cline{2-16}
\ & \multicolumn{5}{l}{$r.s.5$}\\ \cline{2-16}
$m$ & 1.1.5 & 1.2.5 & 1.3.5 & 1.4.5 & 1.5.5 & 2.2.5 & 2.3.5 & 2.4.5 & 2.5.5 & 3.3.5 & 3.4.5 & 3.5.5 & 4.4.5 & 4.5.5 & 5.5.5\\ \hline
0 & 1 & 1 & 1 & 1 & 1 & 1 & 1 & 1 & 1 & 1 & 1 & 1 & 1 & 1 & 1\\
1 & 5 & 10 & 15 & 20 & 25 & 20 & 30 &  40 & 50 & 45 & 60 & 75 & 80 & 100 & 125 \\
2 &  & 20 & 60 & 120 & 200 & 130 & 330 & 620 & 1,000 & 810 & 1,500 & 2,400 & 2,760 & 4,400 & 7,000\\
3&  &  &  60 & 240 & 600 & 320 & 1,680 & 4,800 & 10,400 & 7,590 & 20,520 & 43,200 & 54,240 & 112,800&233,000
\\
4 &  &  &  &  120 & 600 & 260 & 4,140 & 20,040 & 61,400 & 40,500 & 169,920 & 486,000 & 676,200 & 1,881,600&5,159,000\\
5 &  &  &  &  &  120 &  & 4,680 & 45,600 & 211,440 & 126,900 & 891,360 & 3,594,960 & 5,641,920 & 21,612,480 & 80,602,200\\
6 &  &  &  &  &  &  &  1,920 & 54,480 &  421,200 & 232,680 & 3,018,000 & 17,930,400 & 32,423,520 & 176,546,400 & 920,160,000\\
7 &  &  &  &  &  &  &  & 30,720 & 465,600 & 240,840 & 6,605,280 & 60,912,000 & 130,248,960 & 1,045,147,200& 7,845,192,000\\
8 &  &  &  &  &  &  &  & 6,360 & 262,200 & 128,520 & 9,224,280 & 140,826,600 & 367,731,360 & 4,530,640,800 & 50,648,616,000\\
9 &  &  &  &  &  &  &  &  & 63,600 & 27,480 & 7,983,840 &  219,307,800 &  728,440,320 & 14,444,083,200& 249,687,408,000\\
10&  &  &  &  &  &  &  &  & 5,280 & & 4,063,680 & 225,419,040 & 1,004,380,800 & 33,852,910,080 &  944,069,668,800\\
11 &  &  &  &  &  &  &  &  &  & & 1,100,160 & 148,010,400 & 950,238,720 & 58,065,734,400 &  2,741,210,616,000\\
12 &  &  &  &  &  &  &  &  &  & & 120,960 & 59,047,200 & 603,722,880 & 72,278,294,400 &  6,104,066,712,000\\
13 &  &  &  &  &  &  &  &  &  &  &  & 13,284,000 & 249,580,800 & 64,484,985,600 & 10,385,299,320,000\\
14 &  &  &  &  &  &  &  &  &  &  &  & 1,512,000 & 63,884,160 & 40,544,726,400 & 13,420,351,008,000\\
15 &  &  &  &  &  &  &  &  &  &  &  & 66,240 & 9,216,000 & 17,571,260,160 & 13,065,814,483,200\\
16 &  &  &  &  &  &  &  &  &  &  &  &  & 590,400 & 5,099,169,600 & 9,486,099,648,000 \\
17 &  &  &  &  &  &  &  &  &  &  &  &  & &  953,107,200 & 5,073,056,640,000 \\
18 &  &  &  &  &  &  &  &  &  &  &  &  & &  108,288,000 & 1,970,474,400,000 \\
19 &  &  &  &  &  &  &  &  &  &  &  &  & & 6,681,600 & 547,608,096,000\\
20 &  &  &  &  &  &  &  &  &  &  &  &  & & 161,280 & 107,330,054,400\\
21 &  &  &  &  &  &  &  &  &  &  &  &  & &  & 14,667,552,000\\
22 &  &  &  &  &  &  &  &  &  &  &  &  & &  & 1,388,160,000\\
23 &  &  &  &  &  &  &  &  &  &  &  &  & &  & 91,008,000\\
24 &  &  &  &  &  &  &  &  &  &  &  &  & &  & 4,032,000\\
25 &  &  &  &  &  &  &  &  &  &  &  &  & &  & 161,280\\ \hline
Total & 6 & 31 & 136 & 501 & 1,546 & 731 & 12,781 & 162,661 & 1,502,171 & 805,366 & 33,199,561 & 890,442,316 & 4,146,833,121 & 313,185,347,701 & 64,170,718,937,006\\ \hline
\end{tabular}
}}
\end{table}

\begin{table}[h]
{
\caption{Distribution of $\mathcal{R}_{r,s,6}$ according to the size, for $r\leq s\leq 6$ (I).}
\label{table3}
\resizebox{\textwidth}{!}{
\begin{tabular}{rrrrrrrrrrrrrrrr} \hline
\ & \multicolumn{5}{l}{$|\mathcal{R}_{r,s,6:m}|$}\\ \cline{2-16}
\ & \multicolumn{5}{l}{$r.s.6$}\\ \cline{2-16}
$m$ & 1.1.6 & 1.2.6 & 1.3.6 & 1.4.6 & 1.5.6  & 1.6.6 & 2.2.6 & 2.3.6 & 2.4.6 & 2.5.6 & 2.6.6 & 3.3.6 & 3.4.6 & 3.5.6 & 3.6.6\\ \hline
0 & 1 & 1 & 1 & 1 & 1 & 1 & 1 & 1 & 1 & 1 & 1 & 1 & 1 & 1 & 1 \\
1 & 6 & 12 & 18 & 24 & 30 & 36 & 24 & 36 & 48 & 60 & 72  & 54 & 72 & 90 &  108\\
2 &  & 30 & 90 & 180 & 300 & 450 & 192 & 486  & 912 & 1,470 & 2,160 & 1,188 & 2,196 & 3,510 &  5,130\\
3 &  &  & 120 & 480 & 1,200 & 2,400 & 600 & 3,120  & 8,880 & 19,200 & 35,400 & 13,896 & 37,344 & 78,360 &  141,840\\
4 &  &  &  & 360 & 1,800 & 5,400 & 630 & 9,990  & 48,060 & 146,700 & 349,650 & 94,770 & 392,580 & 1,115,100 &  2,547,450 \\
5 &  &  &  &  & 720 & 4,320 &  & 15,120 & 146,880  & 678,240 & 2,168,640 & 389,340 & 2,676,240  & 10,667,160 & 31,419,360\\
6 &  &  &  &  & & 720 & & 8,520 & 245,760 & 1,899,600 & 8,546,880 & 961,380 & 12,082,680 & 70,540,800 & 274,470,480\\
7 &  &  &  &  &  &  &  &   & 204,480 & 3,139,200 & 21,211,200 & 1,375,920 & 36,270,720 & 326,808,000 & 1,727,352,000\\
8 &  &  &  &  &  &  &  &   & 65,160 & 2,881,800 & 32,189,400 & 1,038,960 & 71,633,160 & 1,064,140,200 & 7,893,282,600 \\
9 &  &  &  &  &  &  &  &   &  & 1,303,200 & 28,267,200 & 317,760 & 90,585,600 & 2,422,568,400 & 26,212,965,600 \\
10&  &  &  &  &  &  &  &   &  & 222,480 & 13,063,680 &  & 69,603,840 & 3,803,369,040 &
62,938,898,640\\
11 &  &  &  &  &  &  &  &   &  &  & 2,669,760 & & 29,255,040 & 4,021,099,200 &
108,045,861,120\\
12 &  &  &  &  &  &  &  &   &  &  & 190,800 & & 5,112,000 & 2,756,361,600 &
130,246,779,600\\
13 &  &  &  &  &  &  &  &   &  &  &  &  &  &  1,152,144,000 &
107,367,120,000\\
14 &  &  &  &  &  &  &  &   &  &  &  &  &  & 262,828,800 &
58,252,478,400 \\
15 &  &  &  &  &  &  &  &   &  &  &  &  &  &  24,791,040 &
19,683,613,440\\
16 &  &  &  &  &  &  &  &   &  &  &  &  &  &  &
3,828,798,720\\
17 &  &  &  &  &  &  &  &   &  &  &  &  &  &  & 384,652,800\\
18 &  &  &  &  &  &  &  &   &  &  &  &  &  &  &  15,321,600\\
\hline
Total & 7 & 43 & 229 & 1,045 & 4,051 & 13,327 & 1,447 &  37,273 & 720,181 & 10,291,951 & 108,694,843 & 4,193,269 & 317,651,473 & 15,916,515,301 & 526,905,708,889 \\ \hline
\end{tabular}
}}
\end{table}

\begin{table}[h]
{\scriptsize
\caption{Distribution of $\mathcal{R}_{r,s,6}$ according to the size, for $r\leq s\leq 6$ (II).}
\label{table4}
\resizebox{\textwidth}{!}{
\begin{tabular}{rrrrrrrr} \hline
\ & \multicolumn{5}{l}{$|\mathcal{R}_{r,s,6:m}|$}\\ \cline{2-7}
\ & \multicolumn{5}{l}{$r.s.6$}\\ \cline{2-7}
$m$ & 4.4.6 & 4.5.6 & 4.6.6 & 5.5.6 & 5.6.6 & 6.6.6\\ \hline
0 & 1 & 1 & 1 & 1 & 1 & 1\\
1 & 96 & 120 &144 & 150 & 180 & 216 \\
2 & 4,032 & 6,420  &9,360 & 10,200 & 14,850 & 21,600\\
3 & 98,016 & 203,040  &364,560 & 417,600 & 746,400 & 1,330,200\\
4 & 15,384,24 & 4,245,120 &9,527,220 & 11,532,600 & 25,631,100 & 56,614,950\\
5 & 16,476,480 & 62,189,280 &177,310,080  & 228,154,320 & 639,260,640 & 1,771,796,160\\
6 & 124,148,160 & 660,375,600 &2,434,907,520  & 3,352,566,000 &  12,019,602,000 & 423,57,620,160\\
7 & 669,176,640 & 5,189,068,800 &25,231,996,800 & 37,450,656,000 &  174,585,456,000 & 793,416,600,000 \\
8 & 2,599,625,880 & 30,548,079,000 & 200,165,742,000 & 322,946,451,000 & 1,991,858,418,000 & 11,852,197,317,000\\
9 & 7,281,623,040 & 135,625,603,200 &1,226,542,944,000 & 2,171,483,394,000 &  18,056,836,776,000  & 142,993,809,528,000\\
10& 14,618,868,480 &
455,097,055,680 &5,834,154,055,680 & 11,456,637,616,800 & 131,095,655,863,200  & 1,406,144,941,776,000\\
11 & 20,771,527,680 &  1,152,338,169,600 &21,579,415,960,320 & 47,586,889,008,000 & 766,225,199,808,000  & 11,344,829,123,448,000 \\
12 & 20,451,767,040 & 2,190,542,918,400 &62,007,749,812,800 & 155,763,852,264,000 & 3,616,441,279,056,000  & 75,444,662,621,250,000 \\
13 & 13,491,532,800 & 3,099,028,723,200 &137,935,650,124,800 & 401,342,211,504,000 & 13,801,803,749,280,000 & 414,809,990,051,328,000\\
14 & 5,635,215,360 & 3,221,159,616,000 &236,112,048,230,400  & 811,559,781,792,000 & 42,582,496,312,944,000 & 1,888,965,825,155,136,000\\
15 & 1,337,610,240 &  2,415,807,221,760 &308,313,104,578,560 & 1,281,622,863,052,800 & 106,042,151,250,892,000 & 7,129,083,890,074,291,200\\
16 & 137,116,800 &  1,274,532,969,600 &303,524,671,011,840 & 1,569,898,647,504,000 & 212,529,994,957,440,000 & 22,290,972,757,613,899,200\\
17 & &  455,792,486,400 &221,831,824,435,200 & 1,478,352,018,528,000 & 341,378,166,715,776,000 & 57,672,207,579,205,440,000\\
18 & & 104,134,464,000 &117,967,540,608,000 & 1,058,153,580,288,000 & 437,045,603,416,704,000 & 123,205,370,805,154,944,000\\
19 & & 13,604,889,600 &44,468,899,430,400  & 567,490,862,592,000 &   442,874,461,303,296,000 & 216,689,524,093,737,792,000\\
20 & & 767,854,080 &11,483,903,278,080 & 223,899,017,011,200 & 352,217,521,389,081,000 & 312,570,613,181,156,803,200 \\
21 & &  &1,942,917,304,320 & 63,429,754,752,000 & 217,606,324,462,848,000  & 368,084,100,503,749,939,200\\
22 & &  &202,499,481,600 & 12,467,229,696,000  & 103,166,400,104,064,000  & 351,915,364,298,700,288,000\\
23 & &  &11,670,220,800 & 1,610,606,592,000 & 36,987,139,952,640,000 & 271,409,503,369,430,016,000\\
24 & &  &283,046,400 & 123,628,032,000 & 9,853,601,458,752,000 & 167,607,699,757,168,896,000 \\
25 &  &  &  & 4,356,218,880 & 1,909,729,461,012,480 & 82,187,524,303,374,458,880\\
26 &  &  &  &  &  262,267,391,462,400 & 31,703,766,748,202,926,080\\
27 &  &  &  &  &  24,634,533,888,000 & 9,523,824,649,261,056,000\\
28 &  &  &  &  &  1,496,724,480,000 & 2,204,514,949,427,712,000\\
29 &  &  &  &  &  52,752,384,000 & 389,140,940,150,784,000\\
30 &  &  &  &  &  812,851,200 & 51,905,194,846,617,600\\
31 &  &  &  &  &  & 5,196,712,196,505,600\\
32 &  &  &  &  &  & 389,383,137,792,000 \\
33 &  &  &  &  &  & 21,862,379,520,000\\
34 &  &  &  &  &  & 925,655,040,000\\
35 &  &  &  &  &  & 29,262,643,200\\
36 &  &  &  &  &  & 812,851,200\\ \hline
Total & 87,136,329,169 & 14,554,896,138,901 &
1,474,670,894,380,885 &
7,687,297,409,633,551
 & 2,322,817,844,850,427,451 & 202,7032,853,070,203,981,647\\ \hline
\end{tabular}
}}
\end{table}

\section{Shape, type and structure of partial Latin rectangles}\label{sec:shapes}

The {\em shape} of a partial Latin rectangle $P=(p_{ij})\in\mathcal{R}_{r,s,n}$ is defined as the $r\times s$ binary array $B_P=(b_{ij})$ such that $b_{ij}=1$ if $(i,j,p_{ij})\in E(P)$ and $0$, otherwise. Let $\mathrm{r}_i$, $\mathrm{c}_j$ and $\mathrm{s}_k$ respectively be the number of filled cells in the $i^{th}$ row and $j^{th}$ column of $P$ and the number of occurrences of the symbol $k$ in $P$. According to the terminology exposed by Keedwell \cite{Keedwell1994} and generalized by Bean et al. \cite{Bean2002}, the tuples $R=(\mathrm{r}_1,\ldots,\mathrm{r}_r)$, $C=(\mathrm{c}_1,\ldots,\mathrm{c}_s)$ and $S=(\mathrm{s}_1,\ldots,\mathrm{s}_n)$ determine, respectively, the {\em row}, {\em column} and {\em symbol} types of $P$. The {\em type} of $P$ is then defined as the triple $(R,C,S)$. Thus, for instance, the type of the partial Latin square of Figure \ref{Fig2} is $((2,2,1,0),(2,1,1,1),(2,3,0,0))$. Hereafter, the set of partial Latin rectangles of type $(R,C,S)$ is denoted by $\mathcal{R}_{R,C,S}$.

\vspace{0.1cm}

Let $\mathcal{T}_{n,m}$ be the set of $n$-tuples $T=(t_1,\ldots,t_n)$ of {\em weight} $\sum_{i\leq n}t_i=m$ whose components are non-negative integers. The {\em conjugate} of $T$ is the tuple $T^*=(\mathrm{t}_1^*,\ldots,\mathrm{t}_m^*)$, where each $\mathrm{t}_i^*$ is the number of positive integers $j\leq n$ such that $\mathrm{t}_j\geq i$. If $\overline{T}= (\overline{\mathrm{t}}_1,\ldots,\overline{\mathrm{t}}_n)\in \mathcal{T}_{n,m}$ is obtained after a decreasing rearrangement of the components of $T$, then $T$ is said to be {\em majorized} by a second tuple $T'=(\mathrm{t}'_1,\ldots,\mathrm{t}'_n)\in \mathcal{T}_{n,m}$ if $\sum_{i\leq j} \overline{\mathrm{t}}_i \leq \sum_{i\leq j} \overline{\mathrm{t}'_i}$, for all $j\leq n$. This gives rise to the so-called {\em dominance order} $\preceq$ on $\mathcal{T}_{n,m}$ \cite{Brylawski1973}.

\begin{theo}\label{thmRCS} Let $(R, C, S)\in \mathcal{T}_{r,m}\times \mathcal{T}_{s,m}\times \mathcal{T}_{n,m}$. The set $\mathcal{R}_{R,C,S}$ is non-empty only if $C\preceq R^*$, $S\preceq C^*$ and $R\preceq S^*$.
\end{theo}

{\bf Proof.} The set of shapes of partial Latin rectangles of row type $R$ and column type $C$ is identified with the set of $r\times s$ binary matrices whose row and column sum vectors coincide, respectively, with $R$ and $C$. According to the Gale-Ryser theorem \cite{Ford1962, Gale1957, Ryser1957}, this set is non-empty if and only if $C\preceq R^*$. This constitutes, therefore, a necessary condition for the set $\mathcal{R}_{R,C,S}$ to be non-empty. The result holds then from parastrophism. \hfill $\Box$

\vspace{0.5cm}

The previous result gives a necessary condition to deal with the problem of deciding whether a triple $(R,C,S)\in\mathcal{T}_{r,m}\times\mathcal{T}_{s,m}\times\mathcal{T}_{n,m}$ is the type of a partial Latin rectangle in $\mathcal{R}_{r,s,n:m}$. Nevertheless, this condition is not sufficient because, for instance, $\mathcal{R}_{(3,1,1),(3,1,1),(3,1,1)}=\emptyset$, but $(3,1,1)^*=(3,1,1)$. This problem is equivalent to that of deciding whether a tripartite graph with a given degree sequence has an edge-partition into triangles \cite{Colbourn1984}. Specifically, any partial Latin rectangle $P\in\mathcal{R}_{R,C,S}$ is identified with an edge-partition into triangles of a labeled tripartite graph $(V_1\cup V_2\cup V_3, E_1\cup E_2\cup E_3)$ such that
\begin{enumerate}[a)]
\item $|V_1|=r$, $|V_2|=s$ and $|V_3|=n$.
\item The vertices of $V_1$, $V_2$ and $V_3$ are uniquely and respectively related to the rows, columns and symbols of $P$.
\item The bi-adjacency matrices of the three bipartite graphs $(V_1\cup V_2,E_1)$, $(V_1\cup V_3,E_2)$ and $(V_2\cup V_3,E_3)$ are, respectively, the binary matrices related to the shape of $P$ and that of its two parastrophic partial Latin rectangles $P^{(23)}$ and $P^{(132)}$.
\end{enumerate}

\begin{figure}
\begin{center}
{
\begin{tabular}{ccc}
$\begin{array}{l}
B_P\equiv\left(\begin{array}{cccc}
  1 & 1 & 1 & 0 \\
  0 & 1 & 0 & 1
\end{array}\right)\\ \ \\
B_{P^{(23)}}\equiv\left(\begin{array}{ccc}
  1 & 1 & 1 \\
  1 & 0 & 1
\end{array}\right)\\ \ \\
B_{P^{(132)}}\equiv\left(\begin{array}{ccc}
  1 & 0 & 0 \\
  1 & 1 & 0 \\
  0 & 0 & 1 \\
  0 & 0 & 1
\end{array}\right)
\end{array}$
 & \begin{tabular}{c}  \ \\ \includegraphics[width=0.3\textwidth]{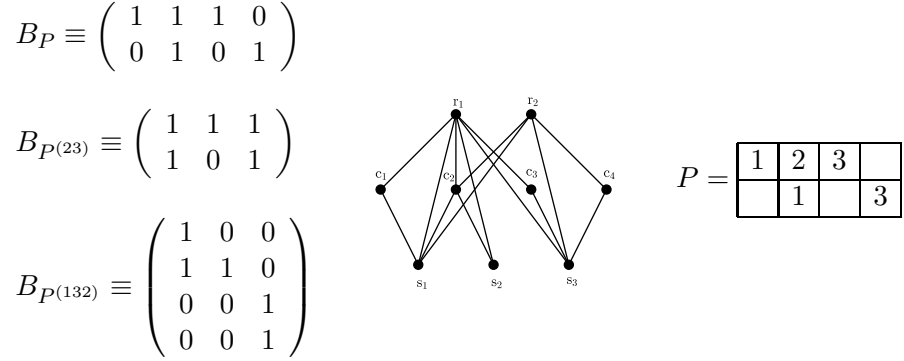}
 \end{tabular} & $P=\begin{array}{|c|c|c|c|} \hline
  1 & 2 & 3 & \\ \hline
   & 1 &  & 3\\ \hline
\end{array}$
\end{tabular}}
\end{center}
\caption{Shapes, tripartite graph and partial Latin rectangle in $\mathcal{R}_{2,4,3}$ related to the type $((3,2), (1,2,1,1), (2,1,2))$.}
\label{fig_rcs}
\end{figure}

This graph satisfies the necessary condition of being {\em uniform} in order to have an edge-partition into triangles. That is, the number of $V_1$-to-$V_2$ edges is equal to that of $V_1$-to-$V_3$ edges and also to that of $V_2$-to-$V_3$ edges. This number coincides with the component of the tuple $R$ (respectively, $C$ and $S$) that is related to that vertex. The partial Latin rectangle $P$ is then uniquely identified with that edge-partition into triangles in which the symbol included in an entry of $P$ is determined by the symbol vertex of the triangle that contains the row and column vertices associated to that cell (see Figure \ref{fig_rcs}).

\vspace{0.5cm}

Computational algebraic geometry can be used to determine explicitly the set $\mathcal{R}_{R,C,S}$. In this regard, the next result indicates those polynomials that have to be added to the set of generators of the ideal $I_{r,s,n}$ in Theorem \ref{thm0} in order to determine the set $\mathcal{R}_{R,C,S}$. Since the constant terms of these new polynomials coincide with the components of the tuples $R$, $C$ and $S$, the order of the base field $\mathbb{F}_2$ in the mentioned theorem is conveniently replaced here by a prime $p\geq 2$. Theorem \ref{thm0} is also valid for this new base field $\mathbb{F}_p$.

\begin{theo} \label{thm_RCS2} Let $R=(\mathrm{r}_1,\ldots,\mathrm{r}_r)$, $C=(\mathrm{c}_1,\ldots,\mathrm{c}_s)$ and $S=(\mathrm{s}_1,\ldots,\mathrm{s}_n)$ be three tuples in $\mathcal{T}_{r,m}$, $\mathcal{T}_{s,m}$ and $\mathcal{T}_{n,m}$, respectively, and let $p$ be the first prime greater than the maximum of all the components of $R$, $C$ and $S$. The set $\mathcal{R}_{R,C,S}$ is identified with the set of zeros of the zero-dimensional radical ideal
$$I_{R,C,S}=I_{r,s,n} + \langle\,\mathrm{r}_i-\sum_{j\leq s,k\leq n} x_{ijk}\colon\, i\leq r\,\rangle +\langle\,\mathrm{c}_j-\sum_{i\leq r,\, k\leq n} x_{ijk}\colon\, j\leq s\,\rangle +$$ $$ \langle\,\mathrm{s}_k-\sum_{i\leq r,j\leq s} x_{ijk}\colon\, k\leq n\,\rangle\subset \mathbb{F}_p[X].$$
Besides, $|\mathcal{R}_{R,C,S}|= \mathrm{dim}_{\mathbb{F}_p}(\mathbb{F}_p[X]/I_{R,C,S})$.
\end{theo}

{\bf Proof.} Since $I_{R,C,S}\subset I_{r,s,n}$, each zero of the ideal $I_{R,C,S}$ is uniquely related to a partial Latin rectangle in $\mathcal{R}_{r,s,n}$. The three subideals that are added to $I_{r,s,n}$ in the definition of $I_{R,C,S}$ involve these partial Latin rectangles to be exactly those ones having $R$, $C$ and $S$ as their row, column and symbol types, respectively. Now, in order to prove the last assertion, observe that the finiteness of $\mathcal{R}_{r,s,n}$ involves $I_{R,C,S}$ to be zero-dimensional and that the intersection between this ideal and the polynomial ring $\mathbb{F}_p[x_{ijk}]$ coincides with the ideal generated by the polynomial $x_{ijk}\left(x_{ijk}-1\right)$, for all $(i,j,k)\in [r]\times [s]\times [n]$. This is contained in $I_{R,C,S}$, which is, therefore, not only zero-dimensional, but also radical. Hence, its number of zeros coincides with $\mathrm{dim}_{\mathbb{F}_p}(\mathbb{F}_p[X]/I_{R,C,S})$. \hfill $\Box$

\vspace{0.5cm}

The {\em structure} of an $n$-tuple $T=(\mathrm{t}_1,\ldots,\mathrm{t}_n)\in \mathcal{T}_{n,m}$ is defined as the expression $z_T=m^{d_m}\ldots 1^{d_1}$, where $d_i$ is the number of occurrences of a given non-negative integer $i$ as a component of $T$. In practice, only those terms $i^{d_i}$ for which $d_i>0$ are written. The {\em length} of the structure $z_T$ is $\sum_{i\leq m}d_i$ and its {\em weight} is $\sum_{i\leq m}i d_i=m$. Hereafter, the set of structures of length $l$ and weight $m$ is denoted by $\mathcal{Z}_{l,m}$. Thus, for instance, the structure of the tuple $(3,1,3,3,1,0)$ is $3^31^2\in \mathcal{Z}_{5,11}$. Isotopisms of partial Latin rectangles preserve the structures of the row, column and symbol types of a partial Latin rectangle. This becomes essential for their enumeration and classification because of the following result.

\begin{lem} \label{lemma_type_1} The number of partial Latin rectangles of a given row, column or symbol type only depends on its structure.
\end{lem}

{\bf Proof.}  Let $T=(\mathrm{t}_1,\ldots,\mathrm{t}_n)\in\mathcal{T}_{n,m}$ and $T'=(\mathrm{t}'_1,\ldots,\mathrm{t}'_{n'})\in\mathcal{T}_{n',m}$ be two tuples with the same structure $z_T=z_{T'}$. Suppose $n\leq n'$. Then, there exists a permutation $\pi$ on $[n]$ such that $\mathrm{t}_i=\mathrm{t}'_{\pi(i)}$ for all $i\leq n$. The rest of components of $T'$ are zeros and do not have any influence on the number of partial Latin rectangles having $T'$ as row, column or symbol type. The same permutation $\pi$ enable us to identify the rows, columns or symbols of two partial Latin rectangles having $T$ and $T'$ as row, column or symbol types, respectively. \hfill $\Box$

\vspace{0.5cm}

Let $P$ be a partial Latin rectangle of type $(R,C,S)\in\mathcal{T}_{r,m}\times \mathcal{T}_{s,m}\times \mathcal{T}_{n,m}$. Its {\em structure} is defined as the triple $(z_R,z_C,z_S)$, where $z_R$, $z_C$ and $z_S$ are called, respectively, the {\em row}, {\em column} and {\em symbol structures} of $P$. Thus, for instance, the partial Latin square of Figure \ref{Fig2} has structure $(2^21,21^3,32)\in \mathcal{Z}_{3,5}\times \mathcal{Z}_{4,5}\times \mathcal{Z}_{2,5}$. Some structures of partial Latin squares have been widely studied in the literature:
\begin{enumerate}[a)]
\item If the empty cells of a partial Latin square of order $n$ are replaced by zeros, then the structure $(k^n,k^n,n^k)$ is related to the set of {\em $F(n;n-k,1^k)$-squares} \cite{Hedayat1970}.
\item The structure $(k^n,k^n,k^n)$ is that of a {\em $k$-plex} \cite{Wanless2002} of order $n$. The case $k=1$ corresponds to a {\em transversal} \cite{Colbourn2007} of a Latin square. Further, every $k$-plex of order $n$, with $k=2<n$ or $k>2$, determines a $k$-regular seminet with $n$ lines in all its parallel classes.
\item The problem of completing partial Latin squares, which is NP-complete \cite{Colbourn1984a}, has dealt with several structures: Ryser \cite{Ryser1951} analyzed the completion of partial Latin squares with pair of row and column structures $(s^r,r^s)$; Andersen and Hilton \cite{Andersen1997} studied those partial Latin squares of structure $((n-k)^n,(n-k)^n,(n-k)^n)$, for $k\in\{1,2\}$; more recently, Adams, Bryant and Buchanan \cite{Adams2008} dealt with the completion of those partial Latin squares with pair of row and column structure $(n^2 2^{n-2},n^22^{n-2})$.
\end{enumerate}

\vspace{0.2cm}

Let $\rho(z_1,z_2,z_3)$ be the number of partial Latin rectangles in $\mathcal{R}_{R,C,S}$ for any type $(R,C,S)\in\mathcal{T}_{r,m}\times \mathcal{T}_{s,m}\times \mathcal{T}_{n,m}$ such that $(z_R,z_C,z_S)=(z_1,z_2,z_3)\in \mathcal{Z}_{r,m}\times \mathcal{Z}_{s,m}\times \mathcal{Z}_{n,m}$.

\begin{theo}\label{thm_cs_0} Let $t$ and $n$ be two positive integers. Then,
$$\frac{n!^t t!^n}{t^{tn}}\leq \rho(t^n,t^n,n^t).$$
\end{theo}

{\bf Proof.}  Let $T=(t,\ldots,t)\in\mathcal{T}_{n,tn}$. Every partial Latin square $P\in\mathcal{R}_{n,n,n}$ of row and column type $T$ can be identified with a proper $n$-edge-colouring of the $t$-regular bipartite graph having the shape of $P$ as bi-adjacency matrix. To this end, an edge $ij$ of this graph is coloured according to a symbol $k$ if and only if $(i,j,k)\in E(P)$. The number of distinct partial Latin squares having $T$ as row and column types coincides, therefore, with that of distinct $n$-edge-colourings over the set of bipartite graphs with bi-adjacency matrix having $T$ as row and column sum vectors. According to Wei \cite{Wei1982}, this set has at least $n!^\mathrm{t}/\mathrm{t}!^n$ bipartite graphs. Further, Corollary 1d in \cite{Schrijver1998} involves every $t$-regular bipartite graph with $2n$ vertices to have at least $t!^{2n}/t^{tn}$ different $t$-edge-colourings. The result follows from combining both inequalities. \hfill $\Box$

\vspace{0.2cm}

\begin{lem} \label{lemma_type_11} Let $r'$, $s'$ and $n'$ be three positive integers greater than or equal to $r$, $s$ and $n$, respectively, and let $(z_1, z_2, z_3)\in \mathcal{Z}_{r',m}\times \mathcal{Z}_{s',m}\times \mathcal{Z}_{n',m}$. Let $(R,C,S)$ be a tuple in $\mathcal{T}_{r,m}\times \mathcal{T}_{s,m}\times \mathcal{T}_{n,m}$ such that $(z_R,z_C,z_S)=(z_1,z_2,z_3)$. Then, $|\mathcal{R}_{R,C,S}|=\rho(z_1,z_2,z_3)$.
\end{lem}

{\bf Proof.}  This result follows straightforward from the fact that the zero components in a tuple do not have any influence on the number of partial Latin rectangles that have this tuple as row, column or symbol type.\hfill $\Box$

\vspace{0.1cm}

\begin{prop}\label{prop_type_1} The next equality holds
$$|\mathcal{R}_{r,s,n:m}|=\sum_{\substack{r'\leq r\\s'\leq s\\n'\leq n}}\sum_{\substack{z_1\in \mathcal{Z}_{r',m}\\z_2\in\mathcal{Z}_{s',m}\\ z_3\in\mathcal{Z}_{n',m}}}\frac {r'!s'!n'!}{\prod_{i,j,k\leq m} d^{z_1}_i!d^{z_2}_j!d^{z_3}_k!}\binom{r}{r'}\binom{s}{s'} \binom{n}{n'}\rho(z_1,z_2,z_3),$$
where $d_i^{z_j}$ is the number of occurrences of the non-negative integer $i\leq m$ in any tuple of structure $z_j$, for each $j\leq 3$.
\end{prop}

{\bf Proof.} The result holds from Lemmas \ref{lemma_type_1} and \ref{lemma_type_11} and the number of tuples with a given structure. \hfill $\Box$

\vspace{0.2cm}

Table \ref{table5} shows the values of $\rho(z_R,z_C,z_S)$ for all $(R,C,S)\in \mathcal{T}_{r,m}\times \mathcal{T}_{s,m}\times \mathcal{T}_{n,m}$ such that $r\leq s\leq n\leq 6$ and $m\leq n$. Parastrophisms involve these values to be preserved under permutations of the components of the triple $(z_R,z_C,z_S)$. The corresponding distribution into isotopism ($IC$) and main ($MC$) classes of $\mathcal{R}_{r,s,n:m}$ is also indicated. The computation of these values has been determined by implementing Theorem \ref{thm_RCS2} in a procedure {\em PLRCS} in {\sc Singular}, which has been included in the previously mentioned library {\em pls.lib}. Proposition \ref{prop_type_1} has then be used to check the data exposed in Tables \ref{table1}--\ref{table4}.

\renewcommand{\tabcolsep}{1pt}
\begin{table}[htb]
{
\caption{Distribution into isotopism and main classes of the set $\mathcal{R}_{R,C,S}$.}
\label{table5}
\begin{center}
\resizebox{\textwidth}{!}{
\begin{tabular}{ccccc}
\begin{tabular}{llllrrr} \hline
$m$ & \, $z_R$ & \, $z_C$ & $z_S$ & $\rho$ & \, {\tiny IC} & \, {\tiny MC}\\ \hline
1 & $1$ & $1$ & $1$ & 1 & 1 & 1\\
2 & $2$ & $1^2$ & $1^2$ & 2 &1 & 1 \\
\ & $1^2$ & $1^2$ & $1^2$ & 4 & 1& 1 \\
3 & $3$ & $1^3$ & $1^3$ & 6 & 1& 1 \\
\ & $21$ & $21$ & $21$ & 1 &1 & 1\\
\ & \ & \  & $1^3$ & 6 & 1& 1\\
\ & & $1^3$ & $1^3$ & 18 & 1& 1\\
\ & $1^3$ & $1^3$ & $1^3$ & 36 & 1& 1 \\
4 & $4$ & $1^4$ & $1^4$ & 24 & 1& 1 \\
\ & $31$ & $21^2$ & $21^2$ & 4 & 1& 1 \\
\ & \ & \ & $1^4$ & 24 &1 & 1\\
\ & \ & $1^4$ & $1^4$ & 96 & 1& 1\\
\ & $2^2$ & $2^2$ & $2^2$ & 2 & 1& 1 \\
\ & \ & \ & $21^2$ & 4 & 1& 1\\
\ & \ & \ & $1^4$ & 24 & 1& 1\\
\ & \ & $21^2$ & $21^2$ & 12 & 2& 2 \\
\ & \ & \ & $1^4$ & 48 & 1& 1\\
\ & \ & $1^4$ & $1^4$ & 144 & 1& 1 \\
\ & $21^2$ & $21^2$ & $21^2$ & 40 & 5 & 3\\
\ & \ & \ & $1^4$ & 120 & 2& 2\\
\ & \ & $1^4$ & $1^4$ & 288 & 1& 1 \\
\ & $1^4$ & $1^4$ & $1^4$ & 576 & 1& 1\\
5 & $5$ & $1^5$ & $1^5$ & 120 & 1& 1 \\
\ & $41$ & $21^3$ & $21^3$ & 18 & 1& 1\\
\ & \ & \ & $1^5$ & 120 & 1& 1\\
\ & \ & $1^5$ & $1^5$ & 600 & 1& 1 \\
\ & $32$ & $2^21$ & $2^21$ & 6 &2 & 2 \\
\ & \ & \ & $21^3$ & 24 &2 & 2\\
\ & \ & \ & $1^5$ & 120 & 1& 1\\
\ & \ & $21^3$ & $21^3$ & 90 &3 & 3 \\
\ & \ & \ & $1^5$ & 360 & 1& 1\\
\ & \ & $1^5$ & $1^5$ & 1,200 & 1& 1 \\
\ & $31^2$ & $31^2$ & $2^21$ & 4 & 1 & 1\\
\ & \ & \ & $21^3$ & 24 & 1& 1\\
\ & \ & \ & $1^5$ & 120 & 1& 1\\
\ & \ & $2^21$ & $2^21$ & 12 & 2& 2 \\
\ & \ & \ & $21^3$ & 60 & 3& 3\\
\ & \ & \ & $1^5$ & 240 & 1& 1\\
\ & \ & $21^3$ & $21^3$ & 252 & 5& 4 \\
\ & \ & \ & $1^5$ & 840 & 2& 2\\
\ & \ & $1^5$ & $1^5$ & 2,400 & 1& 1 \\
\ & $2^21$ & $2^21$ & $2^21$ & 58 & 8 & 4\\
\ & \ & \ & $21^3$ & 180 & 8& 6\\
\ & \ & \ & $1^5$ & 600 & 2& 2\\
\ & \ & $21^3$ & $21^3$ & 504 & 8& 6 \\
\ & \ & \ & $1^5$ & 1,440 & 2& 2\\
\end{tabular}
& \hspace{0.25cm} &
\begin{tabular}{llllrrr} \hline
$m$ & \, $z_R$ & \, $z_C$ & $z_S$ & $\rho$ & \, {\tiny IC} & \, {\tiny MC}\\ \hline
5 & $2^21$ & $1^5$ & $1^5$ & 3,600 & 1& 1\\
\ & $21^3$ & $21^3$ & $21^3$ & 1,296& 8& 4\\
\ & \ & \ & $1^5$ & 3,240& 2& 2\\
\ & \ & $1^5$ & $1^5$ & 7,200& 1& 1 \\
\ & $1^5$ & $1^5$ & $1^5$ & 14,400& 1 & 1\\
6 & $6$ & $1^6$ & $1^6$ & 720& 1& 1\\
\ & $51$ & $21^4$ & $21^4$ & 96& 1& 1 \\
\ & \ & \ & $1^6$ & 720& 1& 1\\
\ & \ & $1^6$ & $1^6$ & 4,320& 1& 1\\
\ & $42$ & $2^21^2$ & $2^21^2$ & 28& 3& 3\\
\ & \ & \ & $21^4$ & 144& 2& 2\\
\ & \ & \ & $1^6$ & 720& 1& 1\\
\ & \ & $21^4$ & $21^4$ & 672& 3& 3 \\
\ & \ & \ & $1^6$ & 2,880& 1& 1\\
\ & \ & $1^6$ & $1^6$ & 10,800& 1& 1\\
\ & $3^2$ & $2^3$ & $2^3$ & 12& 1& 1 \\
\ & \ & \ & $2^21^2$ & 36& 2& 2\\
\ & \ & \ & $21^4$ & 144& 1& 1\\
\ & \ & \ & $1^6$ & 720& 1& 1\\
\ & \ & $2^21^2$ & $2^21^2$ & 88& 5& 4\\
\ & \ & \ & $21^4$ & 336& 3& 3\\
\ & \ & \ & $1^6$ & 1,440& 1& 1\\
\ & \ & $21^4$ & $21^4$ & 1,152& 2& 2 \\
\ & \ & \ & $1^6$ & 4,320& 1& 1\\
\ & \ & $1^6$ & $1^6$ & 14,400& 1& 1\\
\ & $41^2$ & $31^3$ & $2^21^2$ & 24& 1& 1\\
\ & \ & \ & $21^4$ & 144& 1& 1\\
\ & \ & \ & $1^6$ & 720& 1& 1\\
\ & \ & $2^21^2$ & $2^21^2$ & 56& 3& 3\\
\ & \ & \ & $21^4$ & 336& 3& 3\\
\ & \ & \ & $1^6$ & 1,440& 1& 1\\
\ & \ & $21^4$ & $21^4$ & 1,728& 5& 4 \\
\ & \ & \ & $1^6$ & 6,480& 2& 2\\
\ & \ & $1^6$ & $1^6$ & 21,600& 1& 1\\
\ & $321$ & $321$ & $321$ & 1& 1& 1 \\
\ & \ &  & $31^3$ & 6& 1& 1\\
\ & \ &  & $2^3$ & 12& 2& 2\\
\ & \ &  & $2^21^2$ & 40& 10 & 7\\
\ & \ &  & $21^4$ & 168& 7& 5\\
\ & \ &  & $1^6$ & 720& 1& 1\\
\ & \ & $31^3$ & $31^3$ & 36& 1& 1 \\
\ & \ &  & $2^3$ & 36& 1& 1\\
\ & \ &  & $2^21^2$ & 144& 6& 6\\
\ & \ &  & $21^4$ & 576& 5& 5\\
\ & \ & \ & $1^6$ & 2,160& 1& 1\\
\ & \ & \ & $2^3$ & 36& 1& 1\\
\end{tabular}
& \hspace{0.25cm} &
\begin{tabular}{llllrrr} \hline
$m$ & \, $z_R$ & $z_C$ & $z_S$ & $\rho$ & \, {\tiny IC} & \, {\tiny MC}\\ \hline
6 & $321$ & $2^3$ & $2^21^2$ & 156& 7 & 7\\
\ & \ &  & $21^4$ & 576& 4& 4\\
\ & \ &  & $1^6$ & 2,160& 1& 1\\
\ & \ & $2^21^2$ & $2^21^2$ & 512& 33 & 20\\
\ & \ &  & $21^4$ & 1,728& 20& 20\\
\ & \ &  & $1^6$ & 5,760& 3& 3\\
\ & \ & $21^4$ & $21^4$ & 5,280& 15& 10\\
\ & \ &  & $1^6$ & 15,840& 3& 3\\
\ & \ & $1^6$ & $1^6$ & 43,200& 1& 1\\
\ & $2^3$ & $2^3$ & $2^3$ & 144& 2& 2 \\
\ & \ & \ & $31^3$ & 72& 1& 1\\
\ & \ &  & $2^21^2$ & 432& 5& 4\\
\ & \ &  & $21^4$ & 1,296& 2& 2\\
\ & \ &  & $1^6$ & 4,320& 1& 1\\
\ & \ & $31^3$ & $31^3$ & 144& 2& 2 \\
\ & \ &  & $2^21^2$ & 360& 3& 3\\
\ & \ &  & $21^4$ & 1,296& 2& 2\\
\ & \ &  & $1^6$ & 4,320& 1& 1\\
\ & \ & $2^21^2$ & $2^21^2$ & 1,260& 18& 13\\
\ & \ &  & $21^4$ & 3,600& 8& 8\\
\ & \ &  & $1^6$ & 10,800& 2& 2\\
\ & \ & $21^4$ & $21^4$ & 9,504& 4& 4 \\
\ & \ &  & $1^6$ & 25,920& 1& 1\\
\ & \ & $1^6$ & $1^6$ & 64,800& 1& 1\\
\ & $31^3$ & $31^3$ & $31^3$ & 216& 1 & 1\\
\ & \ &  & $2^21^2$ & 576& 5& 4\\
\ & \ & \ & $21^4$ & 2,160& 5& 4\\
\ & \ & & $1^6$ & 7,200& 2& 2\\
\ & \ & $2^21^2$ & $2^21^2$ & 1,344& 16& 11\\
\ & \ &  & $21^4$ & 4,320& 10& 10\\
\ & \ &  & $1^6$ & 12,960& 2& 2\\
\ & \ & $21^4$ & $21^4$ & 12,672& 8& 6\\
\ & \ &  & $1^6$ & 34,560& 2& 2\\
\ & \ & $1^6$ & $1^6$ & 86,400& 1& 1\\
\ & $2^21^2$ & $2^21^2$ & $2^21^2$ & 3,320& 62& 19\\
\ & \ &  & $21^4$ & 8,976& 29& 19\\
\ & \ &  & $1^6$ & 24,480& 5& 4\\
\ & \ & $21^4$ & $21^4$ & 22,464& 15& 11\\
\ & \ &  & $1^6$ & 56,160& 3& 3\\
\ & \ & $1^6$ & $1^6$ & 129,600& 1& 1\\
\ & $21^4$ & $21^4$ & $21^4$ & 52,416& 9& 5\\
\ & \ &  & $1^6$ & 120,960& 2& 2\\
\ & \ & $1^6$ & $1^6$ & 259,200& 1& 1\\
\ & $1^6$ & $1^6$ & $1^6$ & 518,400& 1& 1\\
\ & \ &  &   &  &  &  \\
\ & \ &  &   &  &  &  \\
\end{tabular}
\end{tabular}
}
\end{center}
}
\end{table}

\vspace{0.2cm}

Table \ref{table5} is also used in the next theorem to determine the number of partial Latin rectangles of size up to six. This generalizes a recent result \cite{Falcon2015} in which the case $m\leq 2$ was already exposed. In order to avoid an excessive length of the polynomials that appear in the theorem, the polynomial $\sum_{\sigma\in \mathrm{Sym}(\{a,b,c\})}r^a s^b n^c$ is denoted as $\overline{abc}$, for all $a,b,c\geq 0$, where $\mathrm{Sym}(\{a,b,c\})$ constitutes the set of permutations of the ordered set $\{a,b,c\}$. Thus, for instance, $3\ \overline{211}$ denotes the polynomial $3(r^2sn+rs^2n+rsn^2)$.

\begin{theo}\label{thm_rsnm} The next equalities hold
\begin{enumerate}[a)]
\item $|\mathcal{R}_{r,s,n:0}|=1$.
\item $|\mathcal{R}_{r,s,n:1}|=\overline{111}$.
\item $2!   |\mathcal{R}_{r,s,n:2}|=\overline{111}\ (\overline{111} - \overline{100} +2)$.
\item $3!   |\mathcal{R}_{r,s,n:3}|=\overline{111}\ (\overline{222}-3\ \overline{211}+6\ (\overline{111} + \overline{110}) + 2\ \overline{200}-12\ \overline{100}+14)$.
\item $4!   |\mathcal{R}_{r,s,n:4}|=\overline{111}\  (\overline{333} - 6\ \overline{322} + 12\ \overline{222} + 11\ \overline{311} + 30\ \overline{221} - 60\ \overline{211} - 6\ \overline{300} - 36\ \overline{210} - 28\ \overline{111} +72\ \overline{200} +198\ \overline{110} -228\ \overline{100} +198)$.
\item $5!   |\mathcal{R}_{r,s,n:5}|=\overline{111}\  (\overline{444} -10\ \overline{433} + 20\ \overline{333} +35\ \overline{422} +90\ \overline{332} -180\ \overline{322} -50\ \overline{411} -260\ \overline{321} -460\ \overline{222} +520\ \overline{311} +1,350\ \overline{221} + 24\ \overline{400} +240\ \overline{310} +480\ \overline{220} -320\ \overline{211} -480\ \overline{300} -2,520\ \overline{210} -5,090\ \overline{111} +2,880\ \overline{200} +7,440\ \overline{110} -6,360\ \overline{100} +4512)$.
\item $6!   |\mathcal{R}_{r,s,n:6}|=\overline{111}\ (\overline{555}-15\ \overline{544}+30\ \overline{444}+85\ \overline{533}+210\ \overline{443}-420\ \overline{433} -225\ \overline{522}-1,065\ \overline{432} - 2,150\ \overline{333} +2,130\ \overline{422} +5,310\ \overline{332} +274\ \overline{511}  +2,310\ \overline{421} +4,400\ \overline{331}+4,800\ \overline{322}-4,620\ \overline{411}- 22,170\ \overline{321} -49,500\ \overline{222} - 120\ \overline{500} - 1,800\ \overline{410} - 6,000\ \overline{320}  + 10,460\ \overline{311} +34,980\ \overline{221}
    + 3,600 \ \overline{400}+30,600\ \overline{310} +58,440\ \overline{220}
    +88,710  \ \overline{211} - 34,800\ \overline{300} - 165,480\ \overline{210} -364,268\ \overline{111}  + 140,040\ \overline{200}+344,520\ \overline{110}-240,720\ \overline{100} +146,400)$.
\end{enumerate}
\end{theo}

{\bf Proof.} The first equality is immediate. This is refereed to the partial Latin rectangle without any entry. The other equalities follow from Proposition \ref{prop_type_1} and Table \ref{table5}. We prove here in detail the first three expressions; the rest follows similarly. In the use of Table \ref{table5}, recall that the value $\rho(z_R,z_C,z_S)$ is preserved by parastrophism, that is, the placement of the structures $z_R$, $z_C$ and $z_S$ can be interchanged.

\begin{enumerate}
\item[b)] $|\mathcal{R}_{r,s,n:1}|=rsn\, \rho(1,1,1)=rsn$.
\item[c)] $|\mathcal{R}_{r,s,n:2}|=r\binom{s}{2}\binom{n}{2}\rho(2,1^2,1^2) + s\binom{r}{2} \binom{n}{2}\rho(1^2,2,1^2) + n\binom{r}{2}\binom{s}{2} \rho(1^2,1^2,2)+$ $\binom{r}{2}\binom{s}{2}\binom{n}{2}\rho(1^2,1^2,1^2) =\frac {rsn}2 (rsn-r-s-n+2)$.
\item[d)] $|\mathcal{R}_{r,s,n:3}|=r\binom{s}{3}\binom{n}{3}\rho(3,1^3,1^3) + s \binom{r}{3}\binom{n}{3}\rho(1^3,3,1^3) + n \binom{r}{3}\binom{s}{3} \rho(1^3,1^3,3) + 8\binom{r}{2}\binom{s}{2}\binom{n}{2} \rho(21,21,21) + 4\binom{r}{2}\binom{s}{2}\binom {n}{3} \rho(21,21,1^3) + 4\binom{r}{2}\binom{s}{3}\binom {n}{2} \rho(21,1^3,21) + 4\binom {r}{3}\binom{s}{2}\binom{n}{2} \rho(1^3,21,21) + 2\binom{r}{2} \binom {s}{3}\binom {n}{3} \rho(21,1^3,1^3) + 2\binom{r}{3}\binom {s}{2}\binom {n}{3} \rho(1^3,21,1^3) + 2\binom{r}{3}\binom {s}{3}\binom {n}{2}\ \rho(1^3,1^3,21) + \binom {r}{3}\binom {s}{3}\binom {n}{3} \rho(1^3,1^3,1^3) = \frac {rsn}6 (r^2s^2n^2-3r^2sn-3rs^2n-3rsn^2+6rsn+6rs+6rn+6sn+2r^2+2s^2+ 2n^2-12r- 12s-12n+14)$.

    \hfill $\Box$
\end{enumerate}

\begin{corol}\label{crl_rsnm} Let $n$ be a positive integer. Then
\begin{enumerate}[a)]
\item $|\mathcal{R}_{n,n,n:0}|=1$.
\item $|\mathcal{R}_{n,n,n:1}|=n^3$.
\item $2! \  |\mathcal{R}_{n,n,n:2}|=n^3   (n-1)^2   (n+2)$.
\item $3! \  |\mathcal{R}_{n,n,n:3}|=n^3   (n-1)^2   (n^4+2n^3-6n^2-8n+14)$.
\item $4! \  |\mathcal{R}_{n,n,n:4}|=n^3  (n-1)^2 (n^7+2n^6-15n^5-20n^4+98n^3+36n^2-288n+198)$.
\item $5! \  |\mathcal{R}_{n,n,n:5}|=n^3 (n-1)^2 (n-2)^2 (n^8+6n^7-7n^6-88n^5+6n^4+532n^3-84n^2+1386n+1128)$.
\item $6! \  |\mathcal{R}_{n,n,n:6}|=n^3(n-1)^2(n-2)^2 (n^{11}+6n^{10}-22n^9-168n^8+231n^7+2,022n^6-2,014n^5-12,606n^4+16,168n^3+ 32,250n^2-70,740n+36,600)$.
\end{enumerate}
\end{corol}

{\bf Proof.}  This result follows straightforward from Theorem \ref{thm_rsnm} once we impose $r=s=n$.\hfill $\Box$

\section{Classification of seminets with low point rank}

Every seminet is equivalent to a non-compressible regular partial Latin square \cite{Stojakovic1979}. The next lemma follows straightforward from the definition of compressibility and regularity of partial Latin squares and indicates how both properties can be expressed in terms of types of partial Latin squares.

\begin{lem}\label{lmm1} Let $R=(\mathrm{r}_1,\ldots,\mathrm{r}_n)$, $C=(\mathrm{c}_1,\ldots,\mathrm{c}_n)$ and $S=(\mathrm{s}_1,\ldots,\mathrm{s}_n)$ be three tuples in $\mathcal{T}_{n,m}$ and let $P$ be a partial Latin square in $\mathcal{R}_{R,C,S}$. Then,
\begin{enumerate}
\item $P$ is non-compressible if and only if at least one of its row, column or symbol types does not have zero components.
\item $P$ is regular if and only if the next three conditions hold.
\begin{enumerate}
\item The cell $(i,j)$ of $P$ is empty for all $i,j\leq n$ such that $\mathrm{r}_i=\mathrm{c}_j=1$.
\item $\mathrm{s}_k>1$ for all $i,j\leq n$ such that $\mathrm{r}_i=1$ and $(i,j,k)\in E(P)$.
\item $\mathrm{s}_k>1$ for all $i,j\leq n$ such that $\mathrm{c}_j=1$ and $(i,j,k)\in E(P)$. \hfill $\Box$
\end{enumerate}
\end{enumerate}
\end{lem}

\vspace{0.5cm}

Let $\mathcal{R}^{\text{reg}}_{R,C,S}$ be the set of regular partial Latin squares whose row, column and symbol types coincide, respectively, with $R$, $C$ and $S$. Since regularity is preserved by paratopism of partial Latin squares, the cardinality of this set only depends on the structures of $R$, $C$ and $S$. The next result shows how this cardinality is immediately determined for certain structures. Recall that each exponent $d^z_i$ in the structure $z=m^{d^z_m}\ldots 1^{d^z_1}$ is the number of occurrences of a given non-negative integer $i$ as a component of any tuple of structure $z$.

\begin{prop}\label{prp_reg} Let $z_1$, $z_2$ and $z_3$ be three structures of weight $m$. Then,
\begin{enumerate}[a)]
\item If $d_1^{z_1}=d_1^{z_2}=0$, then every partial Latin square having two of their row, column or symbol structures equal to $z_1$ and $z_2$, respectively, is regular.
\item If $d_1^{z_1} + d_1^{z_2}+ d_1^{z_3}>m$, then no partial Latin square of structure $(z_1, z_2, z_3)$ is regular.
\end{enumerate}
\end{prop}

{\bf Proof.} None partial Latin rectangle in (a) contains a row or a column with exactly one entry. All of them are, therefore, regular. Further, from the definition of regularity, assertion (b) holds because every regular partial Latin rectangle of type $(z_1,z_2,z_3)$ satisfies that $d_1^{z_1}+d_1^{z_2}\leq \sum_{i=2}^m d_i^{z_3} = m - d_1^{z_3}$ and hence, $d_1^{z_1} + d_1^{z_2}+ d_1^{z_3}\leq m$.\hfill $\Box$

\vspace{0.4cm}

The next result indicates how computational algebraic geometry can be used to determine the set $\mathcal{R}^{\text{reg}}_{R,C,S}$.

\begin{theo} \label{thm_reg} Let $R=(\mathrm{r}_1,\ldots,\mathrm{r}_n)$, $C=(\mathrm{c}_1,\ldots,\mathrm{c}_n)$ and $S=(\mathrm{s}_1,\ldots,\mathrm{s}_n)$ be three tuples in $\mathcal{T}_{n,m}$ and let $p$ be the first prime greater than the maximum of all the components of $R$, $C$ and $S$. The set $\mathcal{R}^{\text{reg}}_{R,C,S}$ is identified with the set of zeros of the zero-dimensional radical ideal
$$I^{\text{reg}}_{R,C,S}=I_{R,C,S} + \langle\,x_{ijk}\colon\, i,j,k\leq n, \mathrm{r}_i=\mathrm{c}_j=1\,\rangle +$$ $$\langle\,x_{ijk}\colon\, i,j,k\leq n, \mathrm{r}_i=\mathrm{s}_k=1\,\rangle+\langle\,x_{ijk}\colon\, i,j,k\leq n, \mathrm{c}_j=\mathrm{s}_k=1\,\rangle\subset \mathbb{F}_p[X].$$
Besides, $|\mathcal{R}^{\text{reg}}_{R,C,S}|= \mathrm{dim}_{\mathbb{F}_p}(\mathbb{F}_p[X]/I^{\text{reg}}_{R,C,S})$.
\end{theo}

{\bf Proof.} Since $I^{\text{reg}}_{R,C,S}\subseteq I_{R,C,S}$, each zero of the ideal $I^{\text{reg}}_{R,C,S}$ is uniquely related to a partial Latin square whose row, column and symbol types coincide, respectively, with $R$, $C$ and $S$. The rest of the proof is similar to that of Theorem \ref{thm0} once we observe that the three subideals that are added to $I_{R,C,S}$ in the definition of $I^{\text{reg}}_{R,C,S}$ involve these partial Latin squares to verify, respectively, conditions (2.a), (2.b) and (2.c) of Lemma \ref{lmm1}. \hfill $\Box$

\vspace{0.5cm}

\begin{figure*}[ht]
  \begin{center}
  \begin{tabular}{ccc}
   \begin{tabular}{c} \includegraphics[width=.15\textwidth]{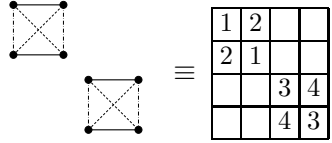}
    \end{tabular} & $\equiv$ & {\small \begin{tabular}{|c|c|c|c|}\hline
   1 & 2 & \ & \ \\ \hline
   2 & 1 & \ & \ \\ \hline
   \ & \ & 3 & 4 \\ \hline
   \ & \ & 4 & 3\\ \hline
\end{tabular}}
  \end{tabular}
  \end{center}
  \caption{Seminet of point rank $8$ that is not a configuration.}  \label{FigSem8}
\end{figure*}

\begin{figure}{\scriptsize
\begin{center}
  \begin{tabular}{cccccccc}
    \begin{tabular}{c}
    \\
    \includegraphics[width=.08\textwidth]{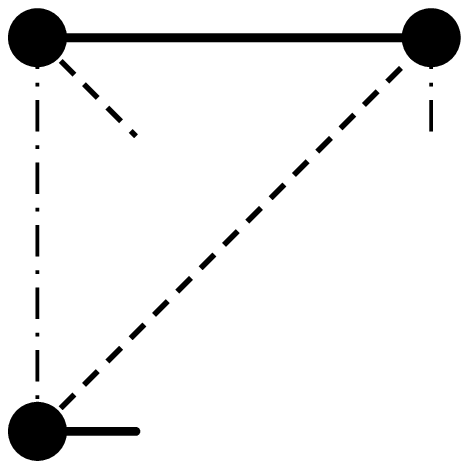}
    \end{tabular} &
    \begin{tabular}{c}
    \\
    \includegraphics[width=.1\textwidth]{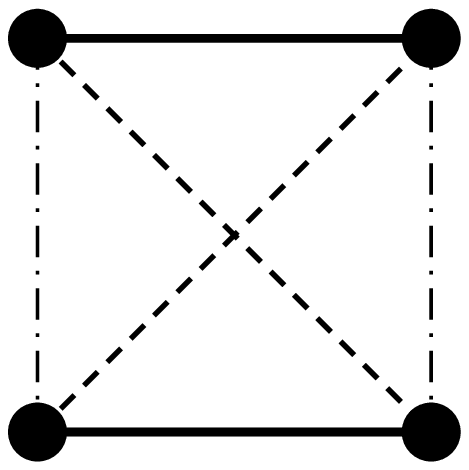}
    \end{tabular} &
    \begin{tabular}{c}
    \\
    \includegraphics[width=.1\textwidth]{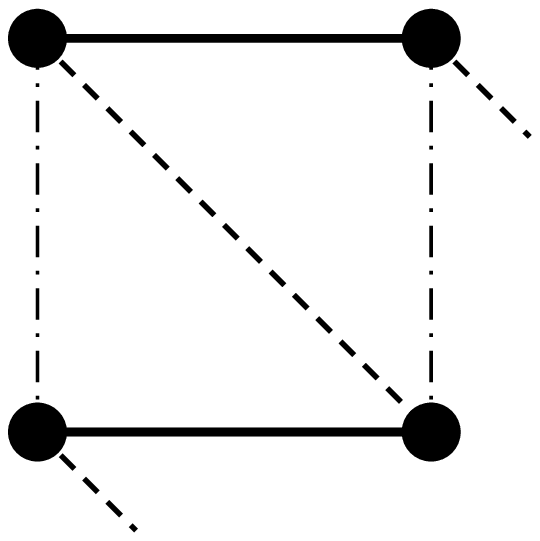}
    \end{tabular} &
    \begin{tabular}{c}
    \\
    \includegraphics[width=.1\textwidth]{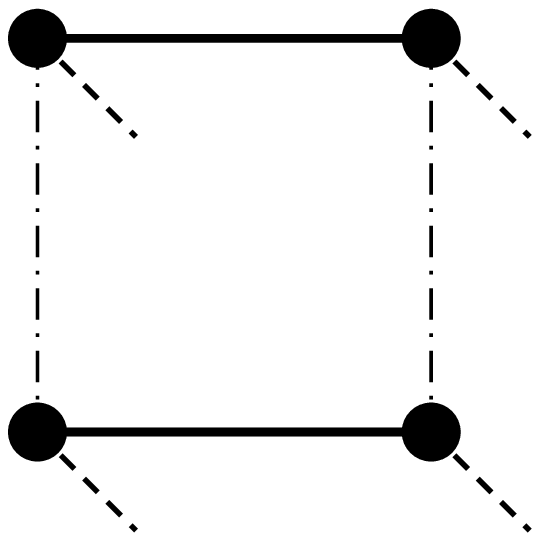}
    \end{tabular} &
    \begin{tabular}{c}
    \\
    \includegraphics[width=.14\textwidth]{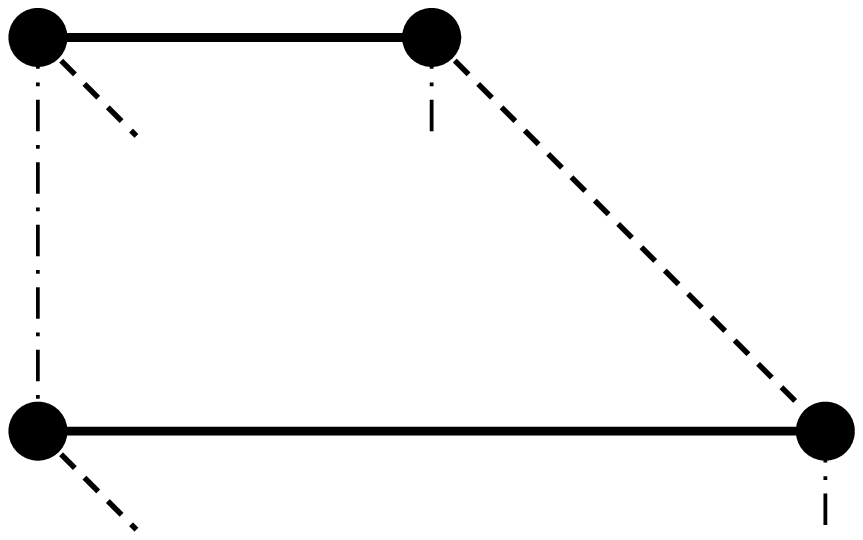}
    \end{tabular} &
    \begin{tabular}{c}
    \\
    \includegraphics[width=.14\textwidth]{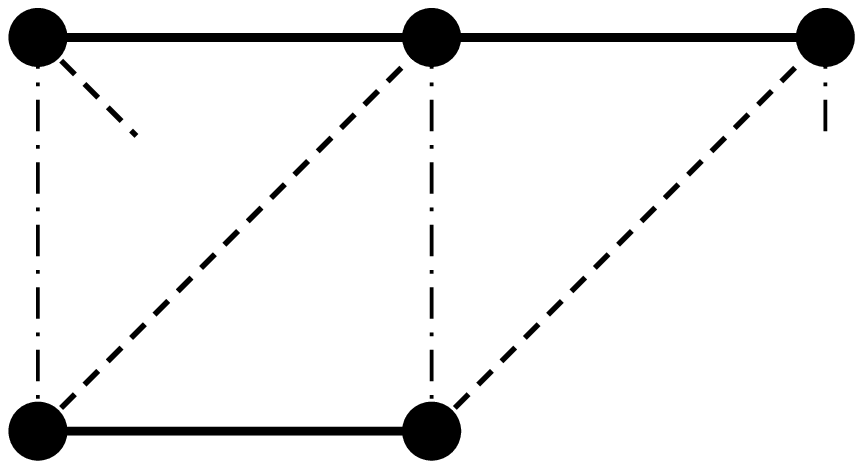}
    \end{tabular}
    \\
    $\mathcal{S}_3$ & $\mathcal{S}_{4,1}$ & $\mathcal{S}_{4,2}$ & $\mathcal{S}_{4,3}$ & $\mathcal{S}_{4,4}$ & $\mathcal{S}_{5,1}$\\
    \begin{tabular}{c}
    \\
    \includegraphics[width=.14\textwidth]{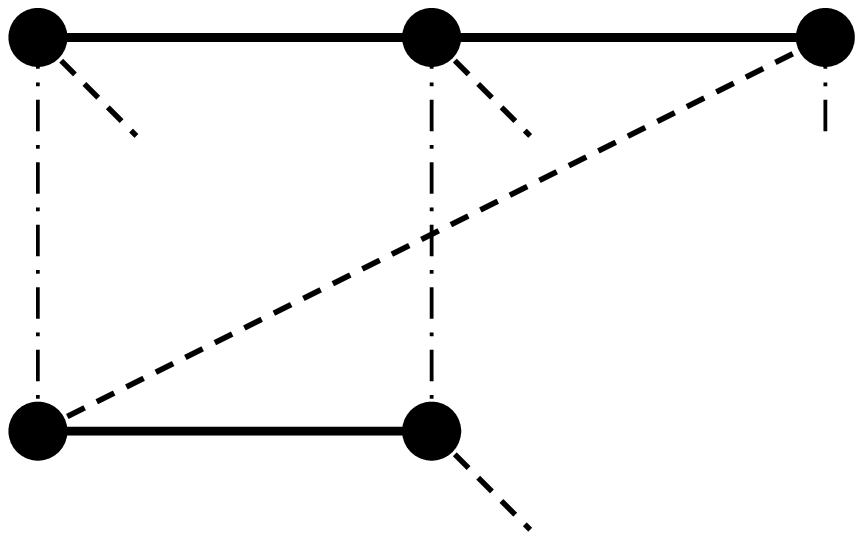}
    \end{tabular} &
    \begin{tabular}{c}
    \\
    \includegraphics[width=.1\textwidth]{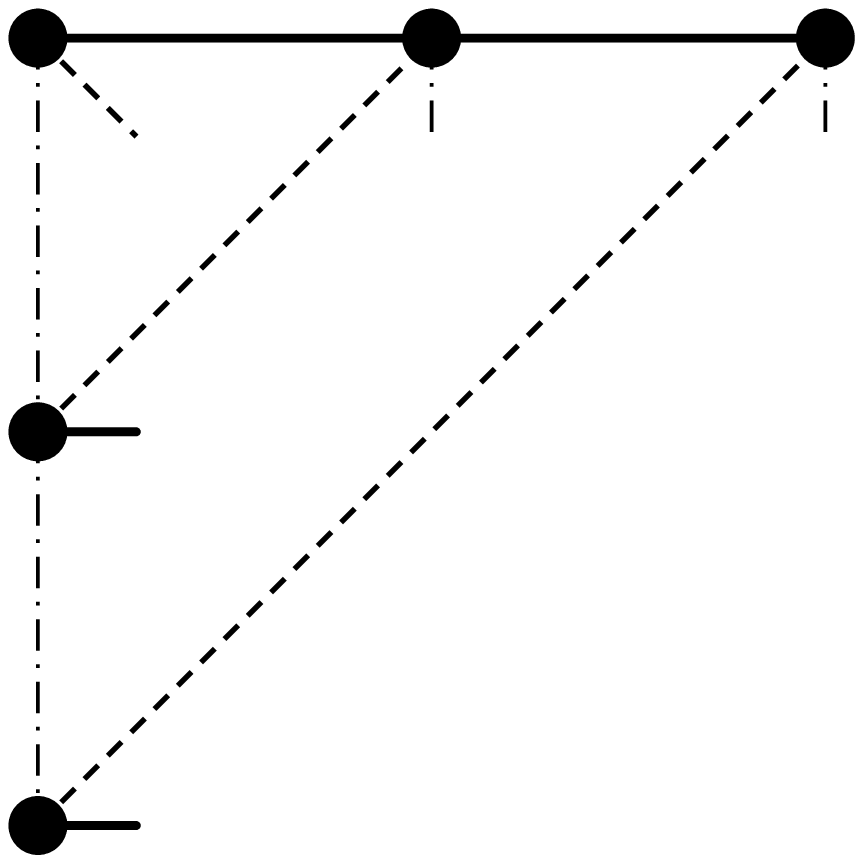}
    \end{tabular} &
    \begin{tabular}{c}
    \\
    \includegraphics[width=.1\textwidth]{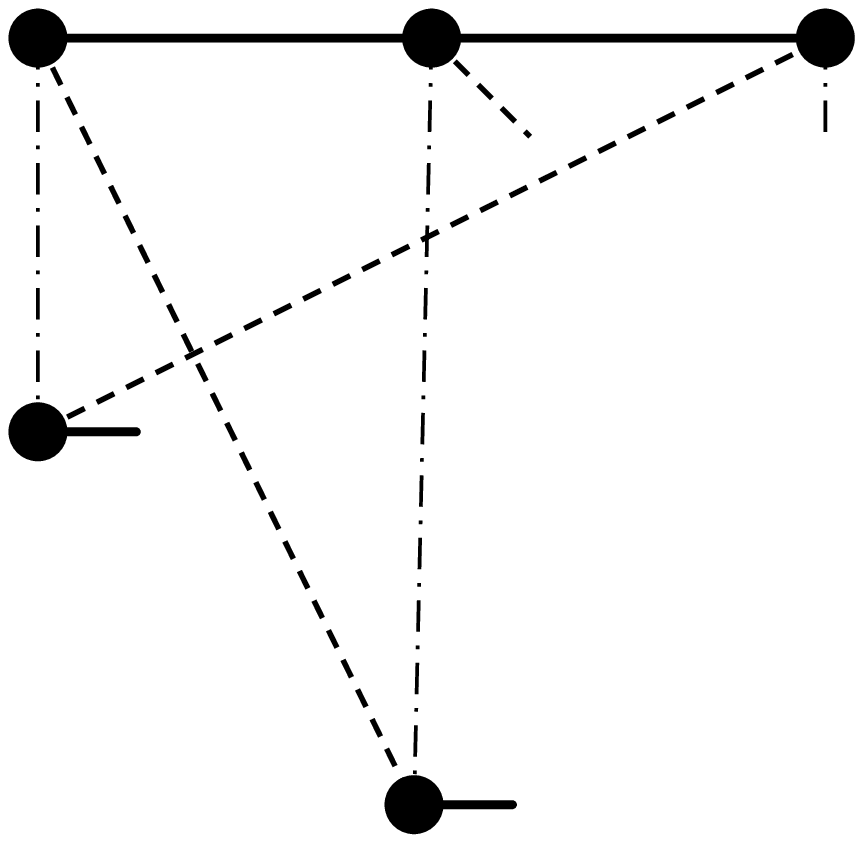}
    \end{tabular} &
    \begin{tabular}{c}
    \\
    \includegraphics[width=.1\textwidth]{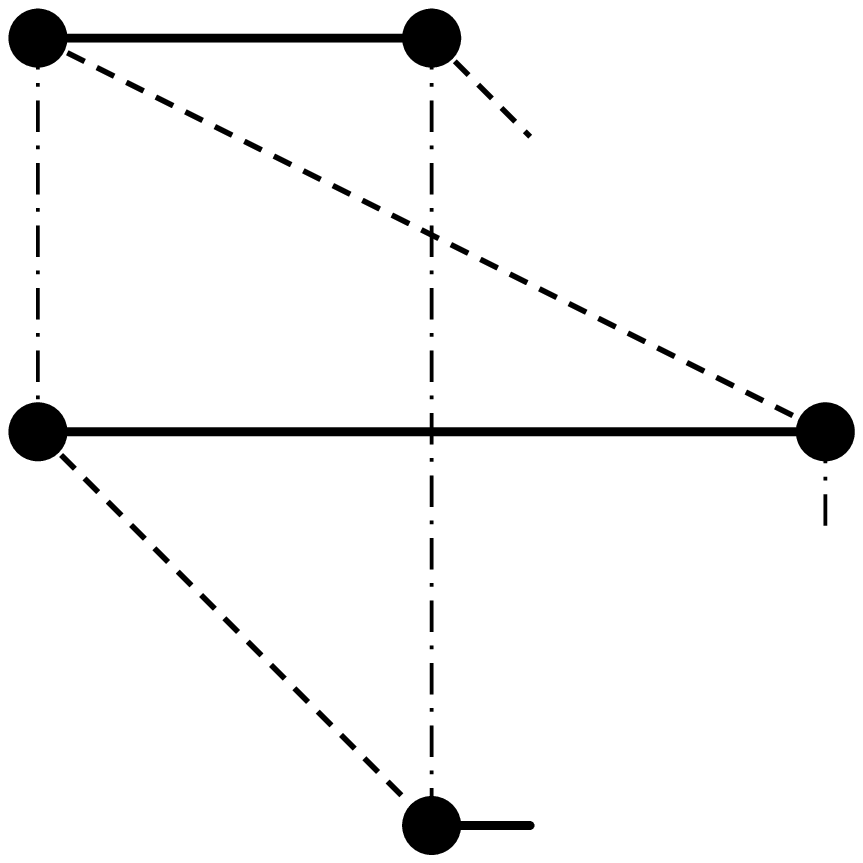}
    \end{tabular} &
    \begin{tabular}{c}
    \\
    \includegraphics[width=.1\textwidth]{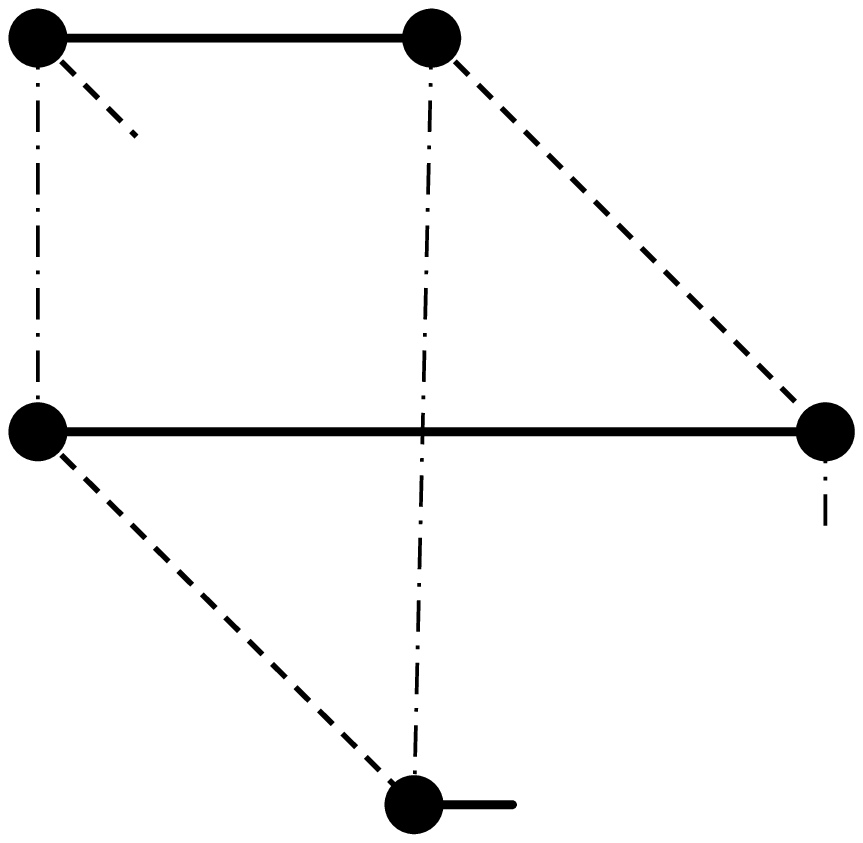}
    \end{tabular} &
    \begin{tabular}{c}
    \\
    \includegraphics[width=.1\textwidth]{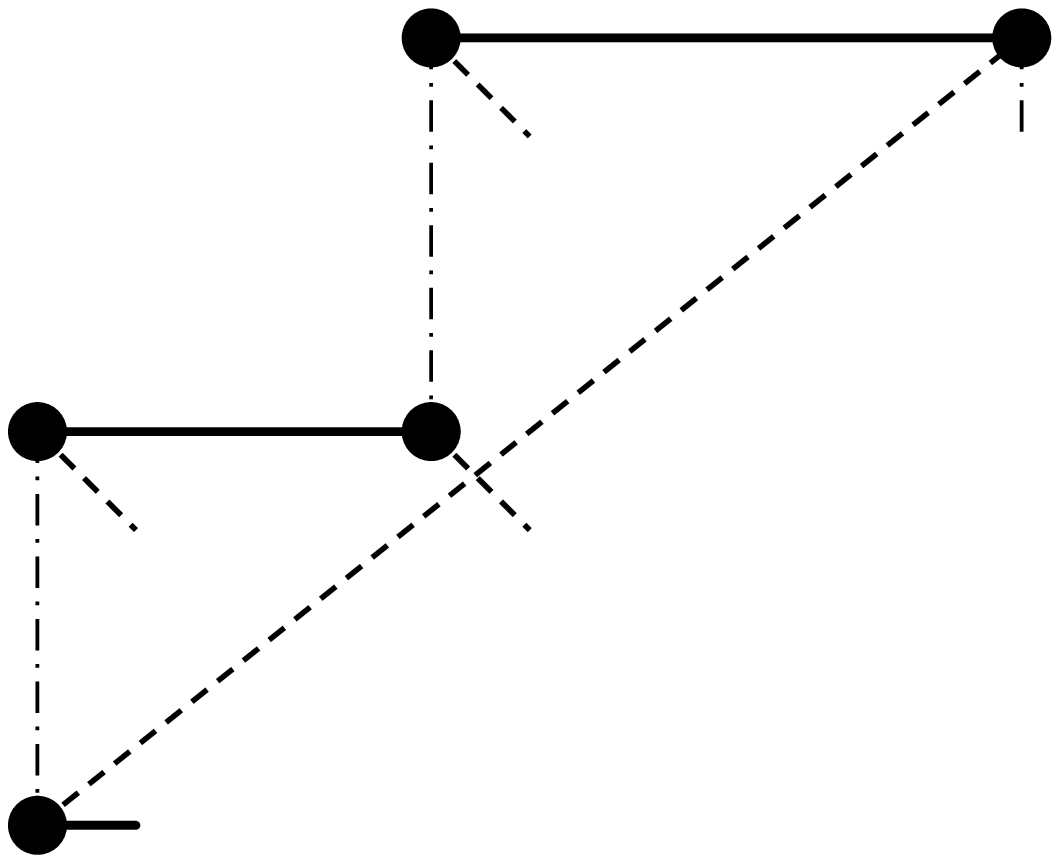}
    \end{tabular}
    \\
    $\mathcal{S}_{5,2}$ & $\mathcal{S}_{5,3}$ & $\mathcal{S}_{5,4}$ & $\mathcal{S}_{5,5}$ & $\mathcal{S}_{5,6}$ & $\mathcal{S}_{5,7}$\\
    \begin{tabular}{c}
    \\
    \includegraphics[width=.14\textwidth]{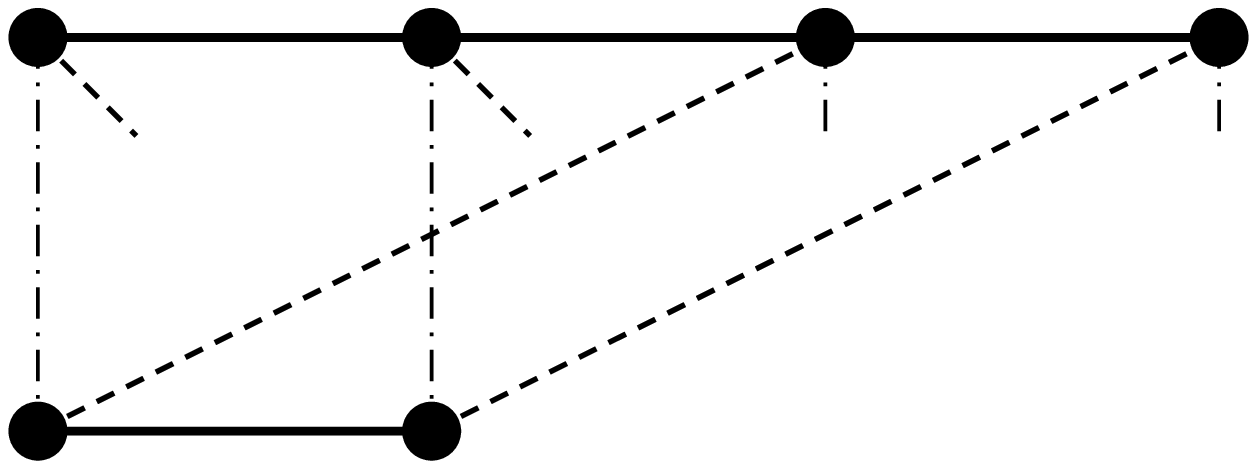}
    \end{tabular} &
    \begin{tabular}{c}
    \\
    \includegraphics[width=.14\textwidth]{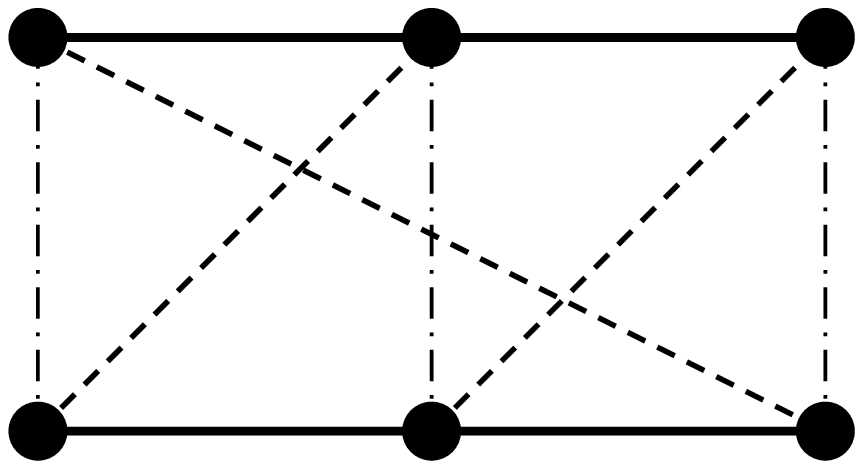}
    \end{tabular} &
    \begin{tabular}{c}
    \\
    \includegraphics[width=.14\textwidth]{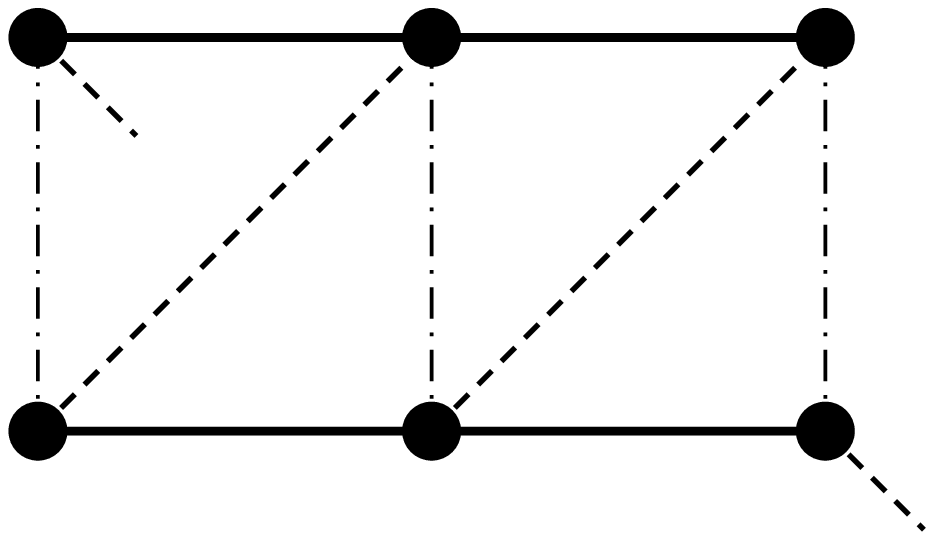}
    \end{tabular} &
    \begin{tabular}{c}
    \\
    \includegraphics[width=.14\textwidth]{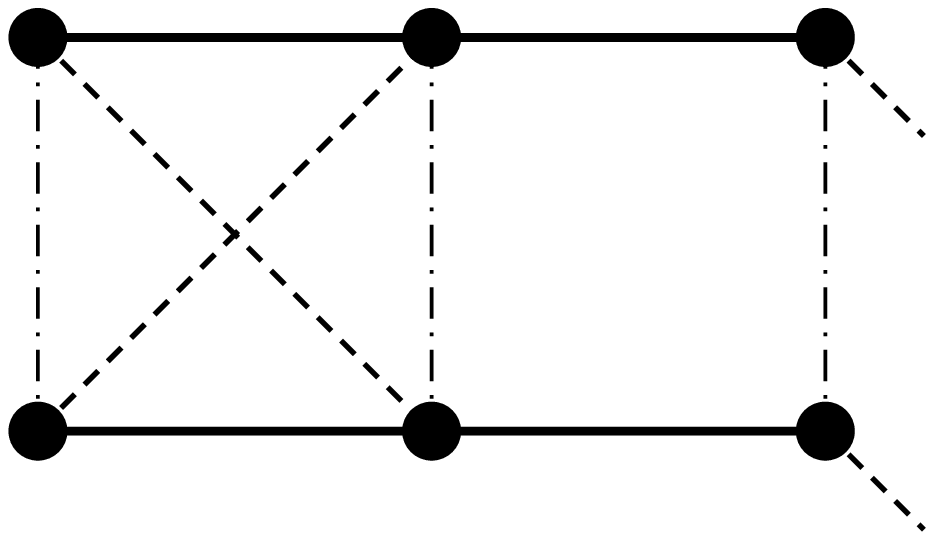}
    \end{tabular} &
    \begin{tabular}{c}
    \\
    \includegraphics[width=.14\textwidth]{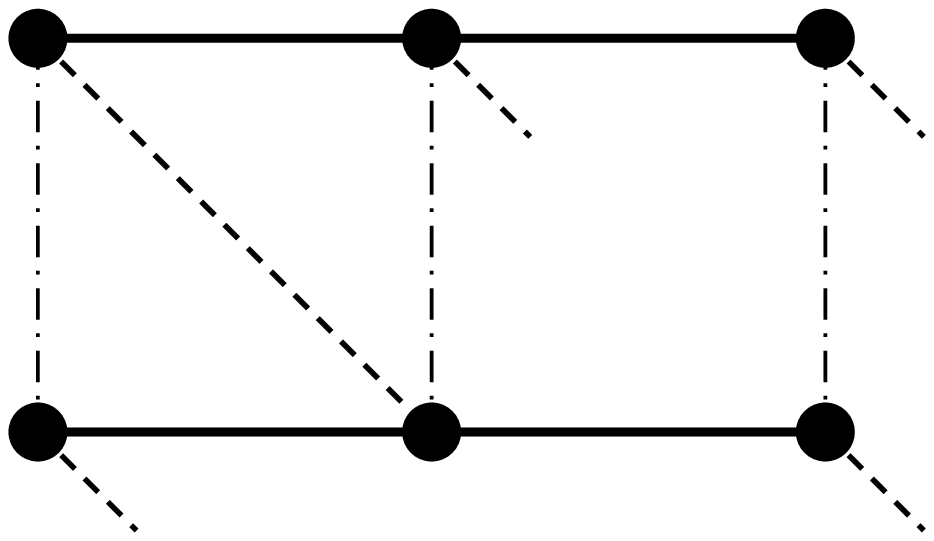}
    \end{tabular} &
    \begin{tabular}{c}
    \\
    \includegraphics[width=.14\textwidth]{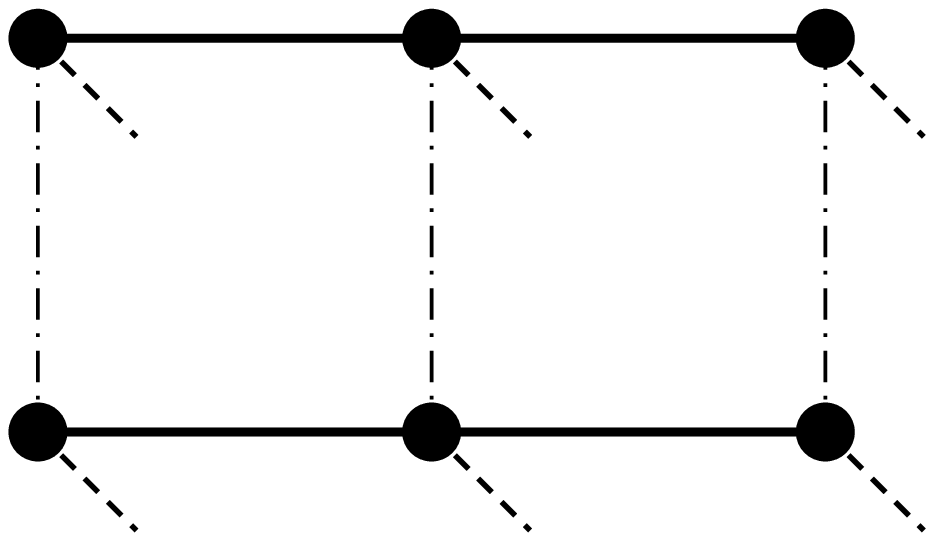}
    \end{tabular}
    \\
    $\mathcal{S}_{6,1}$ & $\mathcal{S}_{6,2}$ & $\mathcal{S}_{6,3}$ & $\mathcal{S}_{6,4}$ & $\mathcal{S}_{6,5}$ & $\mathcal{S}_{6,6}$
    \\
    \begin{tabular}{c}
    \includegraphics[width=.14\textwidth]{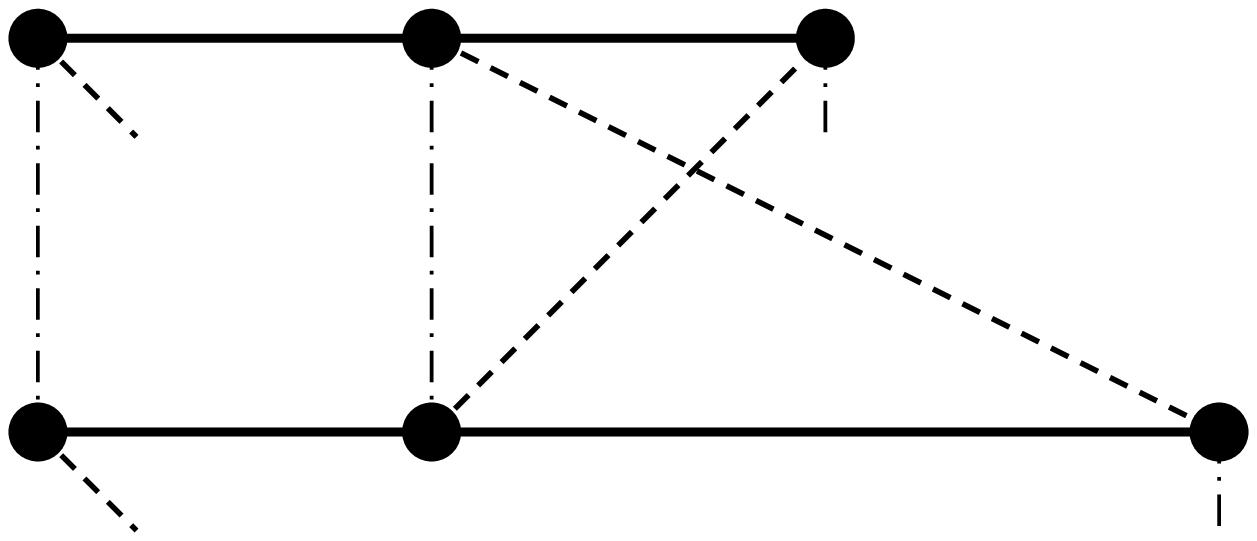}
    \end{tabular} &
    \begin{tabular}{c}
    \\
    \includegraphics[width=.14\textwidth]{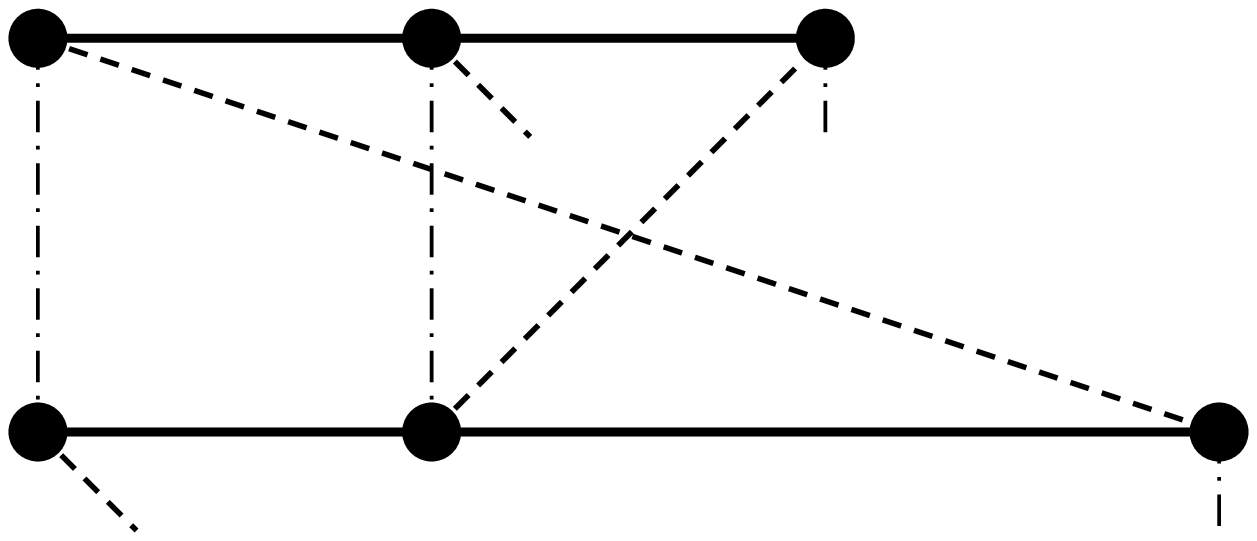}
    \end{tabular} &
    \begin{tabular}{c}
    \\
    \includegraphics[width=.14\textwidth]{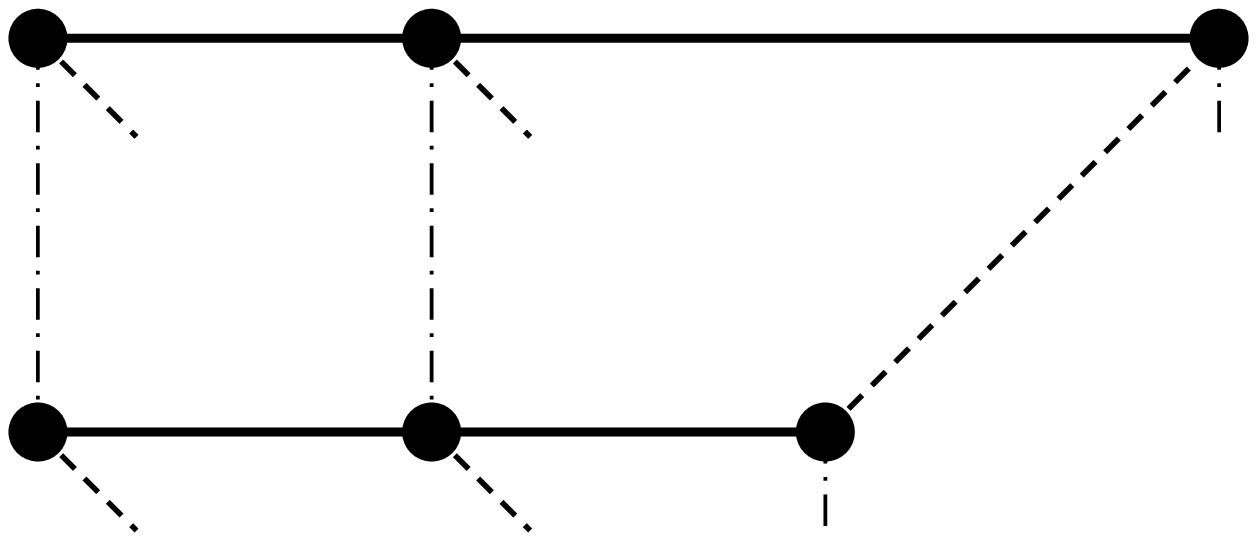}
    \end{tabular} &
    \begin{tabular}{c}
    \\
    \includegraphics[width=.14\textwidth]{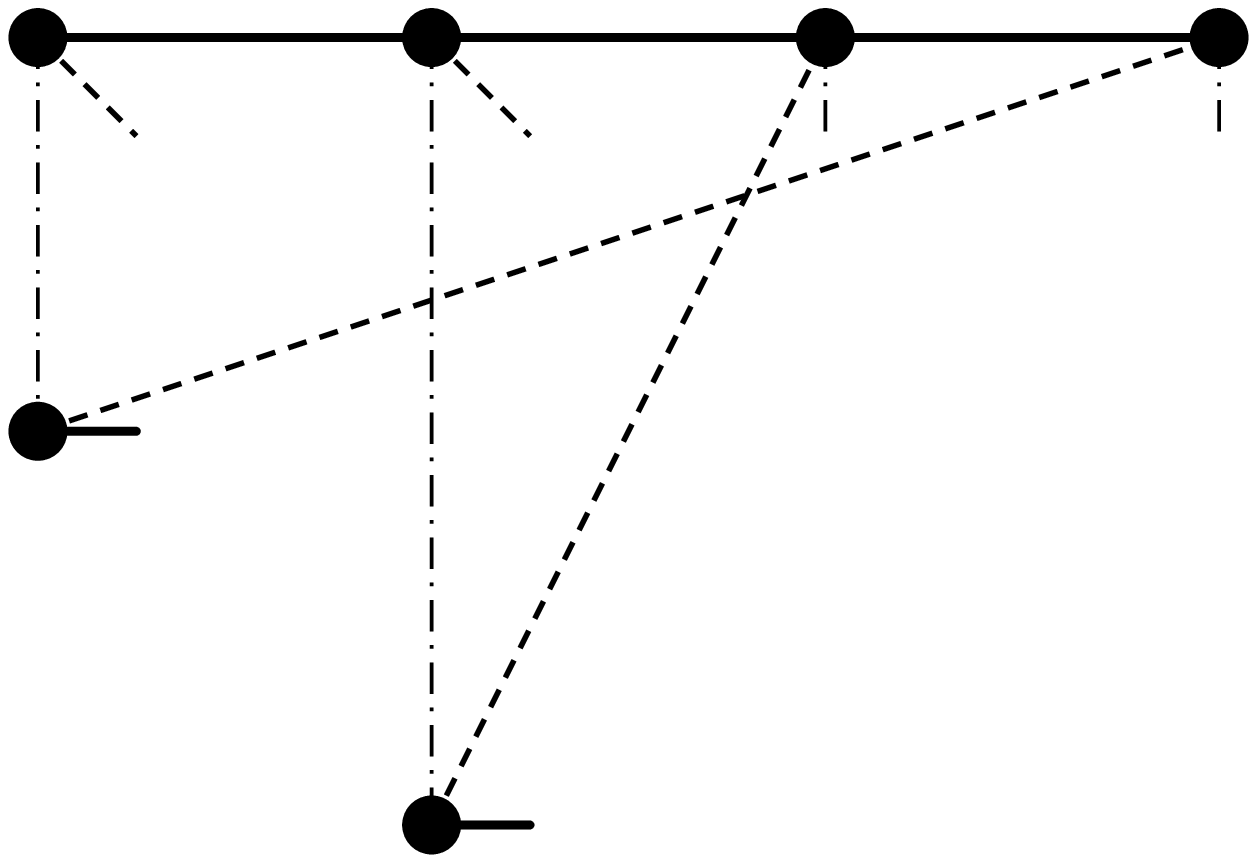}
    \end{tabular} &
    \begin{tabular}{c}
    \includegraphics[width=.11\textwidth]{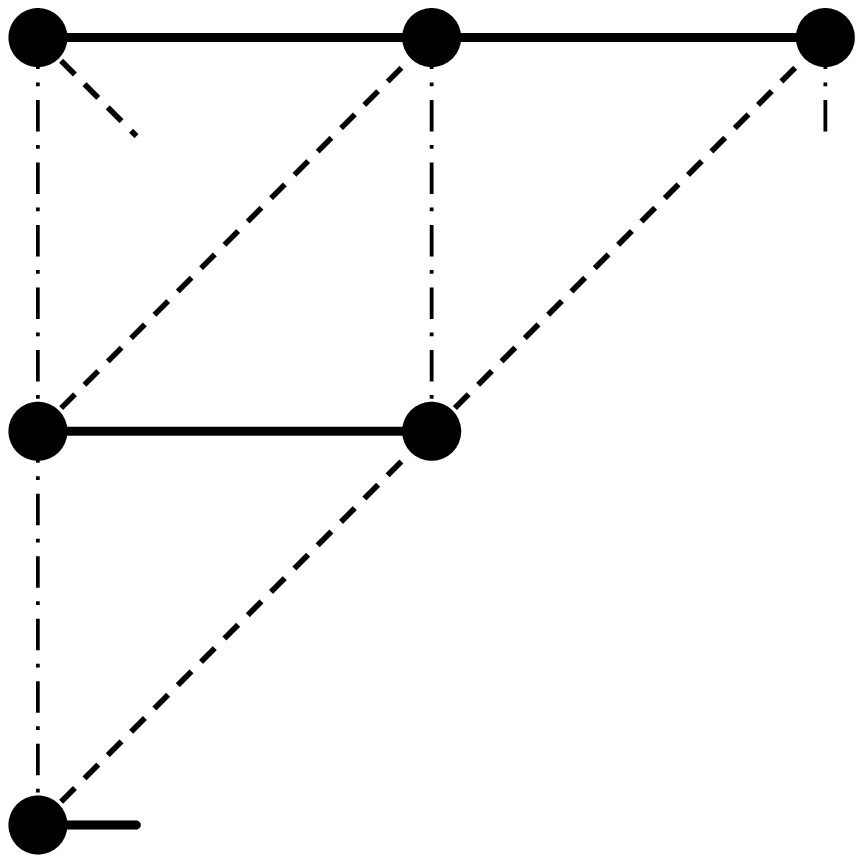}
    \end{tabular} &
    \begin{tabular}{c}
    \\
    \includegraphics[width=.11\textwidth]{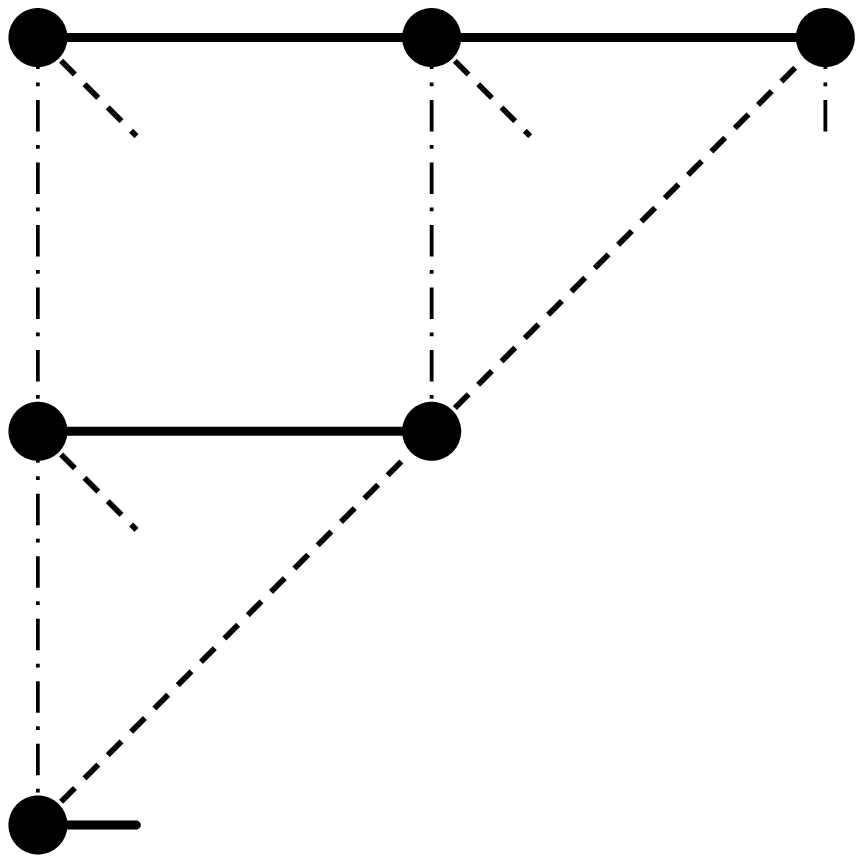}
    \end{tabular} &
    \\
    $\mathcal{S}_{6,7}$ & $\mathcal{S}_{6,8}$ & $\mathcal{S}_{6,9}$ & $\mathcal{S}_{6,10}$ & $\mathcal{S}_{6,11}$ & $\mathcal{S}_{6,12}$\\
    \begin{tabular}{c}
    \\
    \includegraphics[width=.11\textwidth]{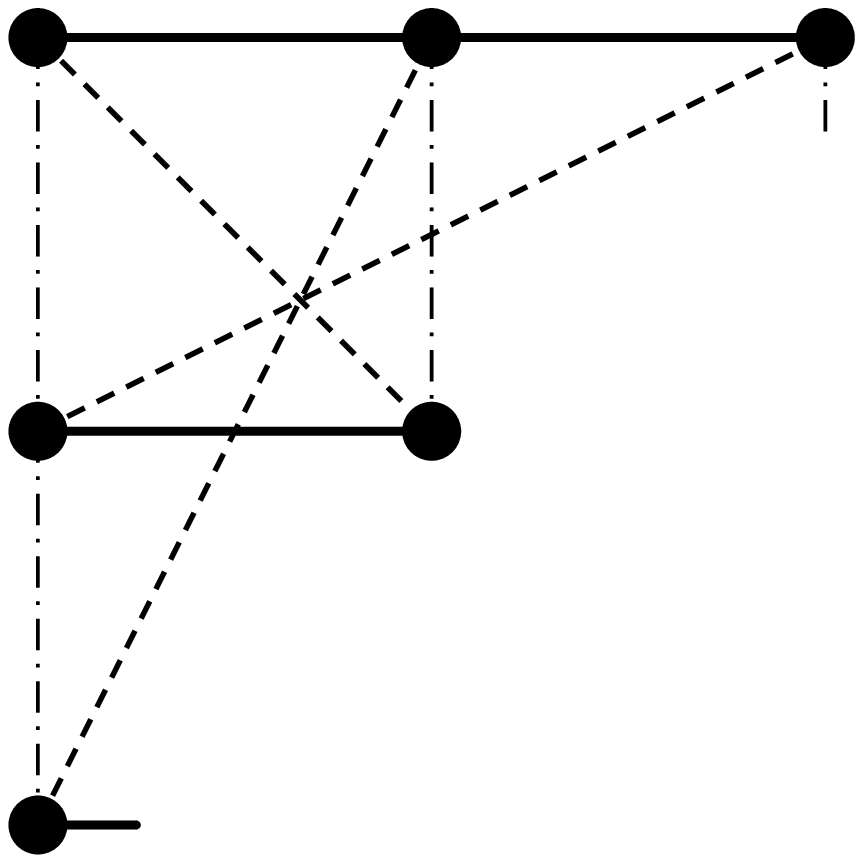}
    \end{tabular} &
    \begin{tabular}{c}
    \\
    \includegraphics[width=.11\textwidth]{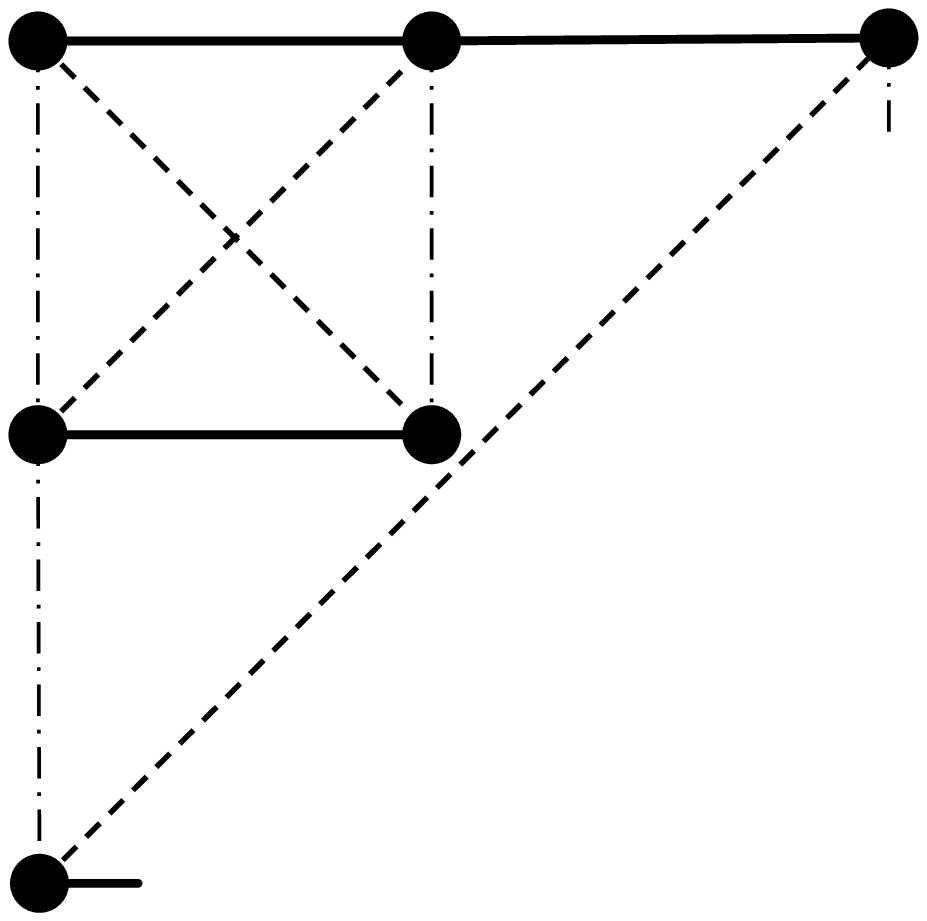}
    \end{tabular} &
    \begin{tabular}{c}
    \includegraphics[width=.11\textwidth]{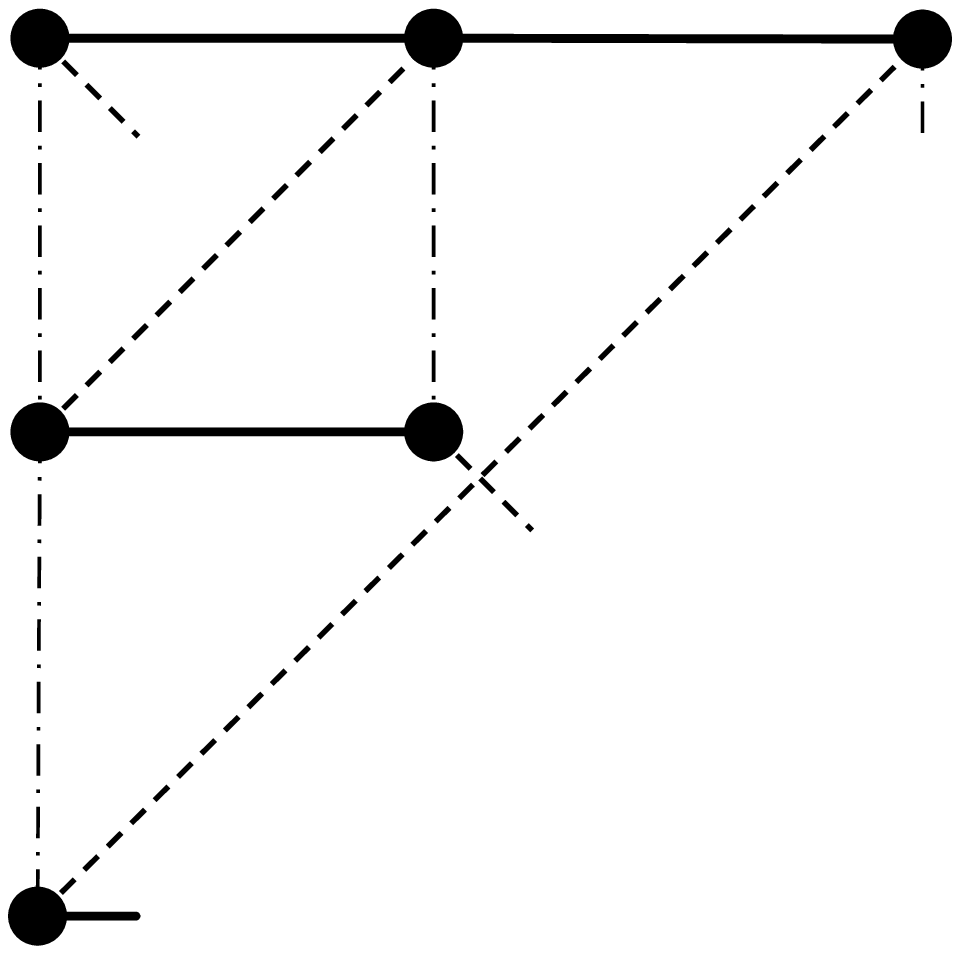}
    \end{tabular} &
    \begin{tabular}{c}
    \\
    \includegraphics[width=.11\textwidth]{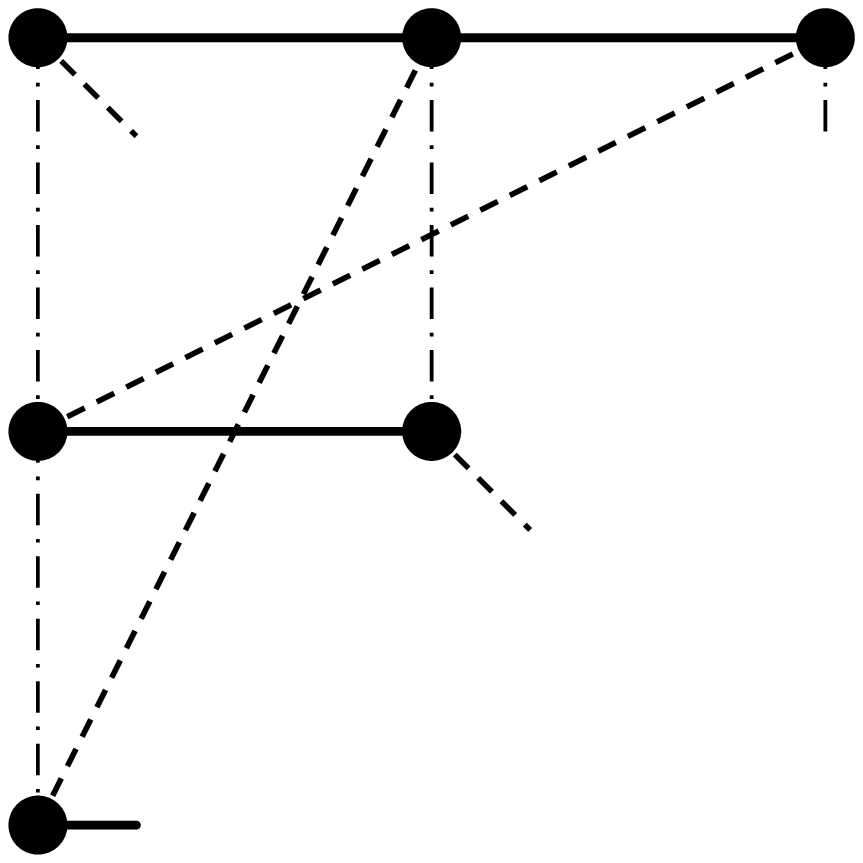}
    \end{tabular} &
    \begin{tabular}{c}
    \\
    \includegraphics[width=.11\textwidth]{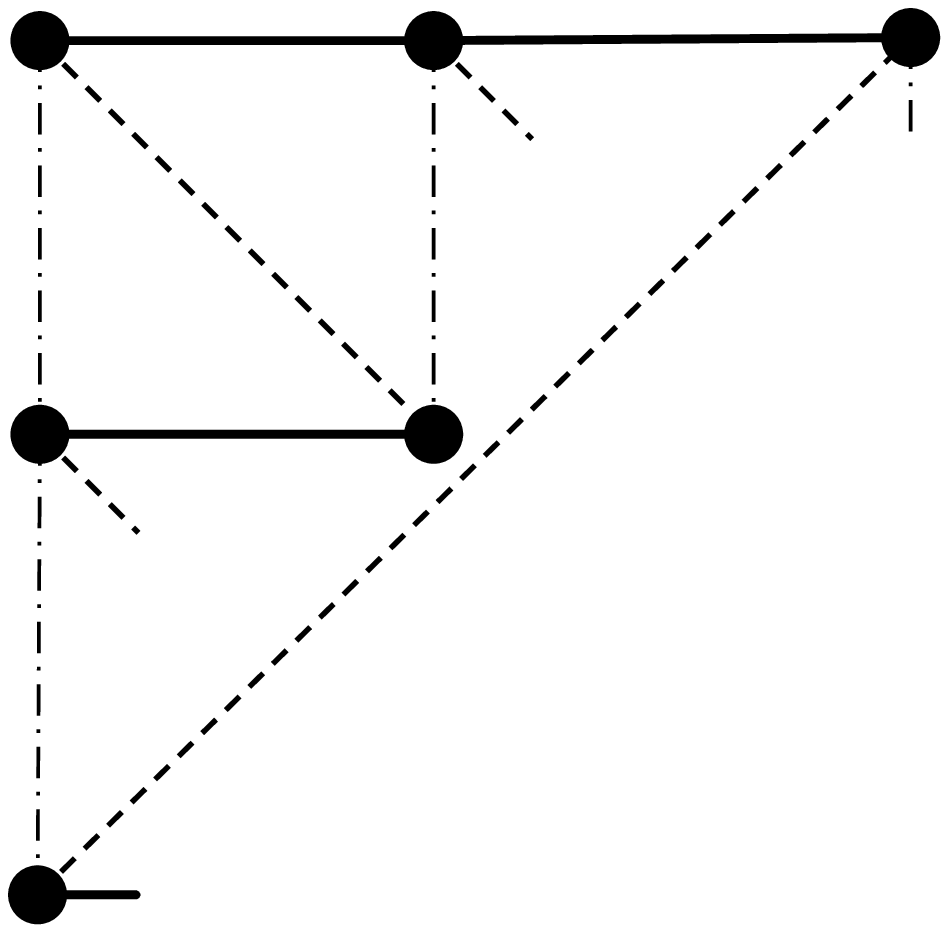}
    \end{tabular} &
    \begin{tabular}{c}
    \\
    \includegraphics[width=.11\textwidth]{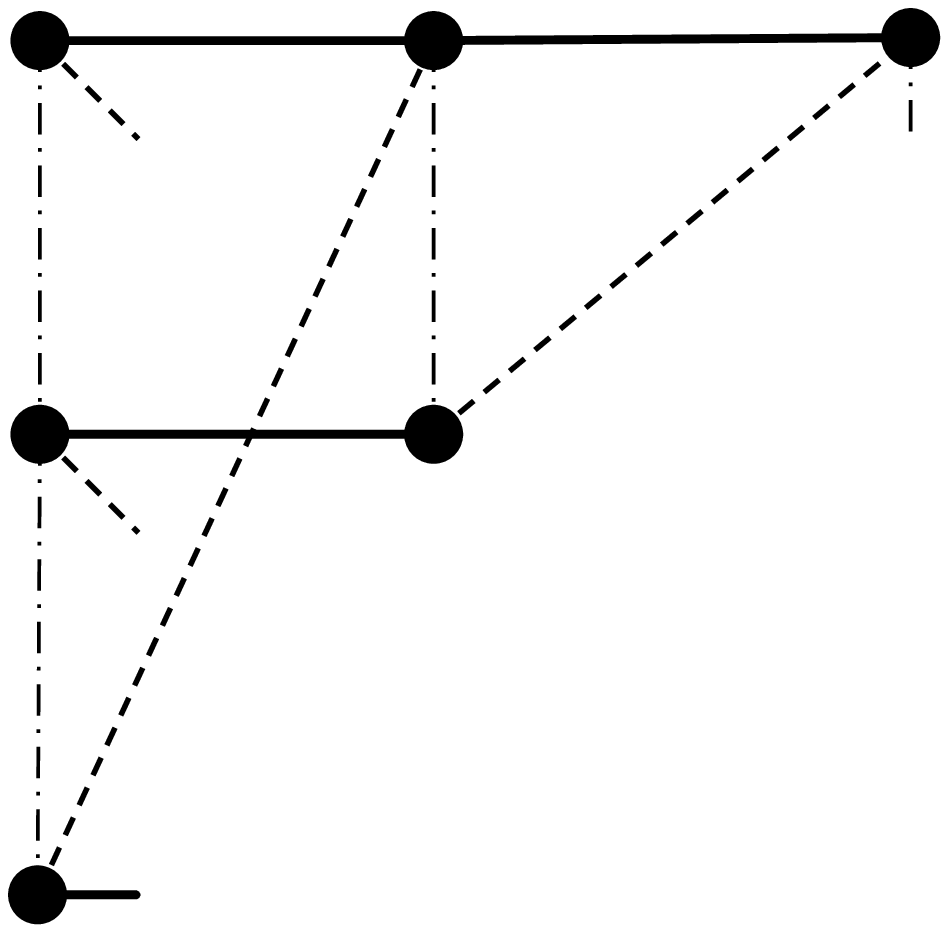}
    \end{tabular}
    \\
    $\mathcal{S}_{6,13}$ & $\mathcal{S}_{6,14}$ & $\mathcal{S}_{6,15}$ & $\mathcal{S}_{6,16}$ & $\mathcal{S}_{6,17}$ & $\mathcal{S}_{6,18}$
    \\
    \begin{tabular}{c}
    \includegraphics[width=.11\textwidth]{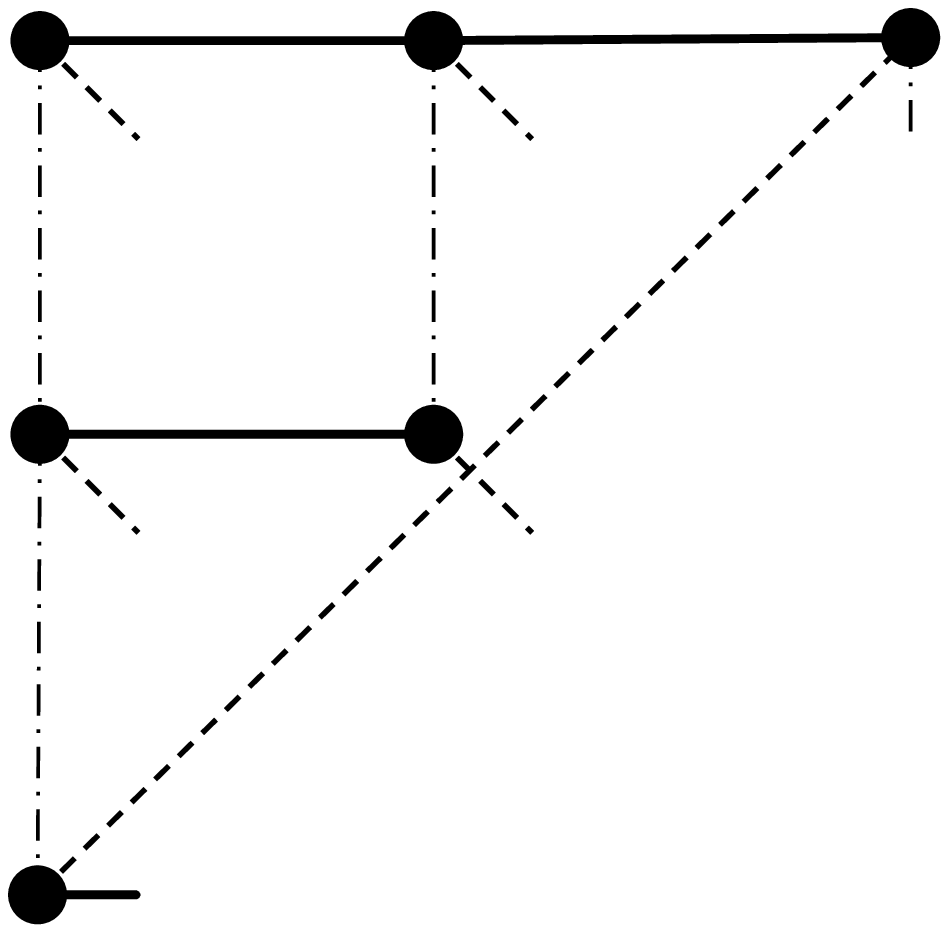}
    \end{tabular} &
    \begin{tabular}{c}
    \\
    \includegraphics[width=.11\textwidth]{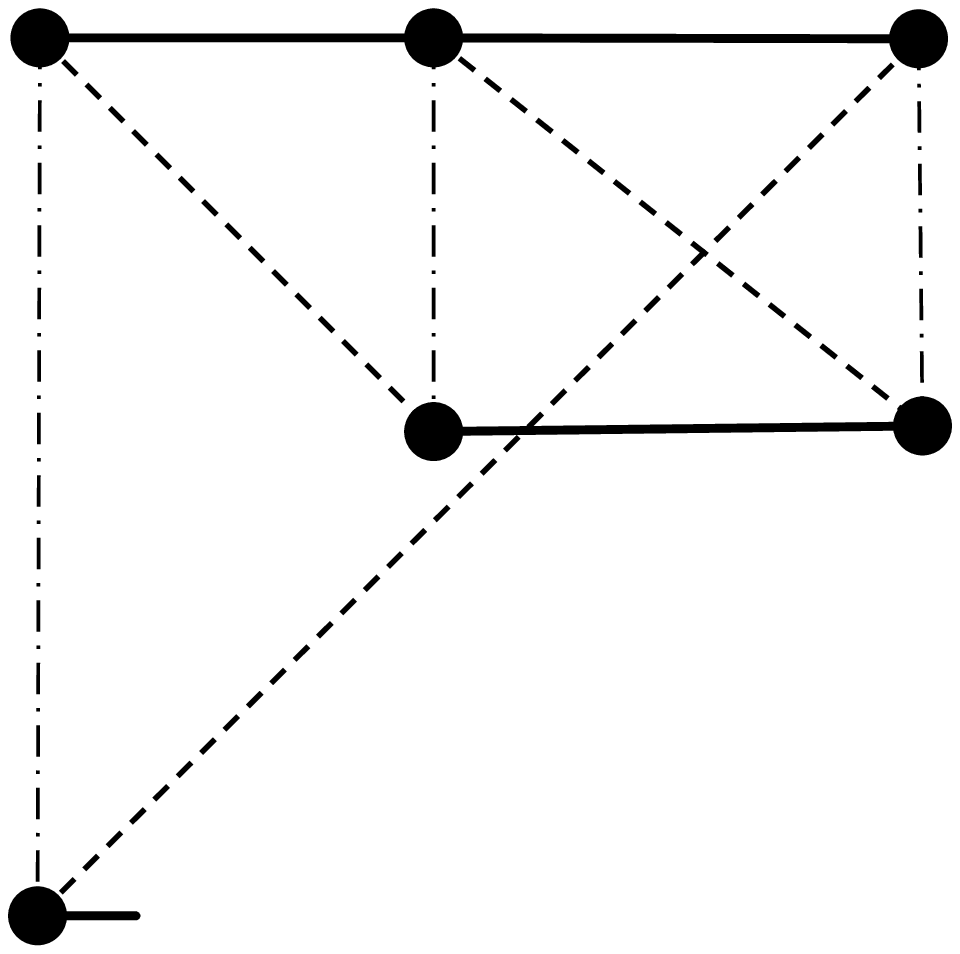}
    \end{tabular} &
    \begin{tabular}{c}
    \\
    \includegraphics[width=.11\textwidth]{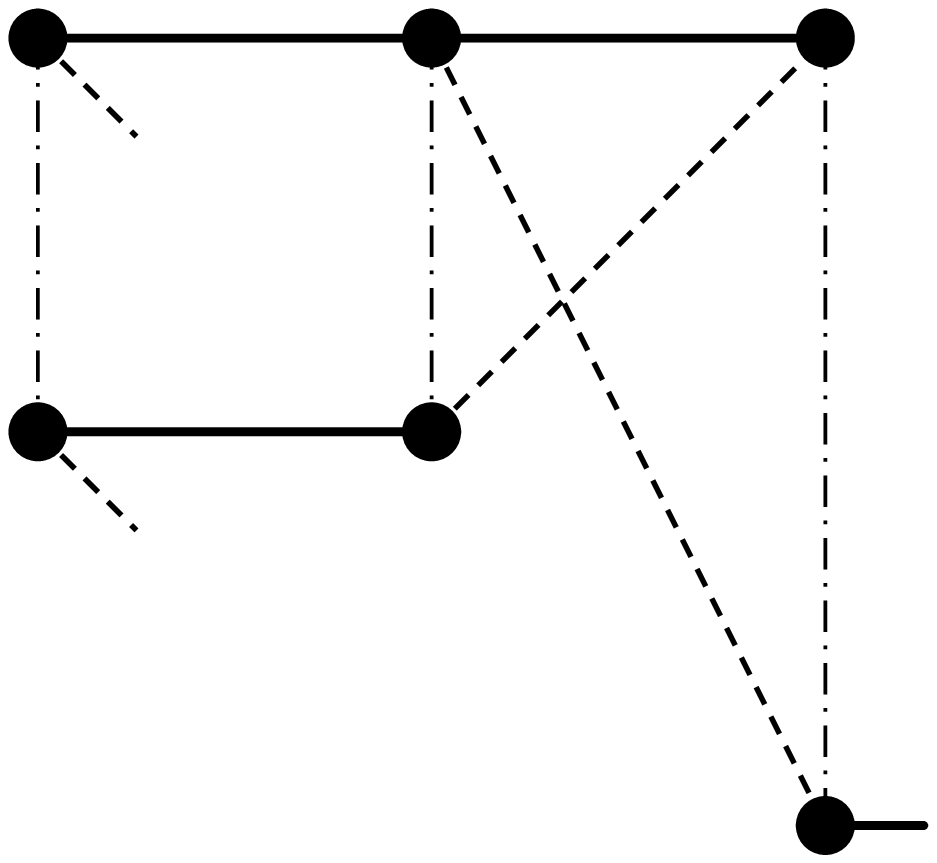}
    \end{tabular} &
    \begin{tabular}{c}
    \\
    \includegraphics[width=.11\textwidth]{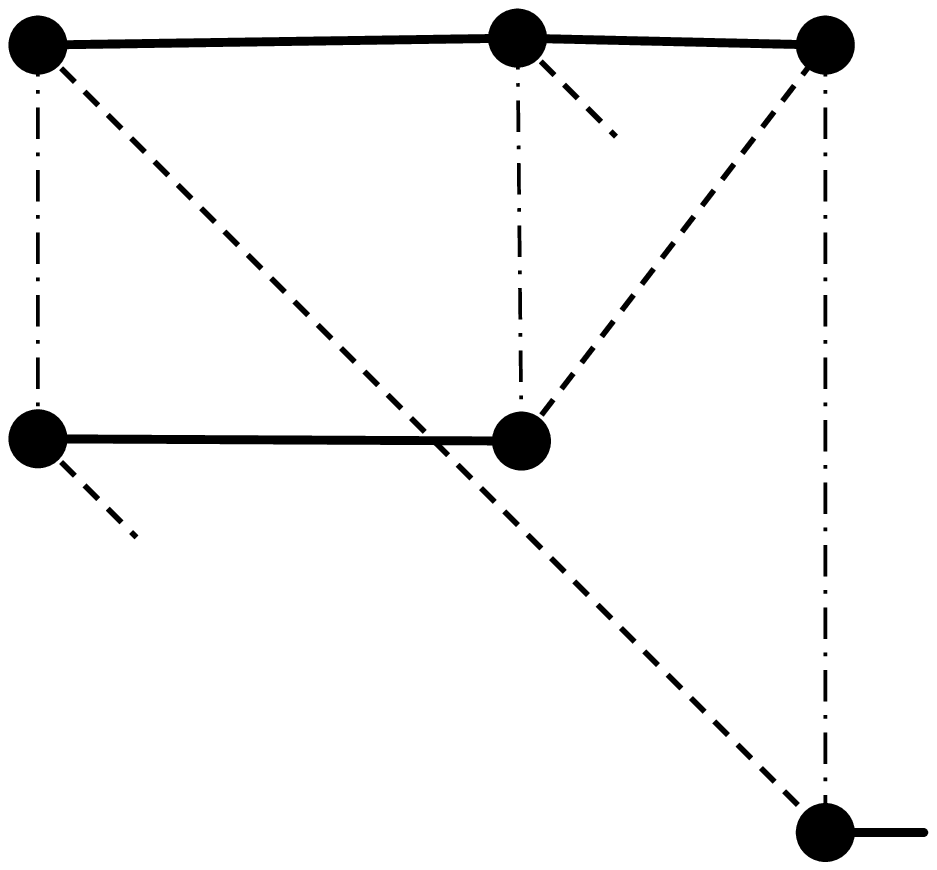}
    \end{tabular} &
    \begin{tabular}{c}
    \includegraphics[width=.11\textwidth]{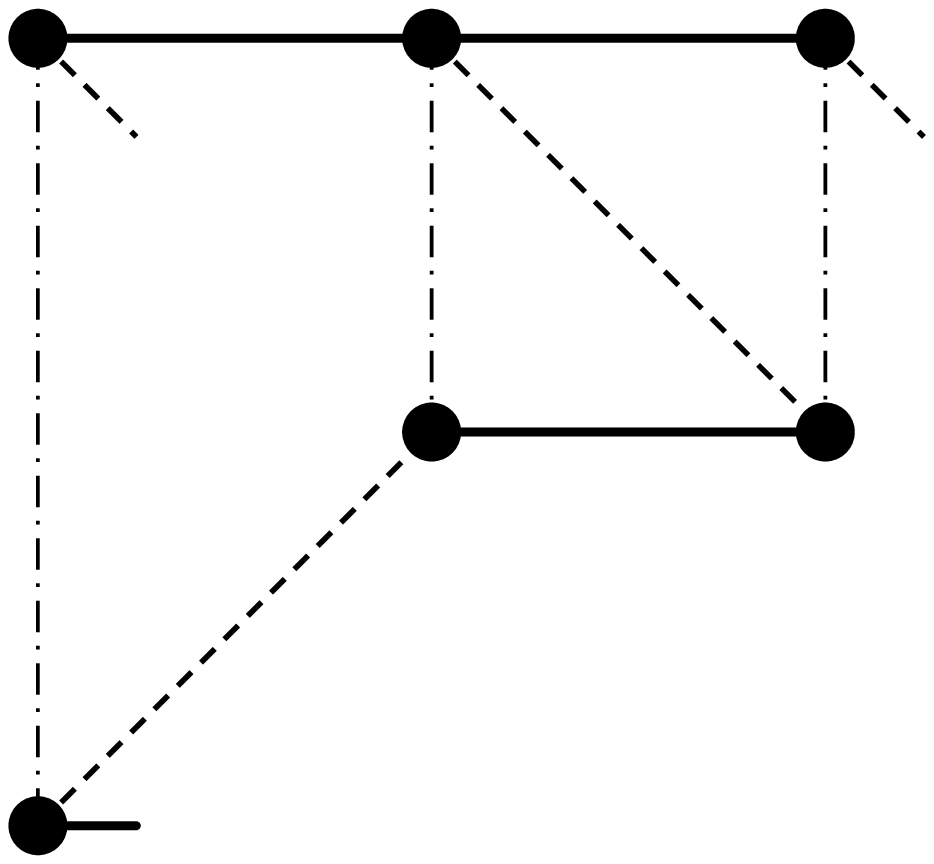}
    \end{tabular} &
    \begin{tabular}{c}
    \\
    \includegraphics[width=.11\textwidth]{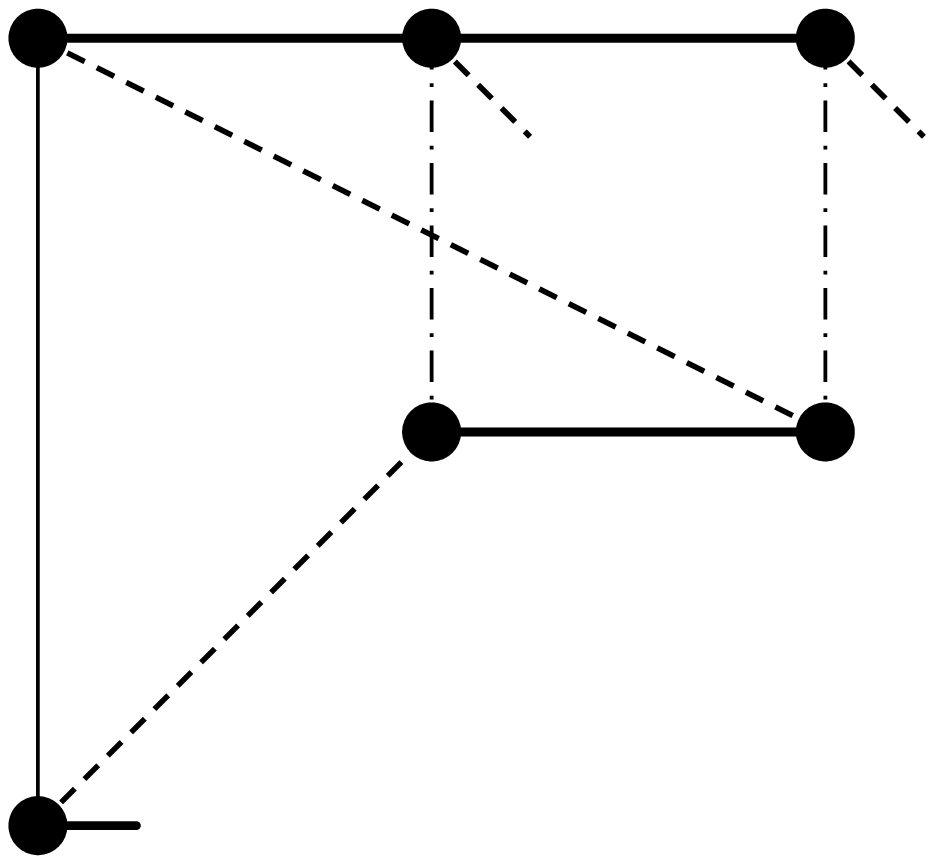}
    \end{tabular} &
    \\
    $\mathcal{S}_{6,19}$ & $\mathcal{S}_{6,20}$ & $\mathcal{S}_{6,21}$ & $\mathcal{S}_{6,22}$ & $\mathcal{S}_{6,23}$ & $\mathcal{S}_{6,24}$\\
    \begin{tabular}{c}
    \\
    \includegraphics[width=.11\textwidth]{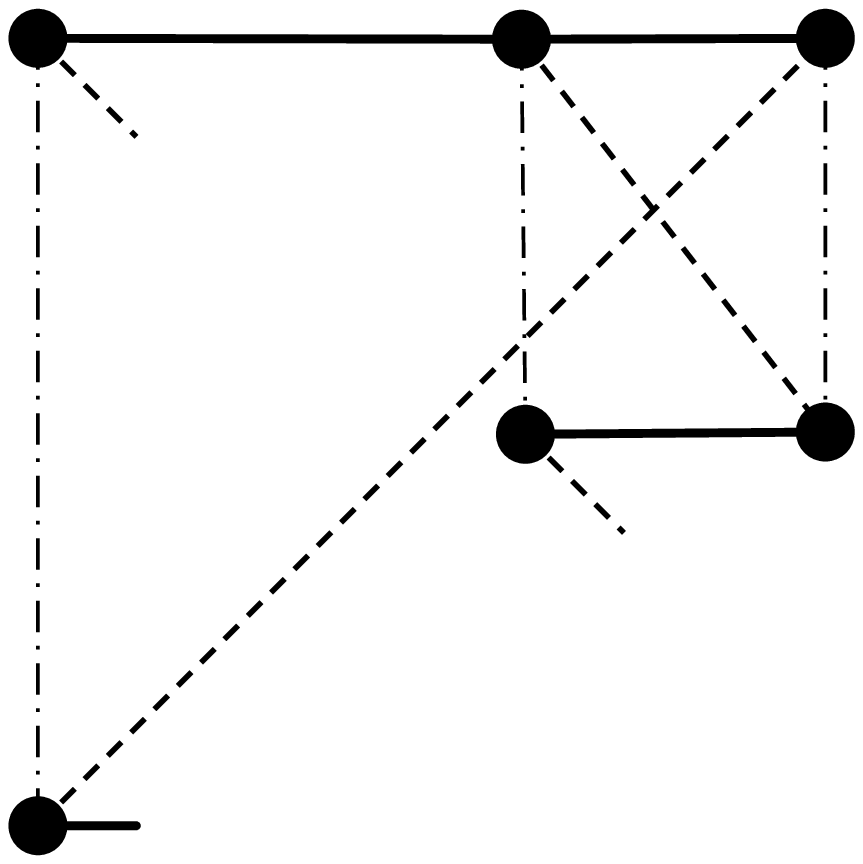}
    \end{tabular} &
    \begin{tabular}{c}
    \\
    \includegraphics[width=.11\textwidth]{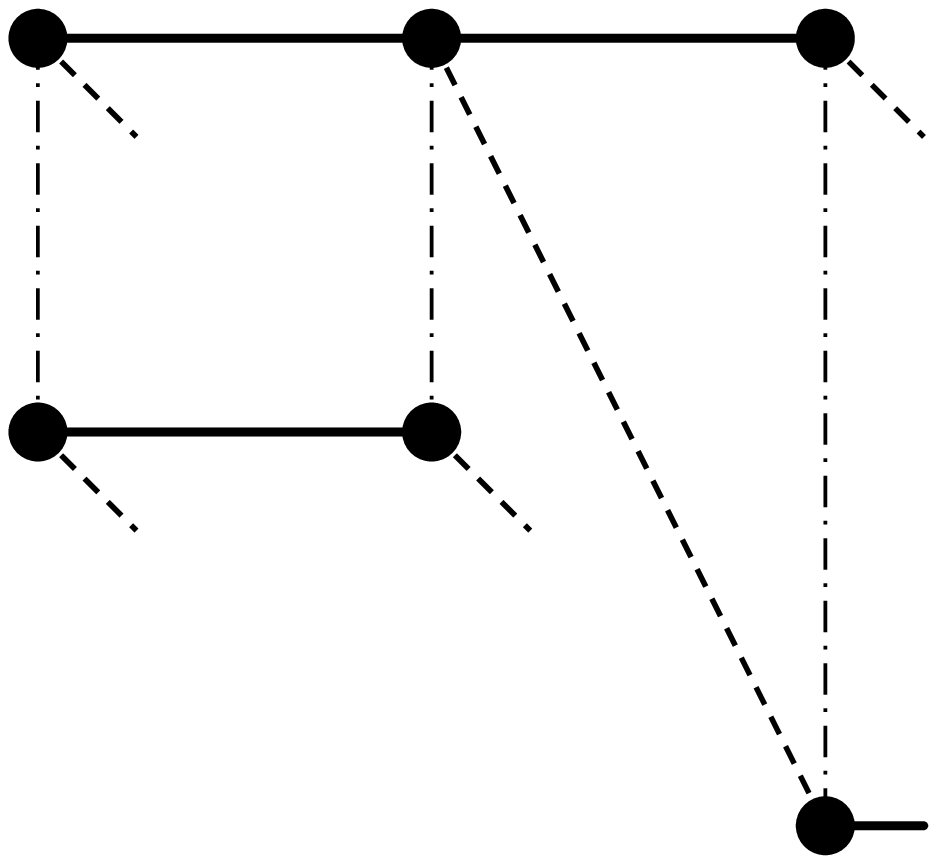}
    \end{tabular} &
    \begin{tabular}{c}
    \includegraphics[width=.11\textwidth]{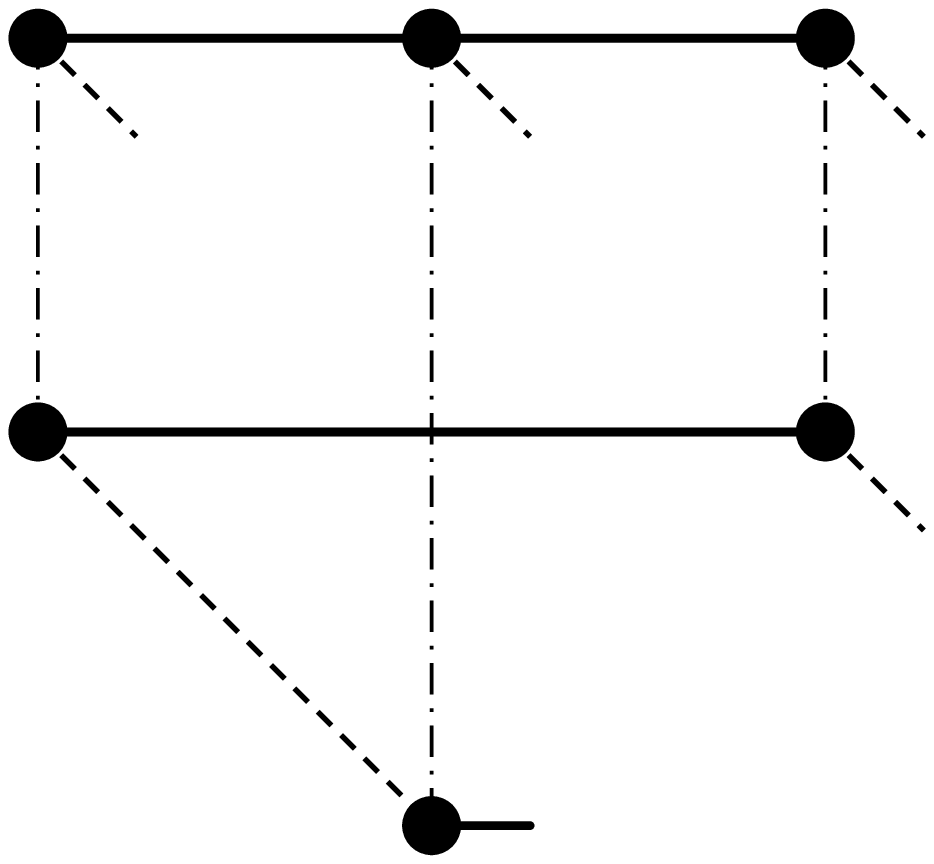}
    \end{tabular} &
    \begin{tabular}{c}
    \\
    \includegraphics[width=.11\textwidth]{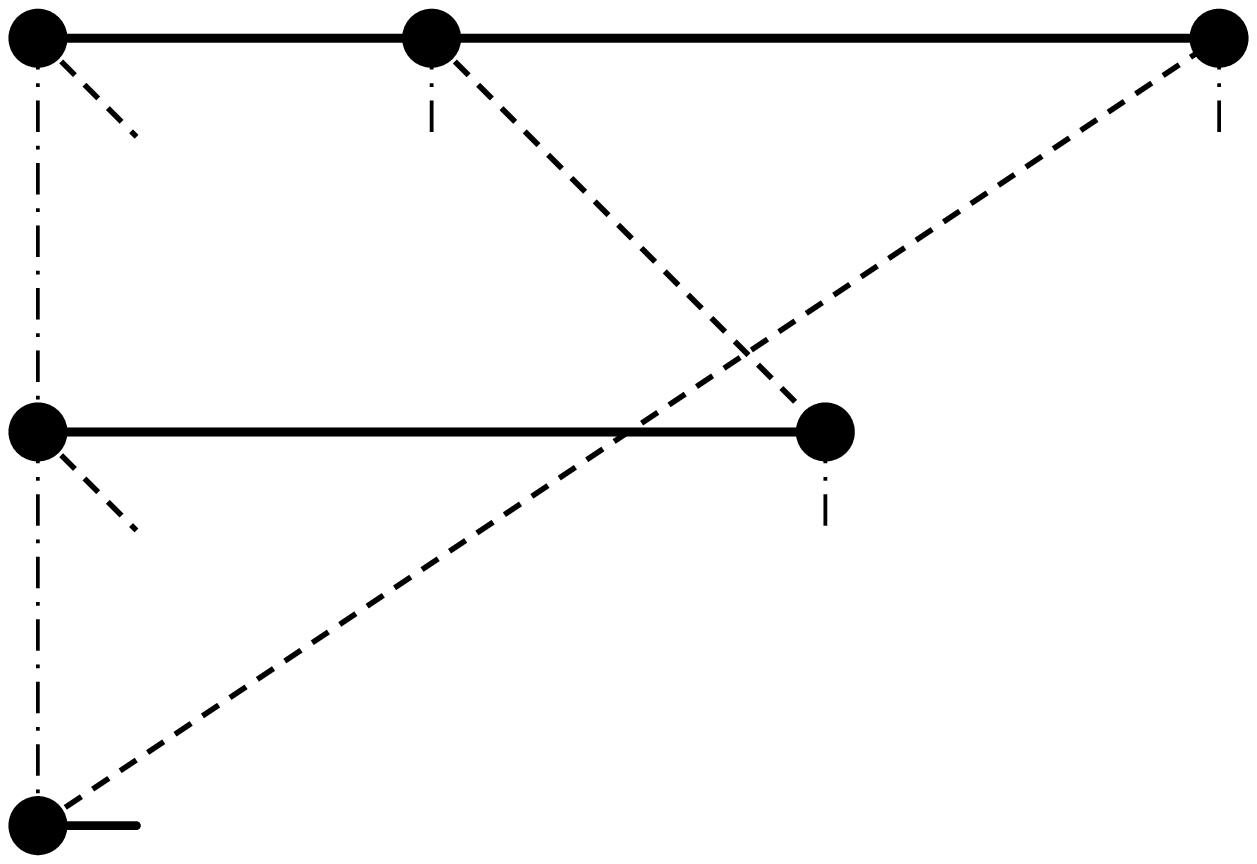}
    \end{tabular} &
    \begin{tabular}{c}
    \\
    \includegraphics[width=.11\textwidth]{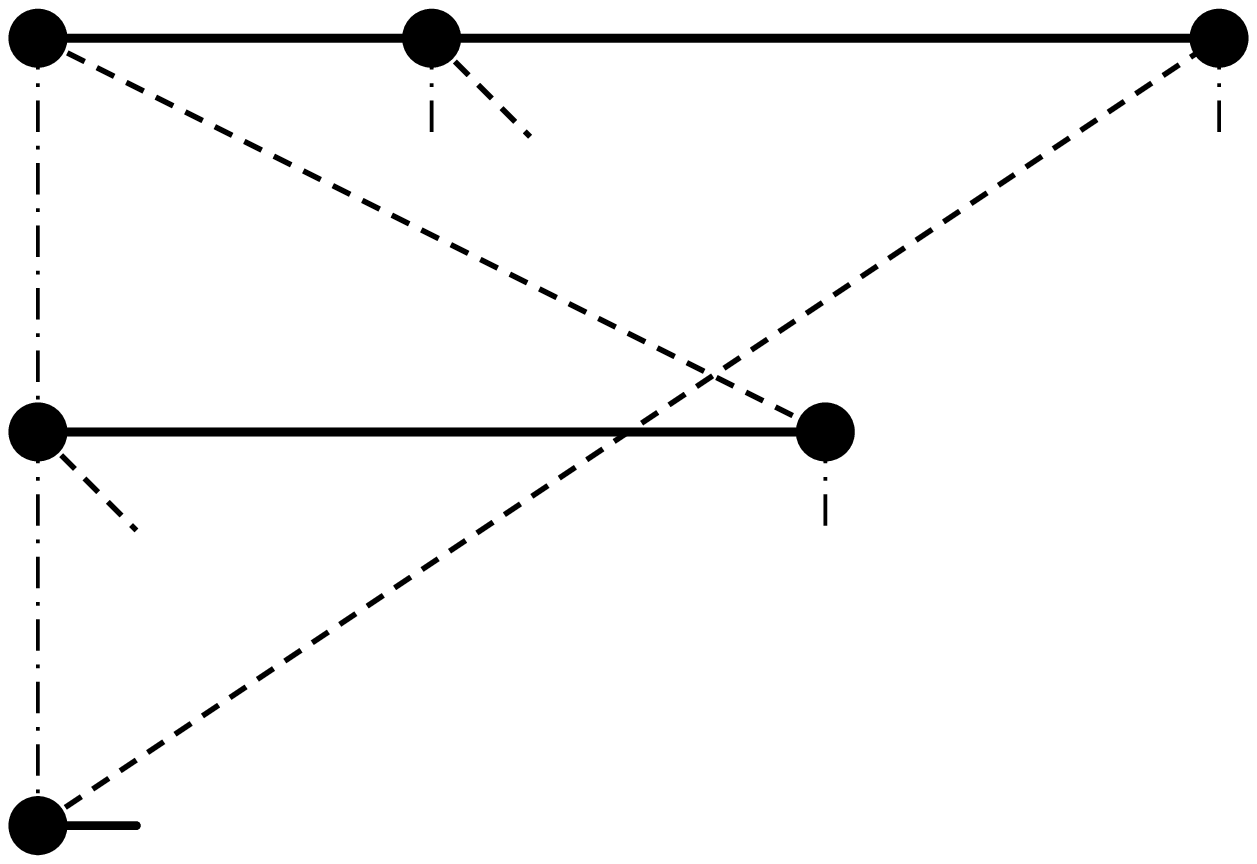}
    \end{tabular} &
    \begin{tabular}{c}
    \\
    \includegraphics[width=.11\textwidth]{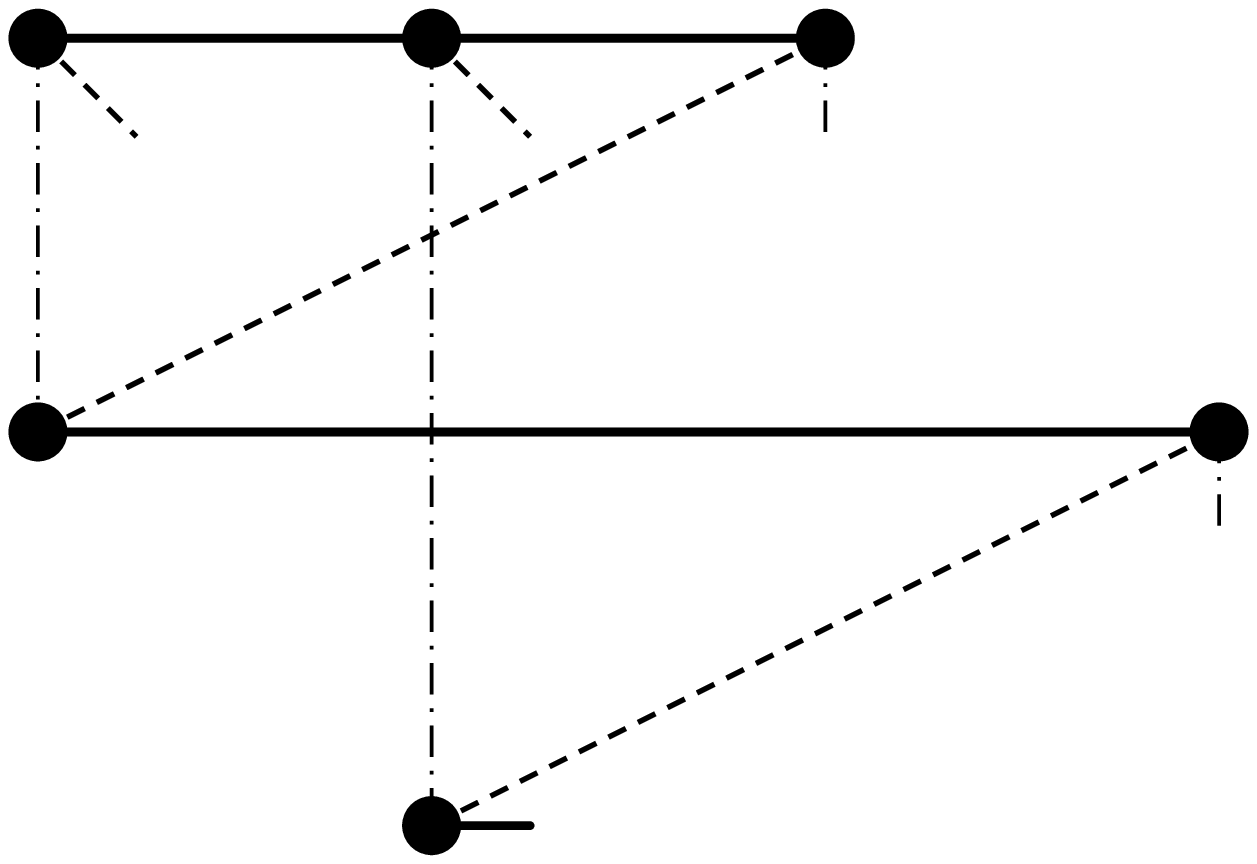}
    \end{tabular}
    \\
    $\mathcal{S}_{6,25}$ & $\mathcal{S}_{6,26}$ & $\mathcal{S}_{6,27}$ & $\mathcal{S}_{6,28}$ & $\mathcal{S}_{6,29}$ & $\mathcal{S}_{6,30}$\\
    \begin{tabular}{c}
    \\
    \includegraphics[width=.11\textwidth]{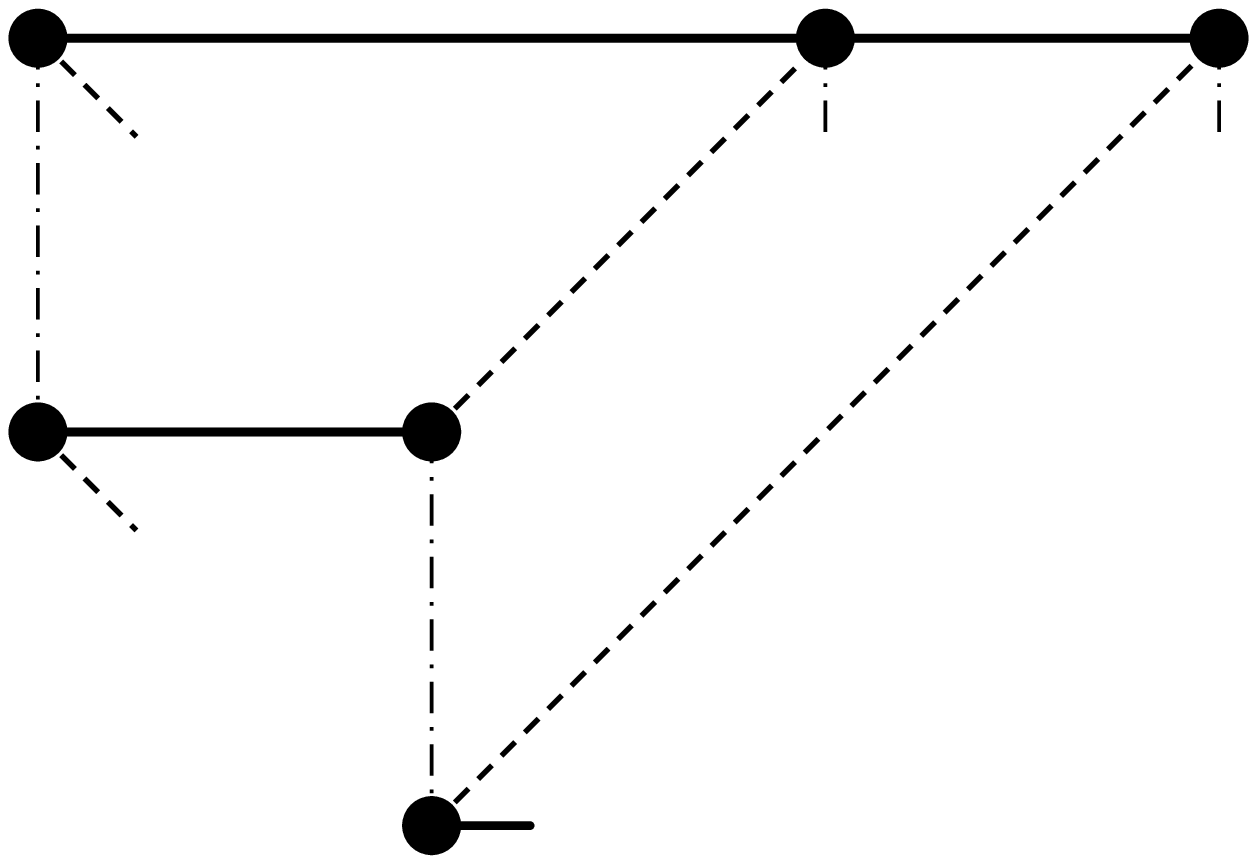}
    \end{tabular} &
    \begin{tabular}{c}
    \\
    \includegraphics[width=.11\textwidth]{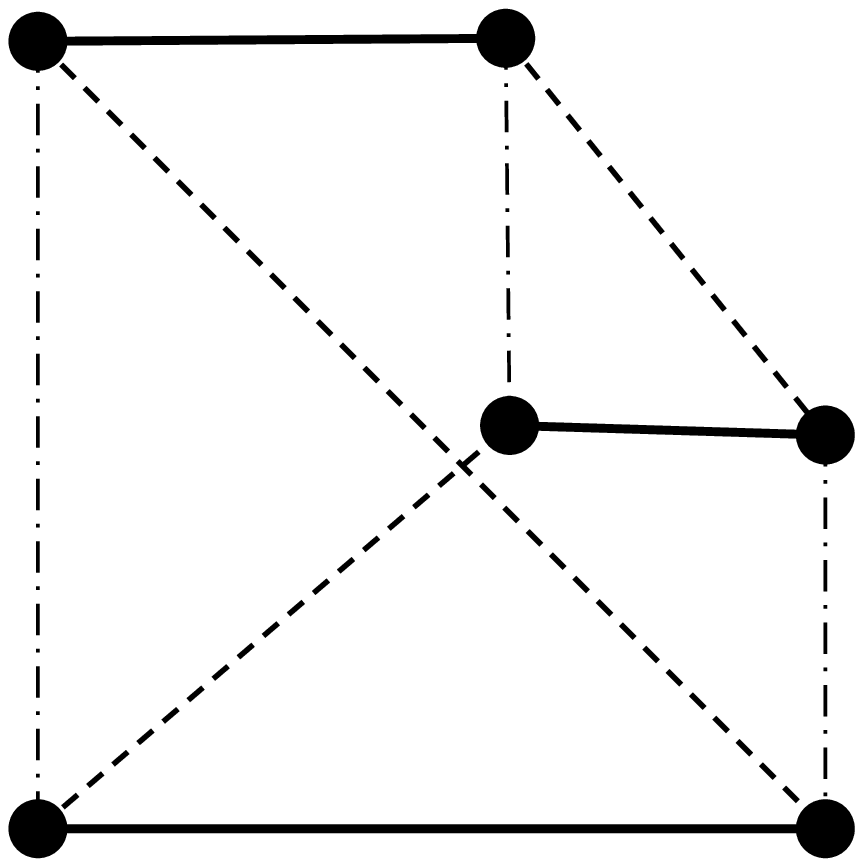}
    \end{tabular} &
    \begin{tabular}{c}
    \\
    \includegraphics[width=.11\textwidth]{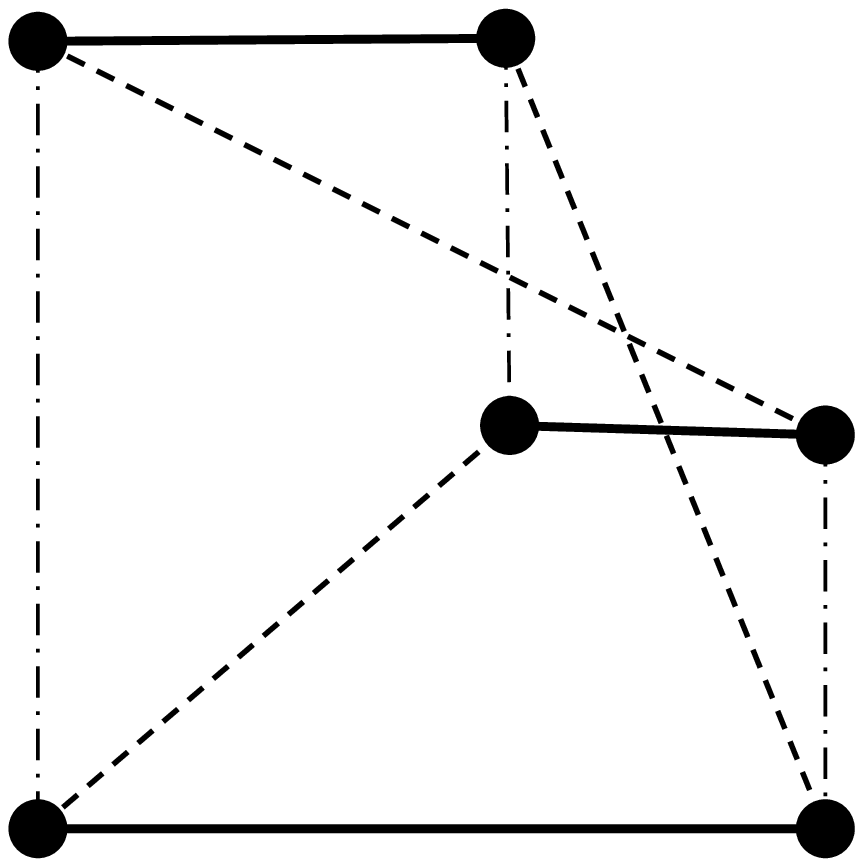}
    \end{tabular} &
    \begin{tabular}{c}
    \\
    \includegraphics[width=.11\textwidth]{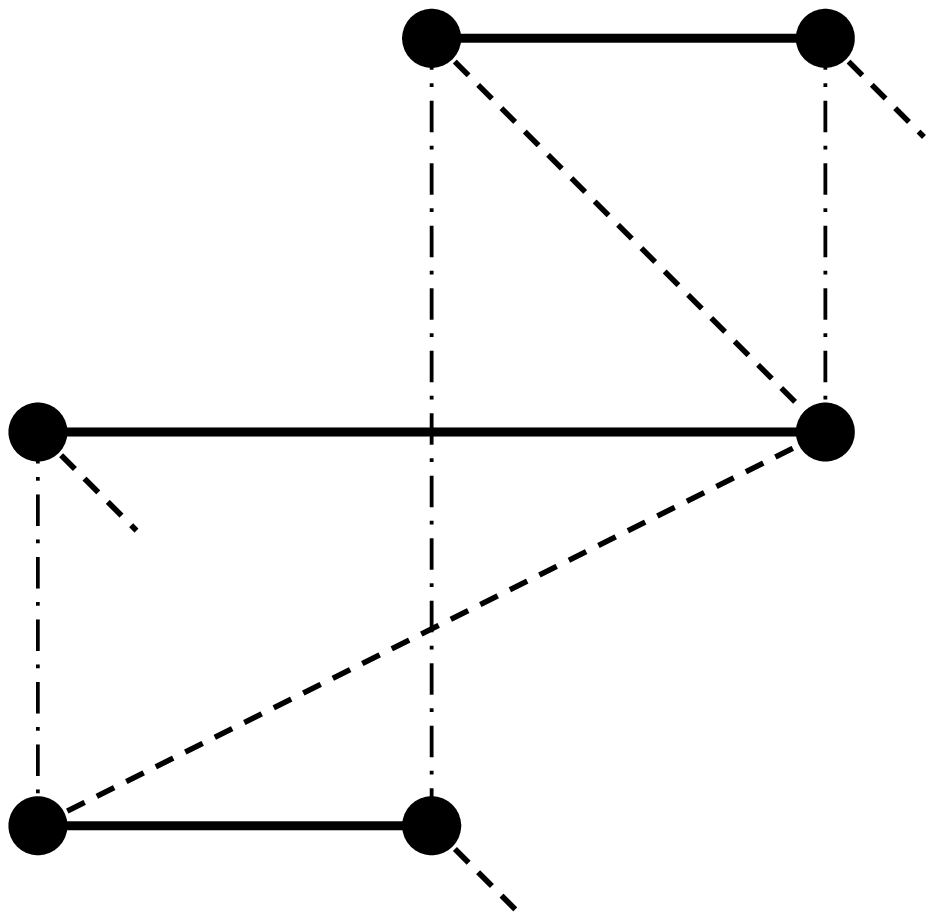}
    \end{tabular} &
    \begin{tabular}{c}
    \\
    \includegraphics[width=.11\textwidth]{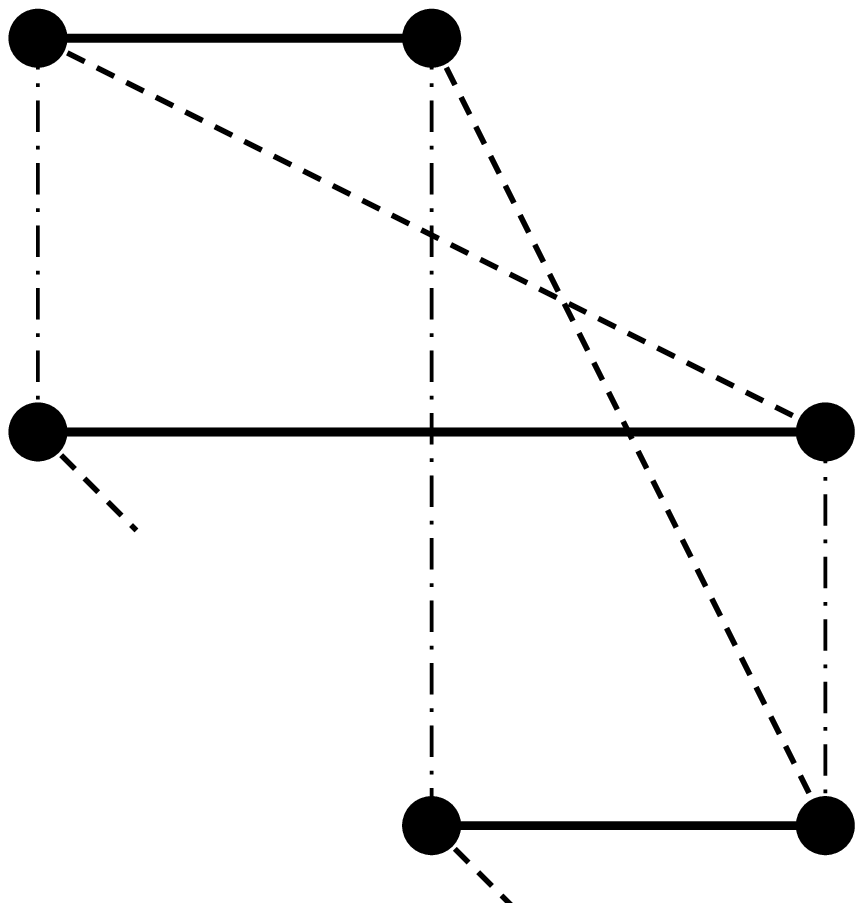}
    \end{tabular} &
    \begin{tabular}{c}
    \\
    \includegraphics[width=.11\textwidth]{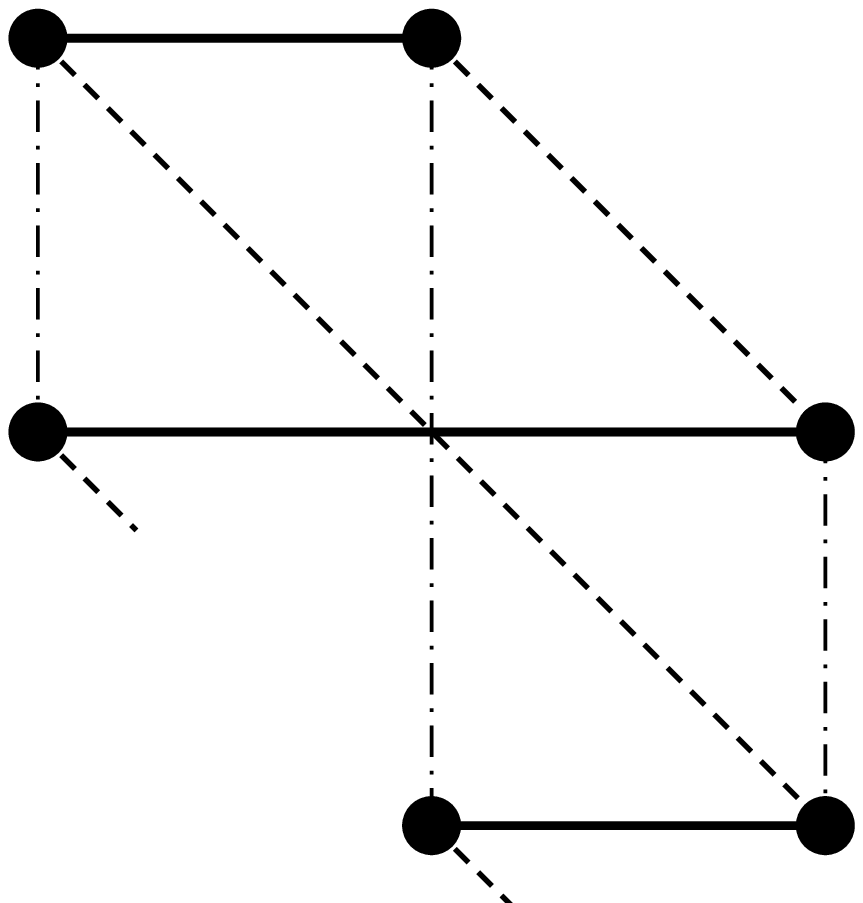}
    \end{tabular}
    \\
    $\mathcal{S}_{6,31}$ & $\mathcal{S}_{6,32}$ & $\mathcal{S}_{6,33}$ & $\mathcal{S}_{6,34}$ & $\mathcal{S}_{6,35}$ & $\mathcal{S}_{6,36}$\\
    \begin{tabular}{c}
    \\
    \includegraphics[width=.11\textwidth]{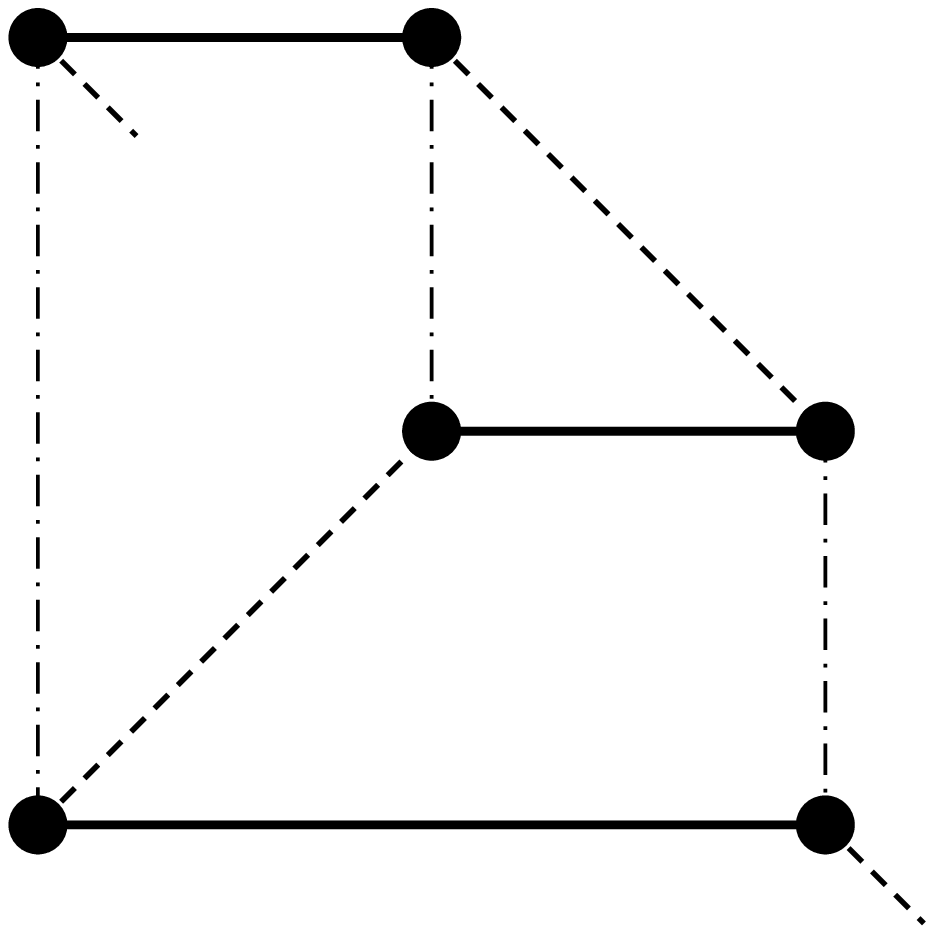}
    \end{tabular} &
    \begin{tabular}{c}
    \\
    \includegraphics[width=.11\textwidth]{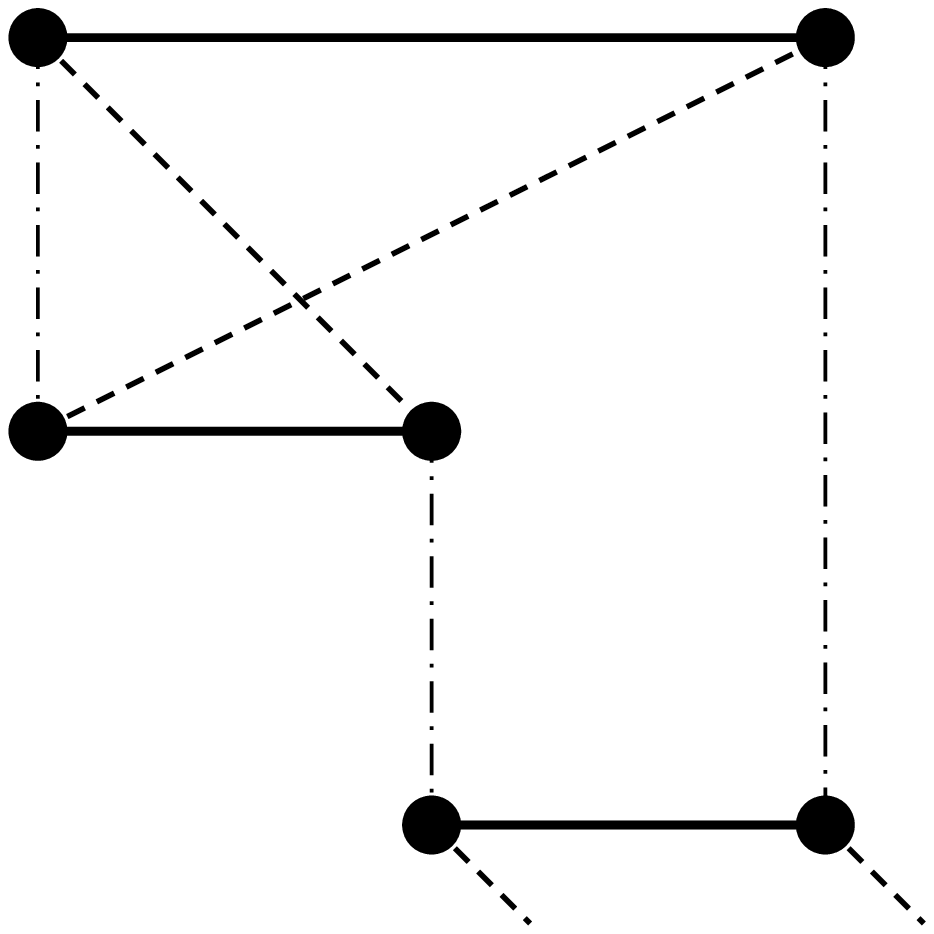}
    \end{tabular} &
    \begin{tabular}{c}
    \\
    \includegraphics[width=.11\textwidth]{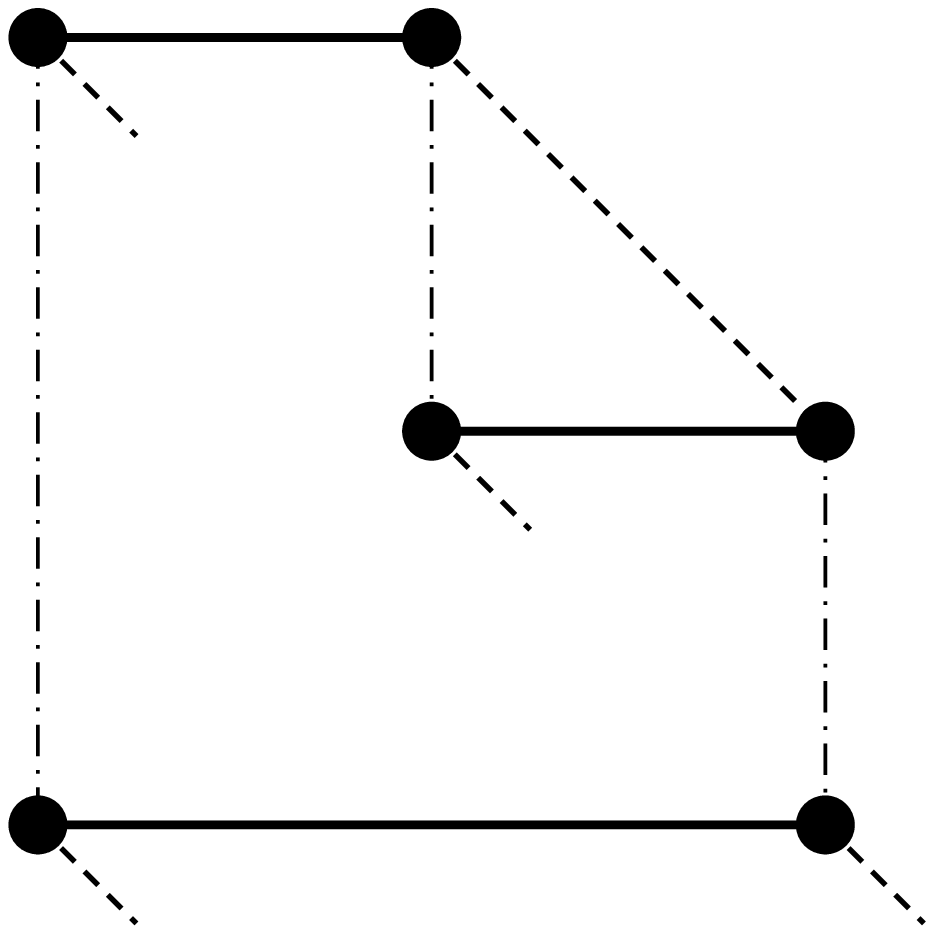}
    \end{tabular} &
    \begin{tabular}{c}
    \\
    \includegraphics[width=.11\textwidth]{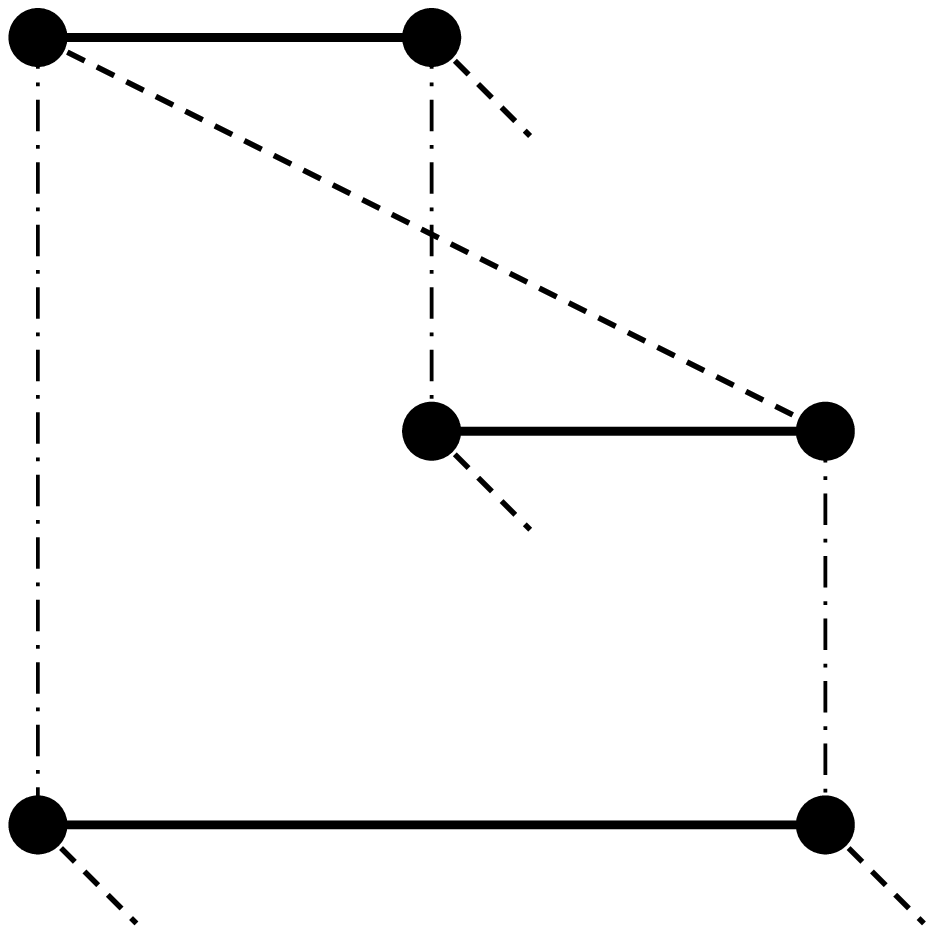}
    \end{tabular} &
    \begin{tabular}{c}
    \\
    \includegraphics[width=.11\textwidth]{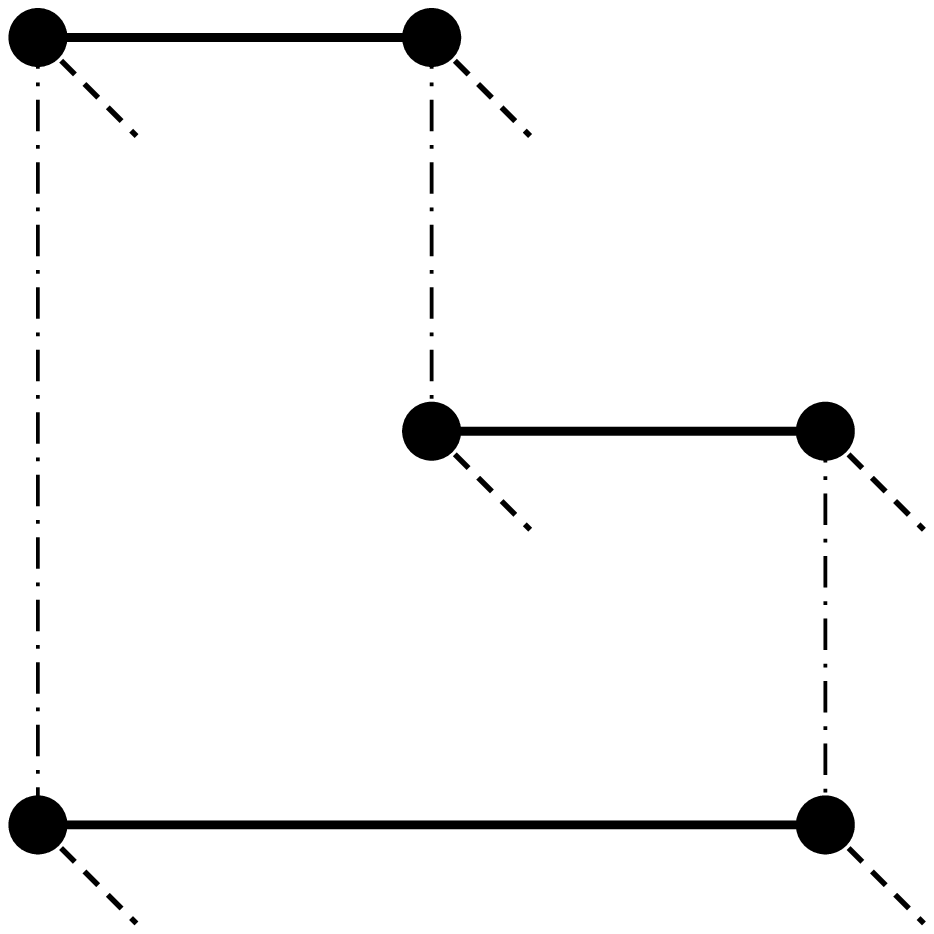}
    \end{tabular} &
    \begin{tabular}{c}
    \\
    \includegraphics[width=.11\textwidth]{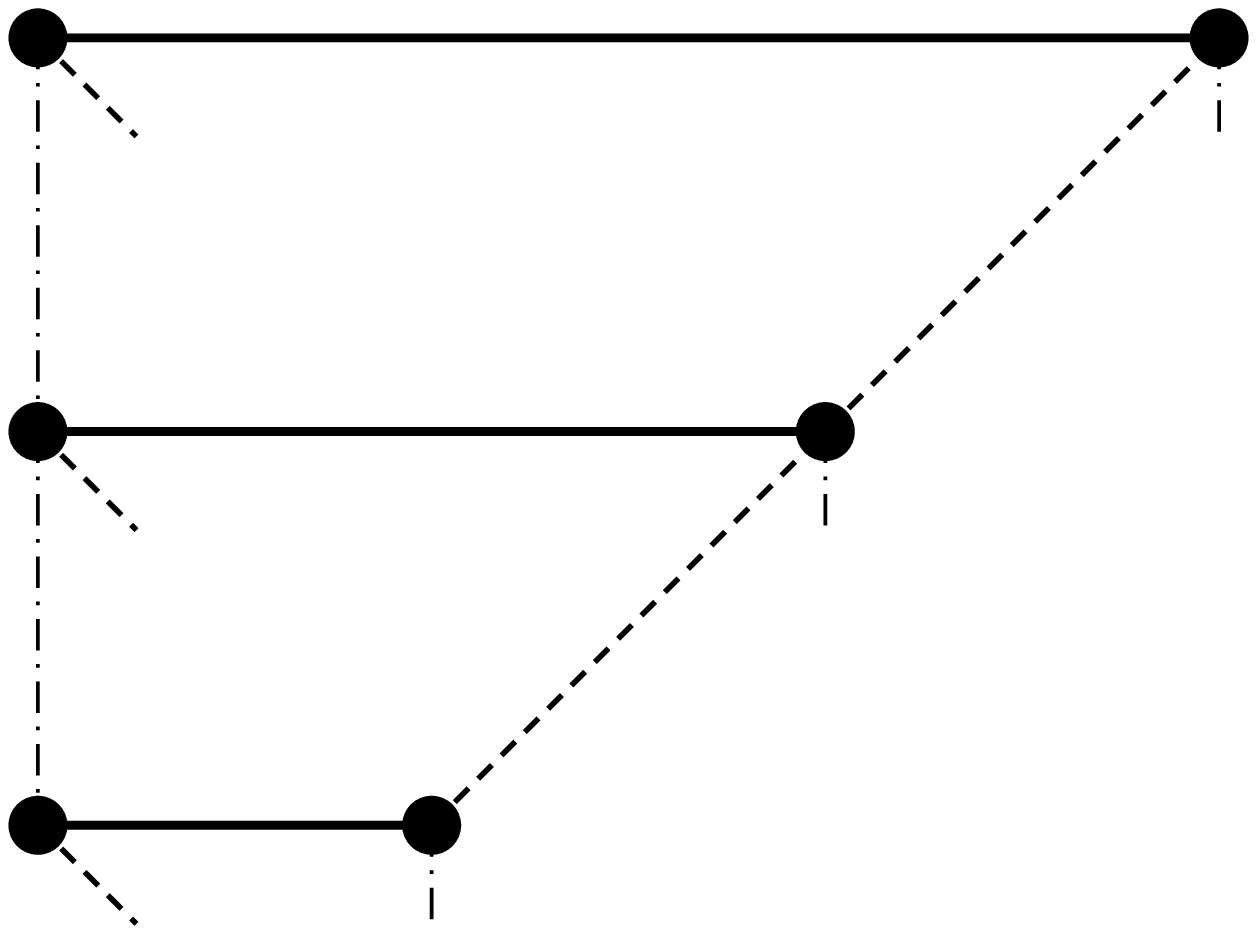}
    \end{tabular}
    \\
    $\mathcal{S}_{6,37}$ & $\mathcal{S}_{6,38}$ & $\mathcal{S}_{6,39}$ & $\mathcal{S}_{6,40}$ & $\mathcal{S}_{6,41}$ & $\mathcal{S}_{6,42}$ \\
    \begin{tabular}{c}
    \\
    \includegraphics[width=.11\textwidth]{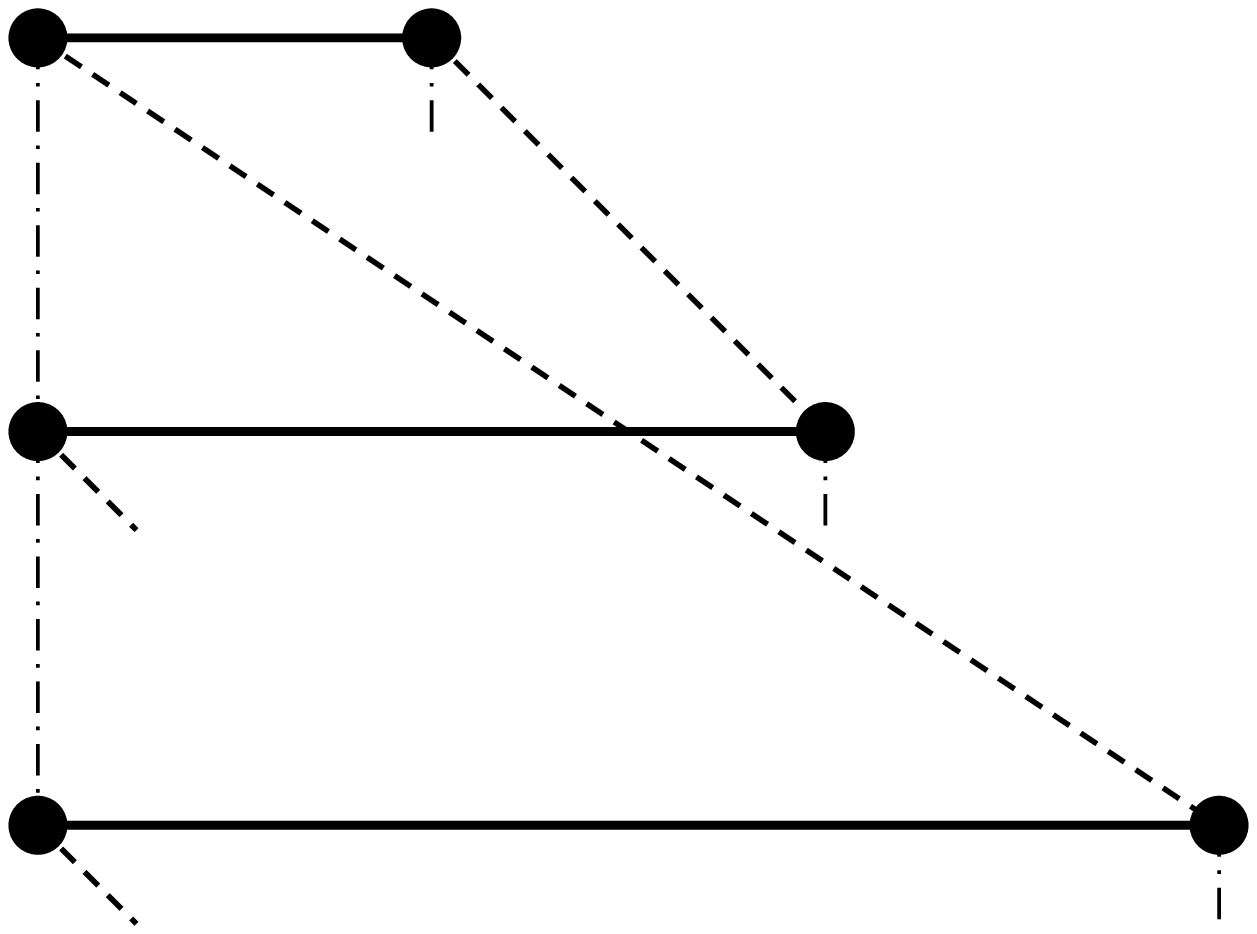}
    \end{tabular}
    &
    \begin{tabular}{c}
    \\
    \includegraphics[width=.11\textwidth]{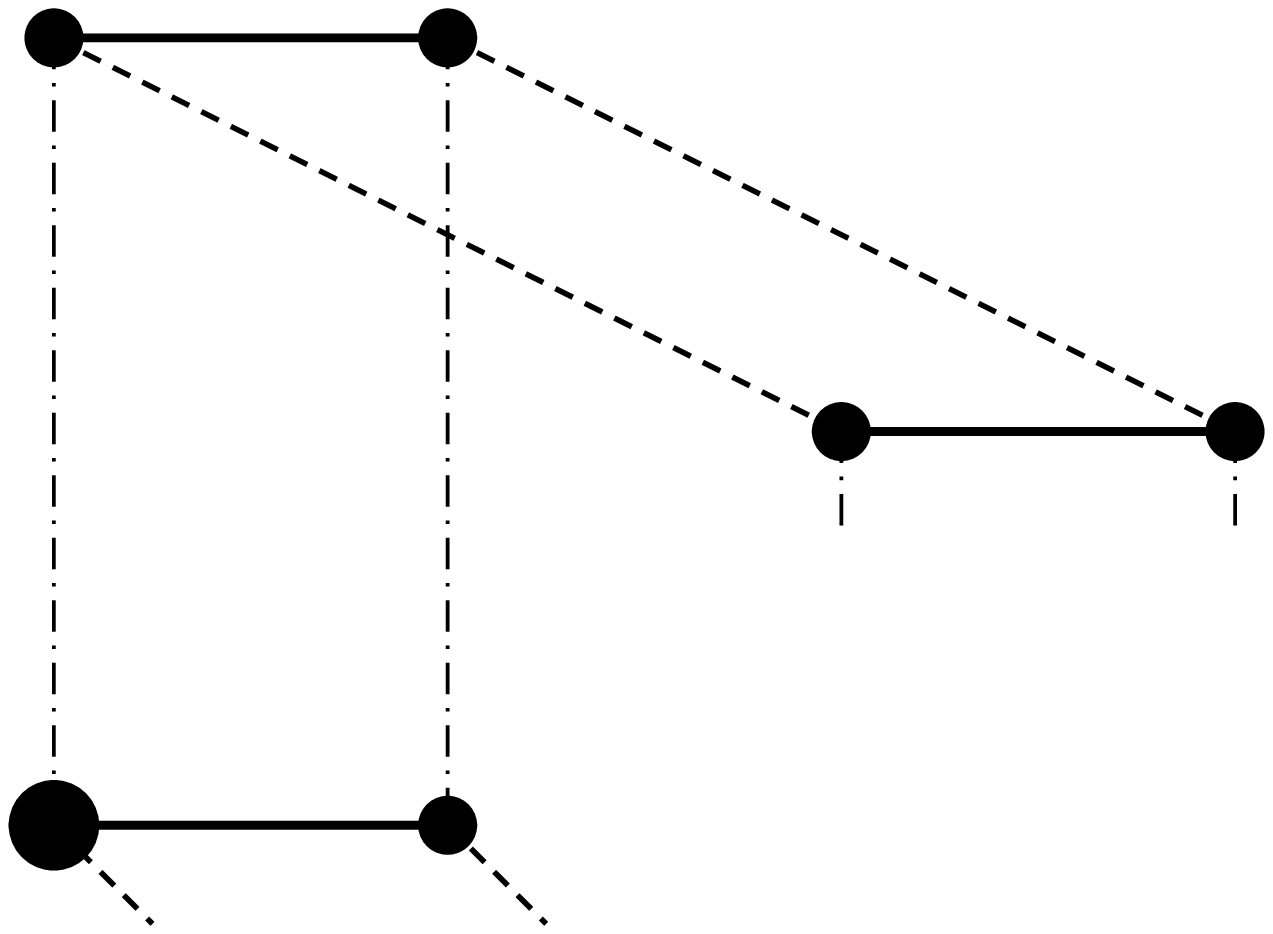}
    \end{tabular} &
    \begin{tabular}{c}
    \\
    \includegraphics[width=.11\textwidth]{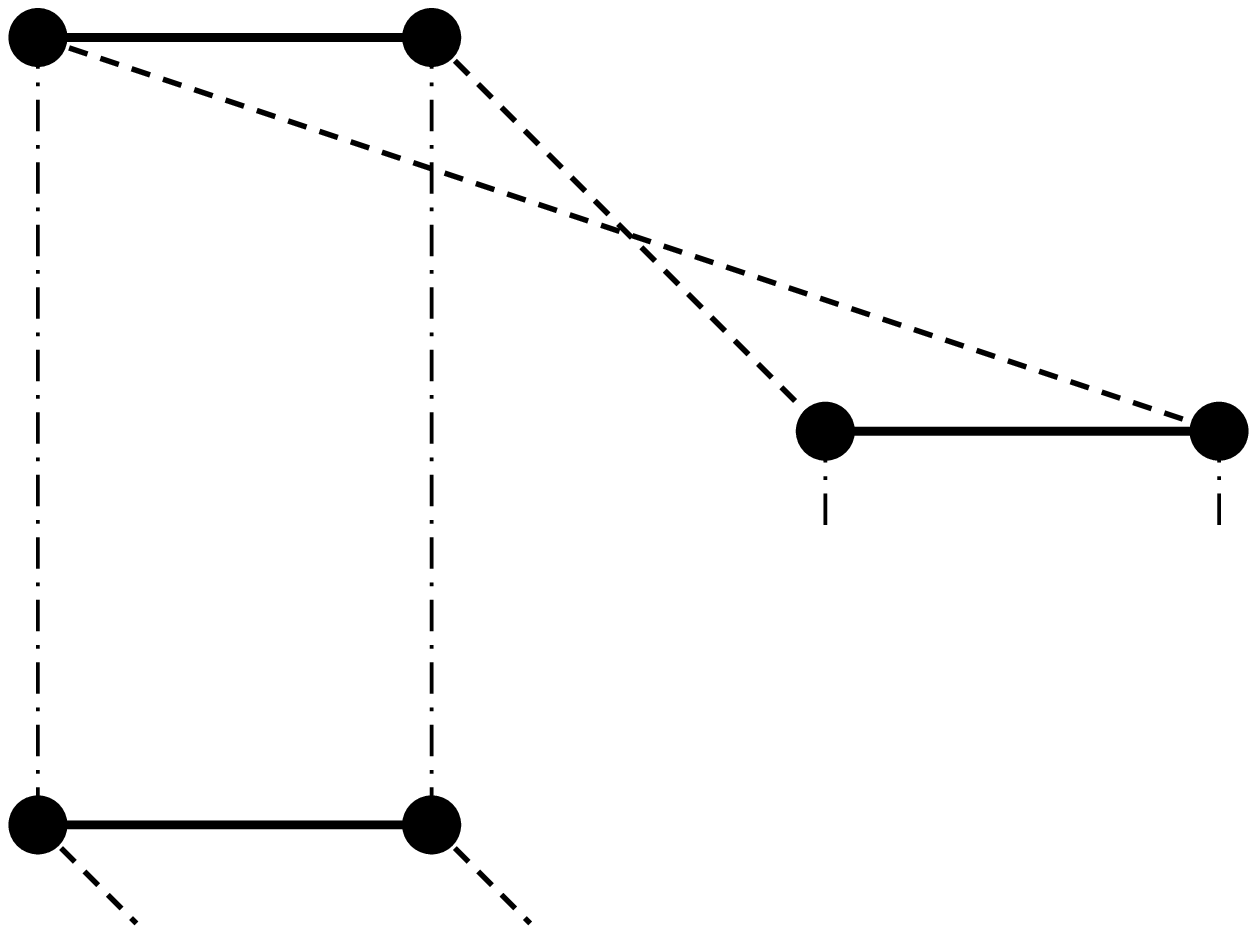}
    \end{tabular} &
    \begin{tabular}{c}
    \\
    \includegraphics[width=.11\textwidth]{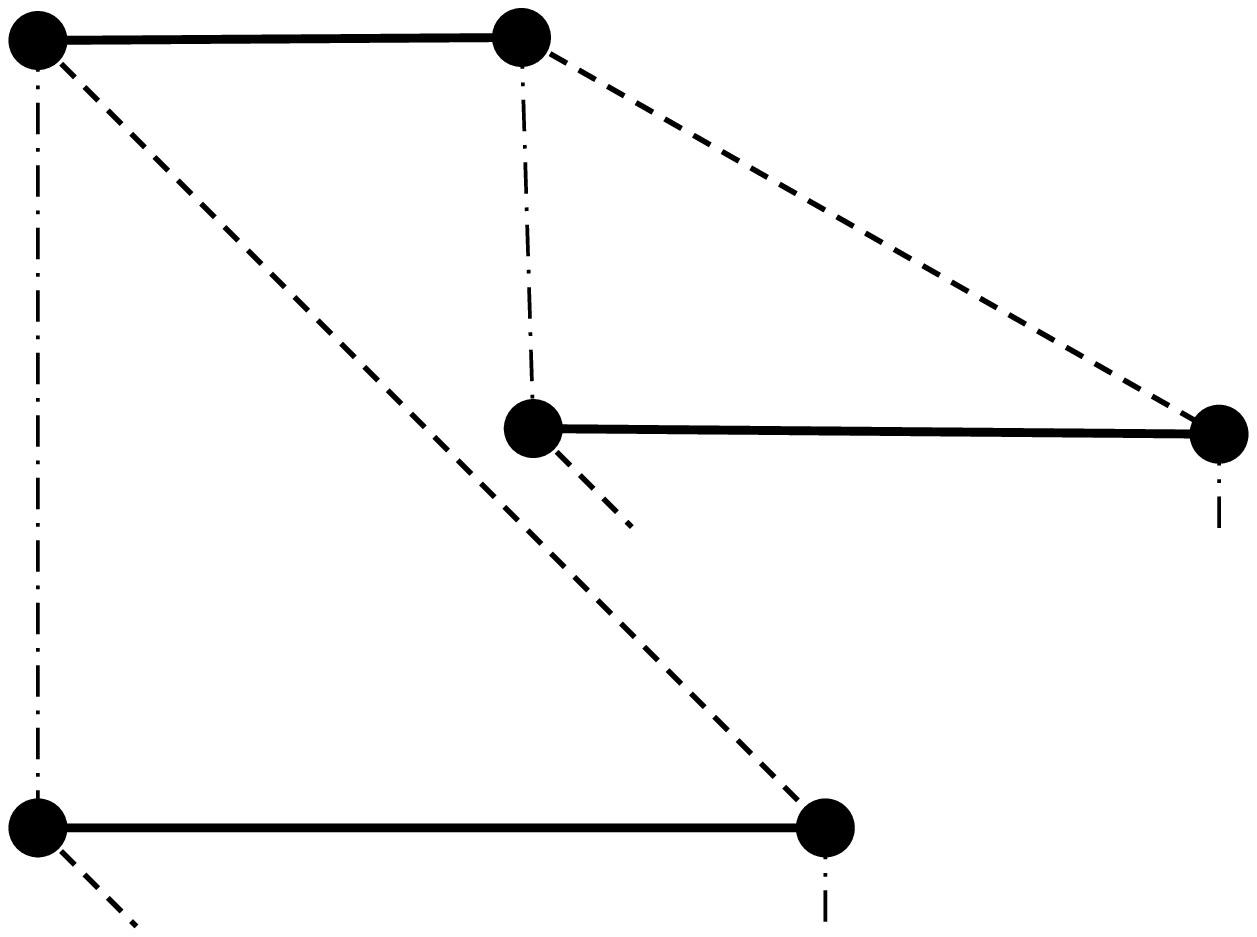}
    \end{tabular} &
    \begin{tabular}{c}
    \\
    \includegraphics[width=.11\textwidth]{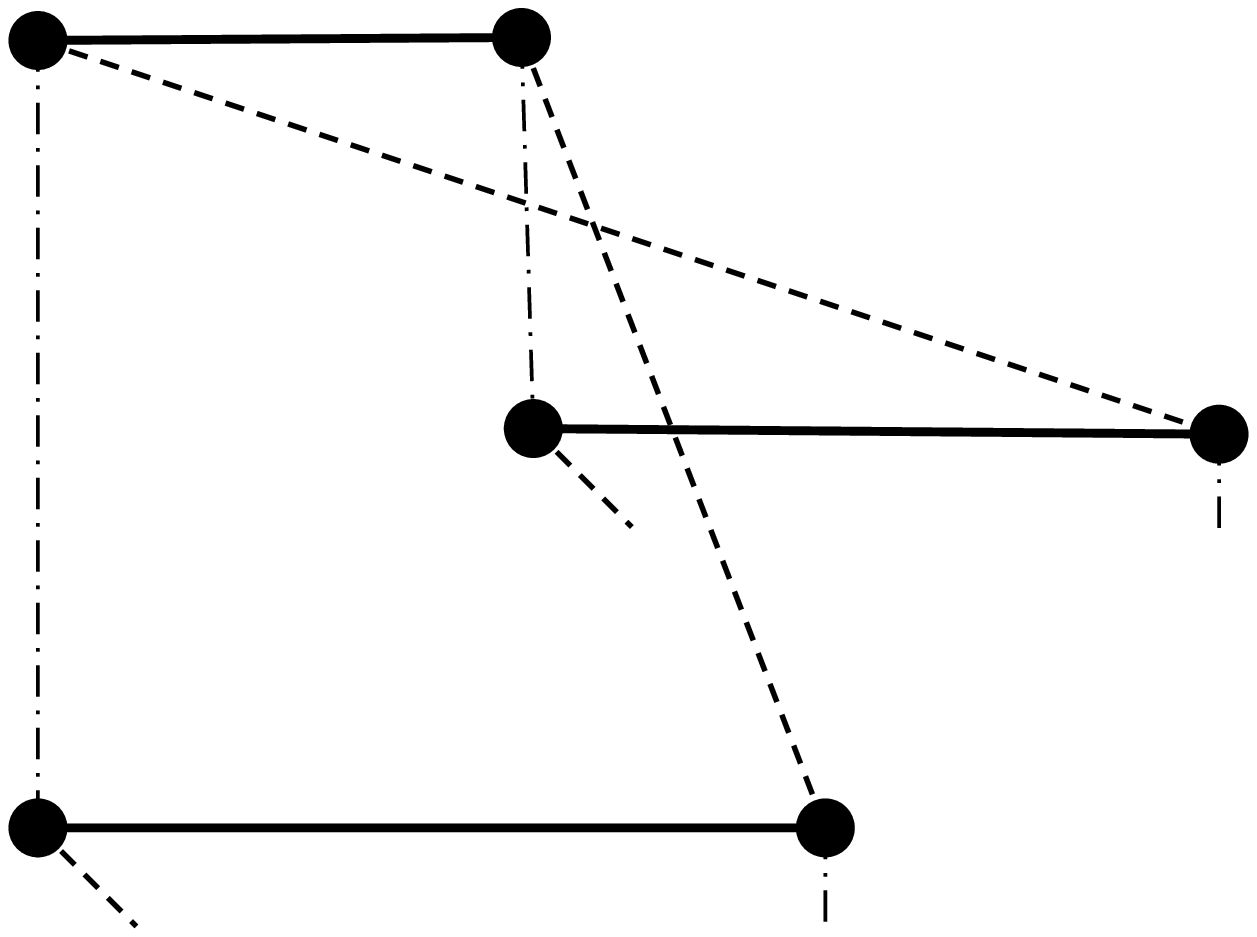}
    \end{tabular} &
    \begin{tabular}{c}
    \\
    \includegraphics[width=.11\textwidth]{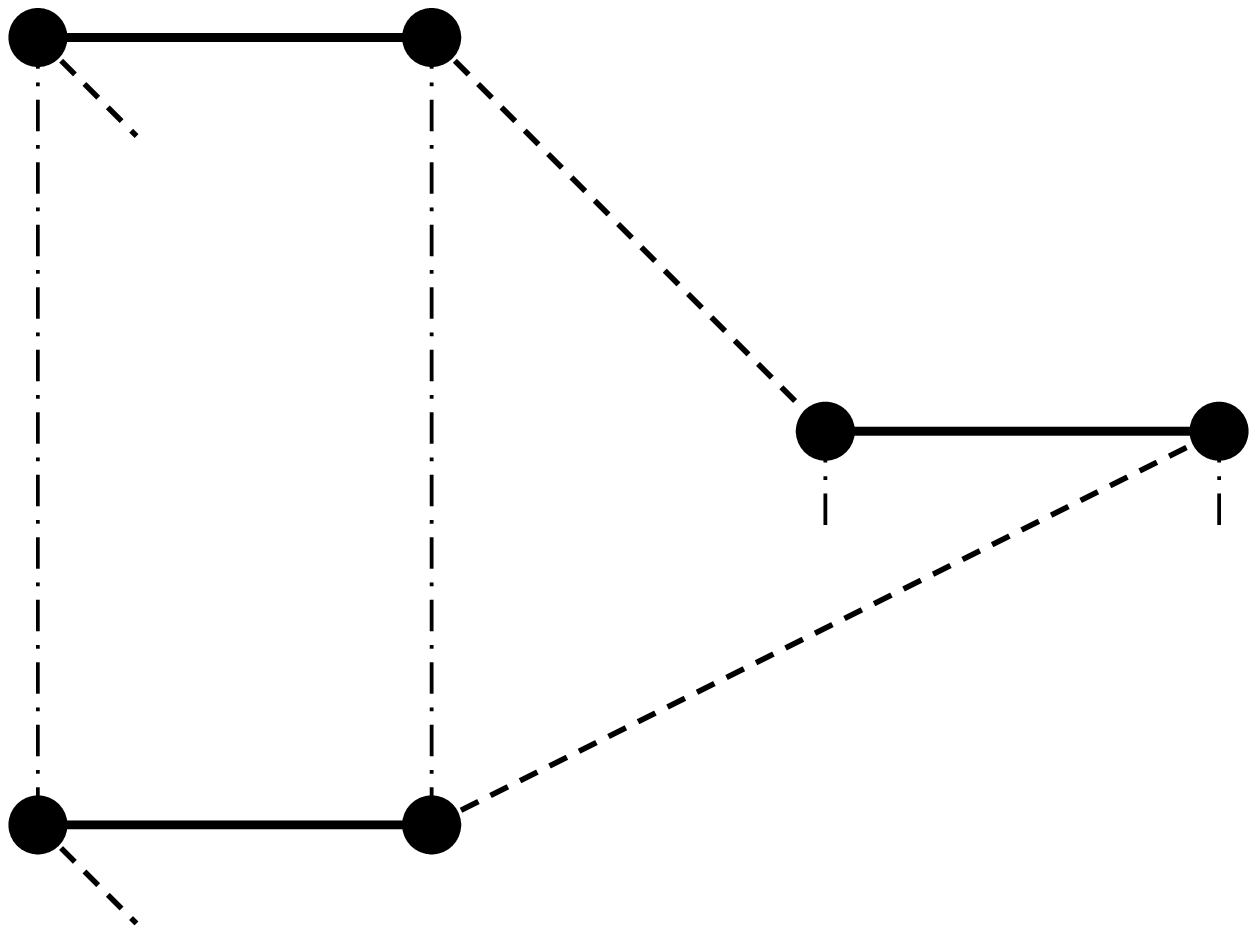}
    \end{tabular}
    \\
    $\mathcal{S}_{6,43}$ & $\mathcal{S}_{6,44}$ & $\mathcal{S}_{6,45}$ & $\mathcal{S}_{6,46}$ & $\mathcal{S}_{6,47}$ & $\mathcal{S}_{6,48}$\\
  \end{tabular}
  \caption{Classification of seminets with point rank up to six (I).}\label{Fig3}
\end{center}}
\end{figure}

\begin{figure}{\scriptsize
\begin{center}
  \begin{tabular}{cccccccc}
    \begin{tabular}{c}
    \\
    \includegraphics[width=.11\textwidth]{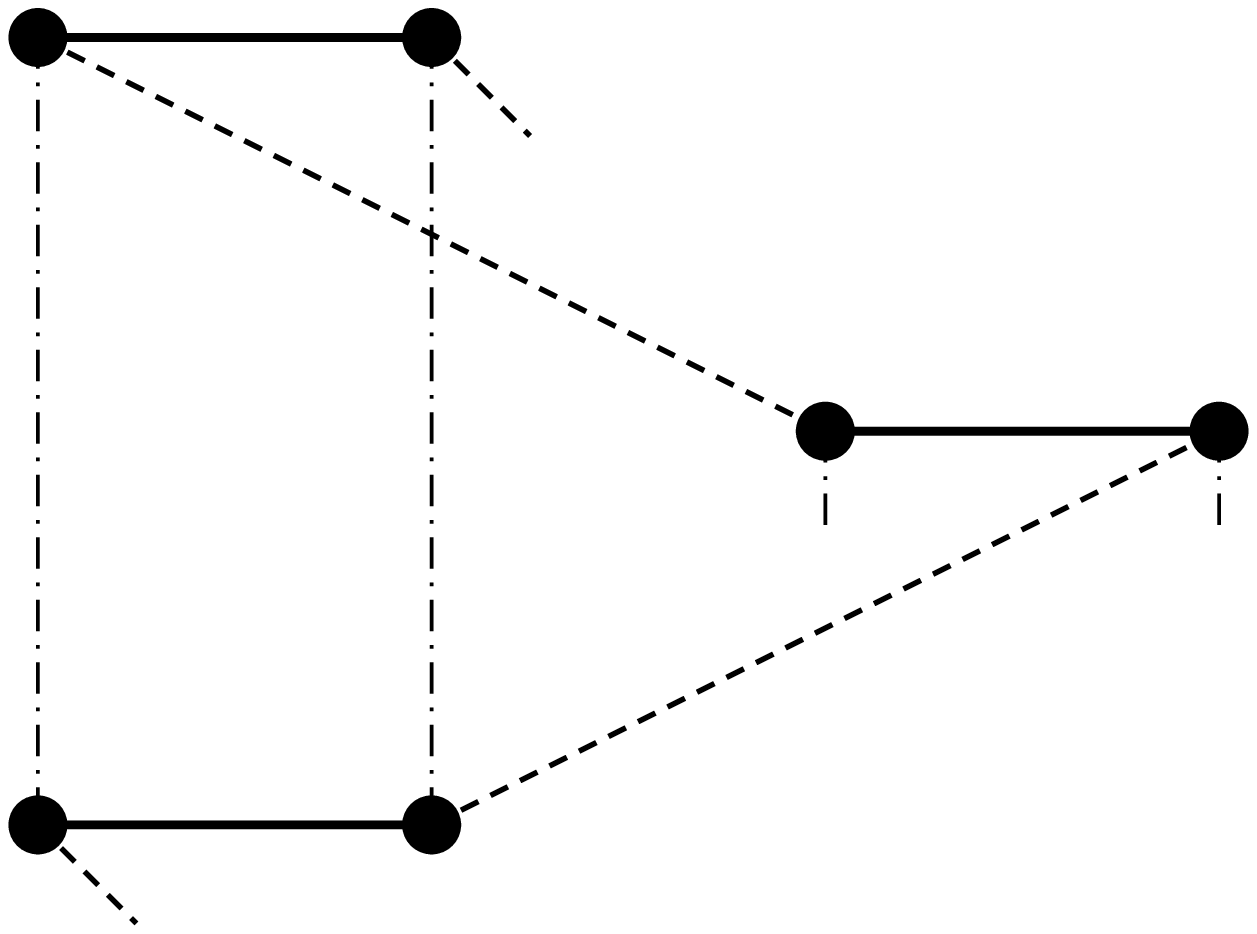}
    \end{tabular} &
    \begin{tabular}{c}
    \\
    \includegraphics[width=.11\textwidth]{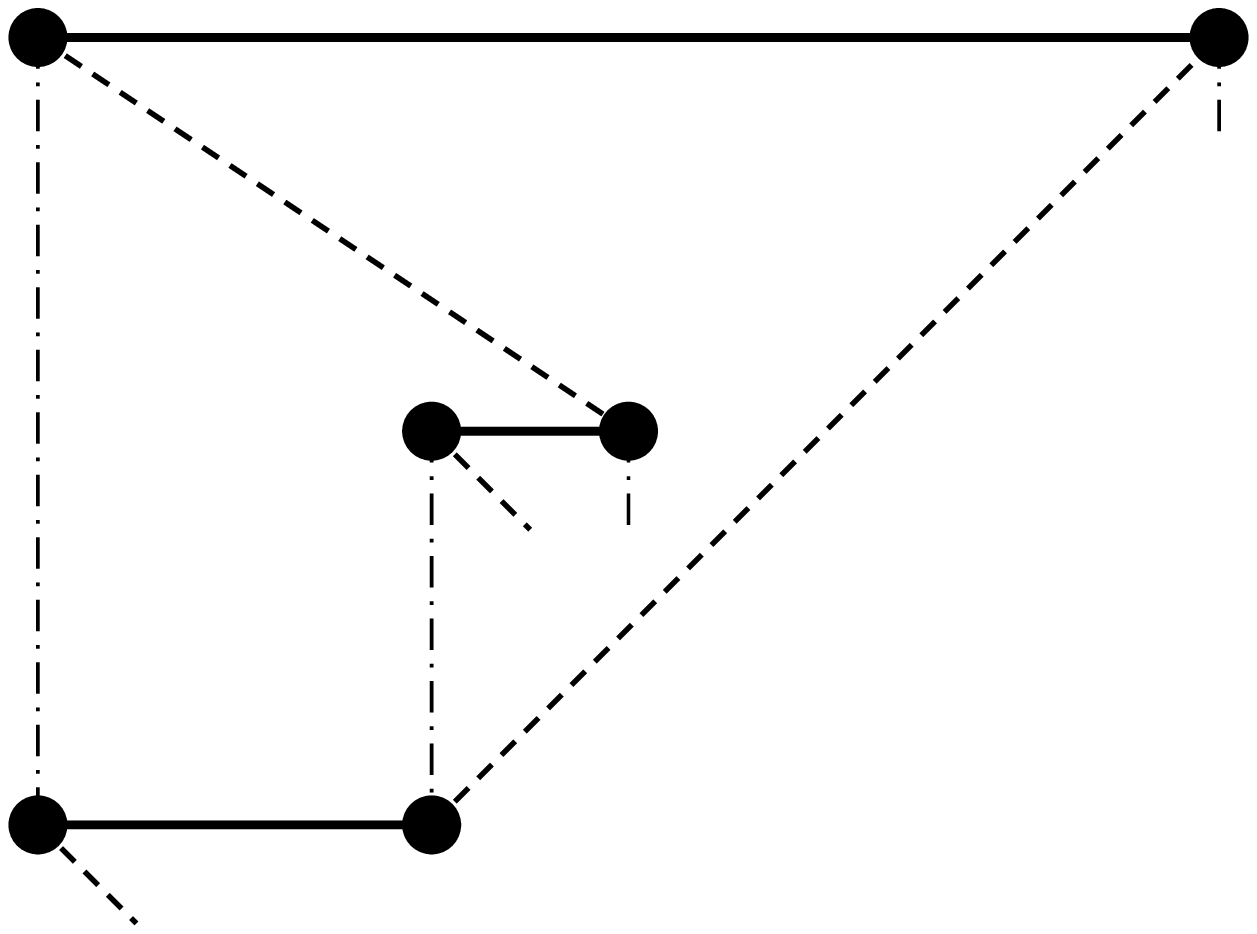}
    \end{tabular} &
    \begin{tabular}{c}
    \\
    \includegraphics[width=.11\textwidth]{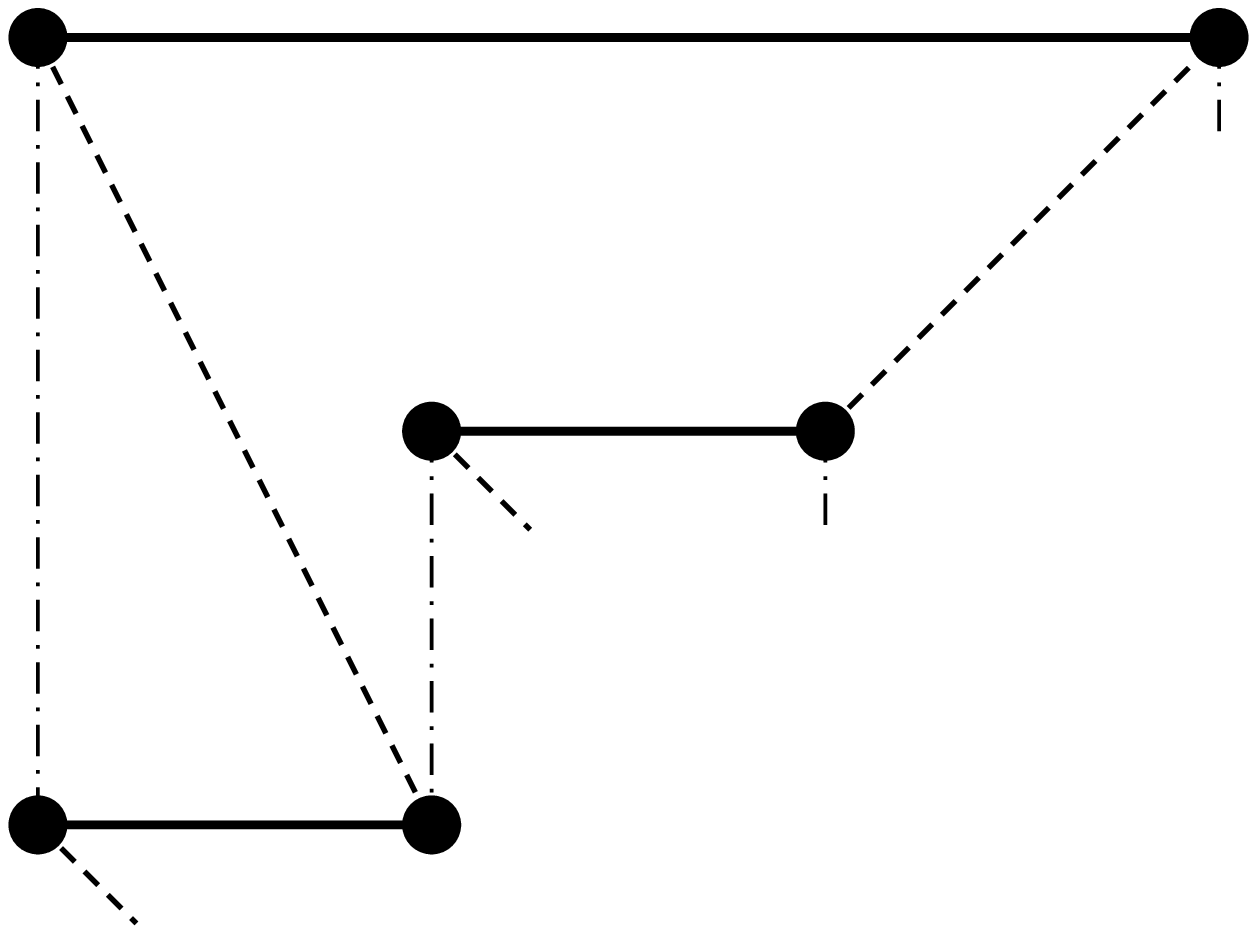}
    \end{tabular} &
    \begin{tabular}{c}
    \\
    \includegraphics[width=.11\textwidth]{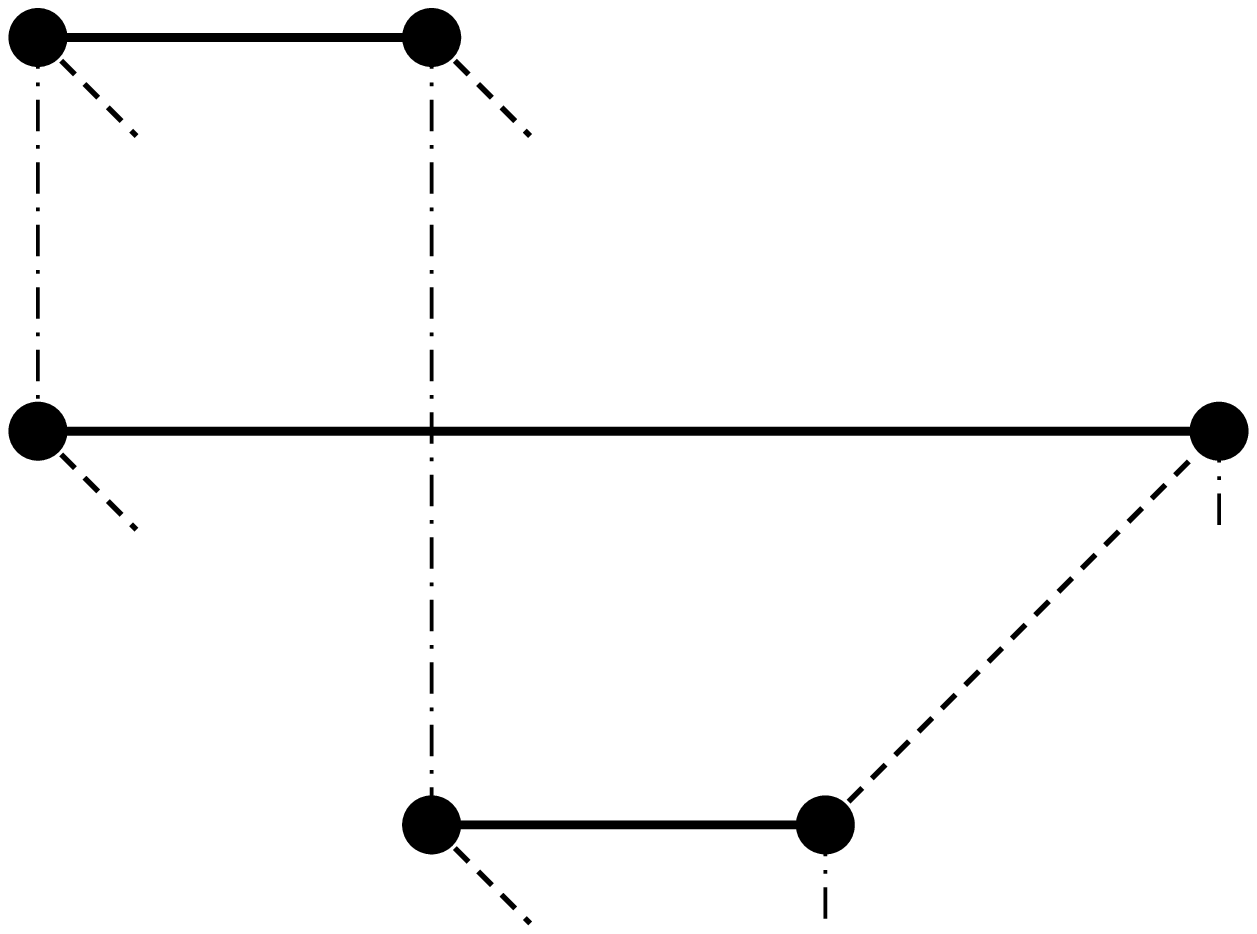}
    \end{tabular} &
    \begin{tabular}{c}
    \\
    \includegraphics[width=.11\textwidth]{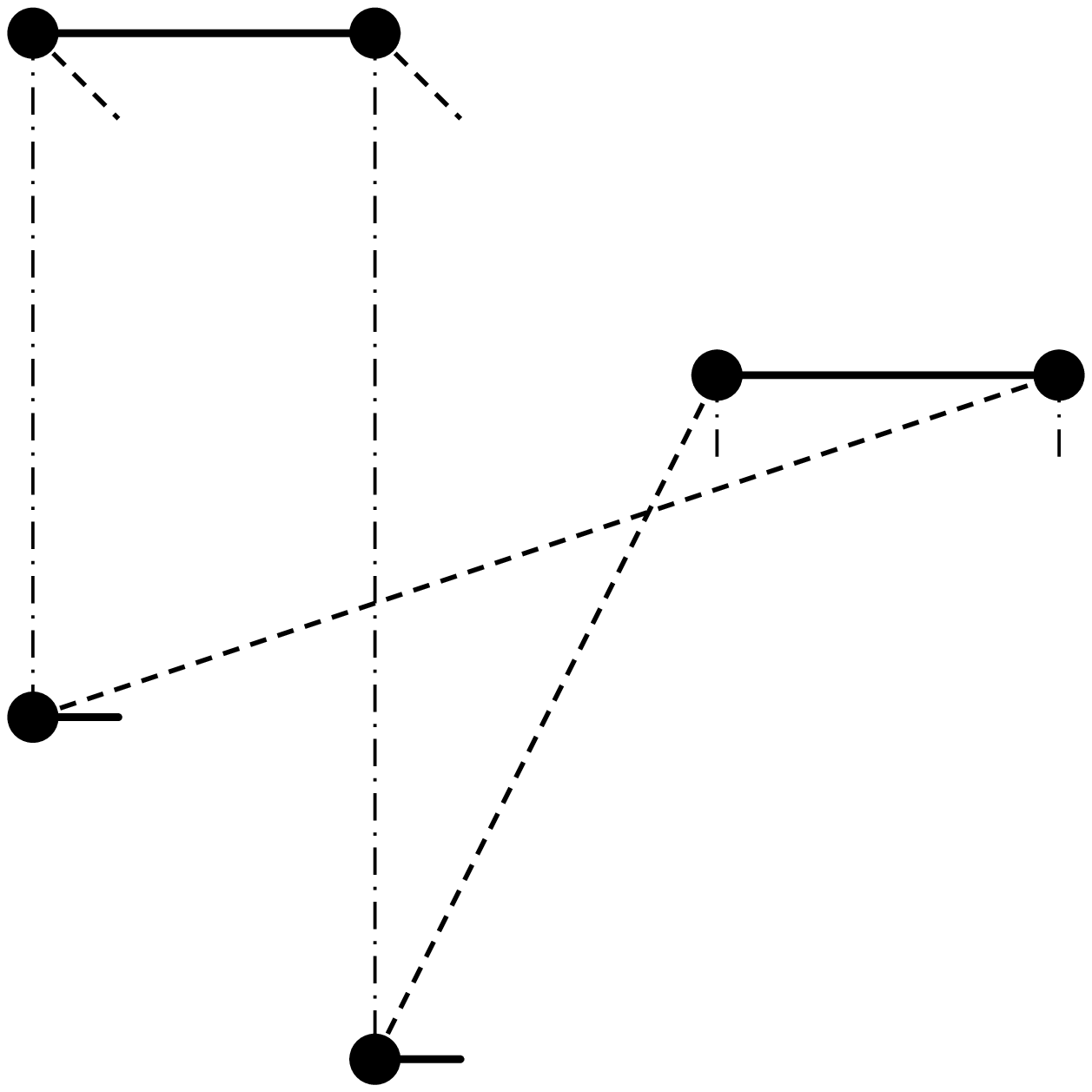}
    \end{tabular} &
    \begin{tabular}{c}
    \\
    \includegraphics[width=.11\textwidth]{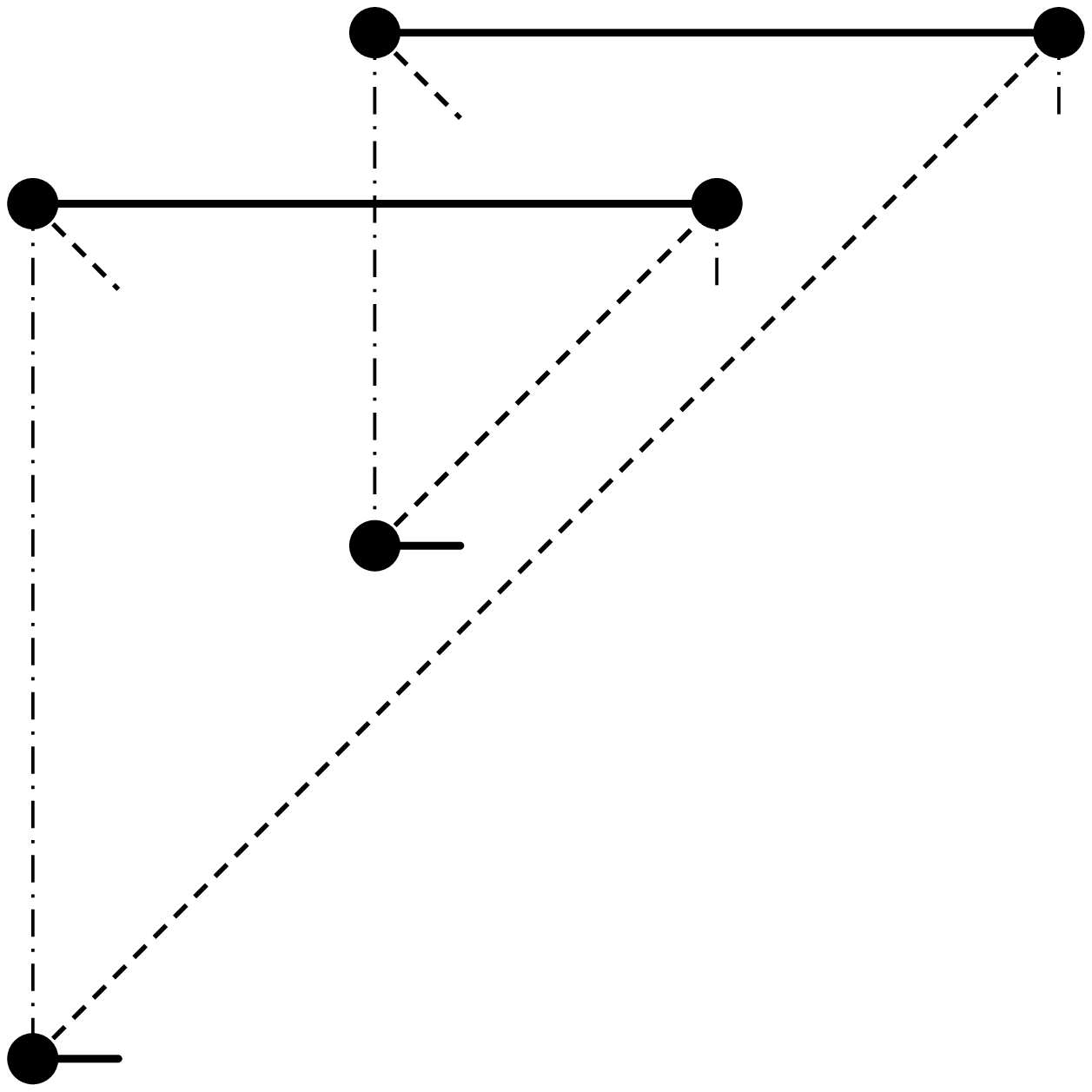}
    \end{tabular} &
    \begin{tabular}{c}
    \\
    \includegraphics[width=.11\textwidth]{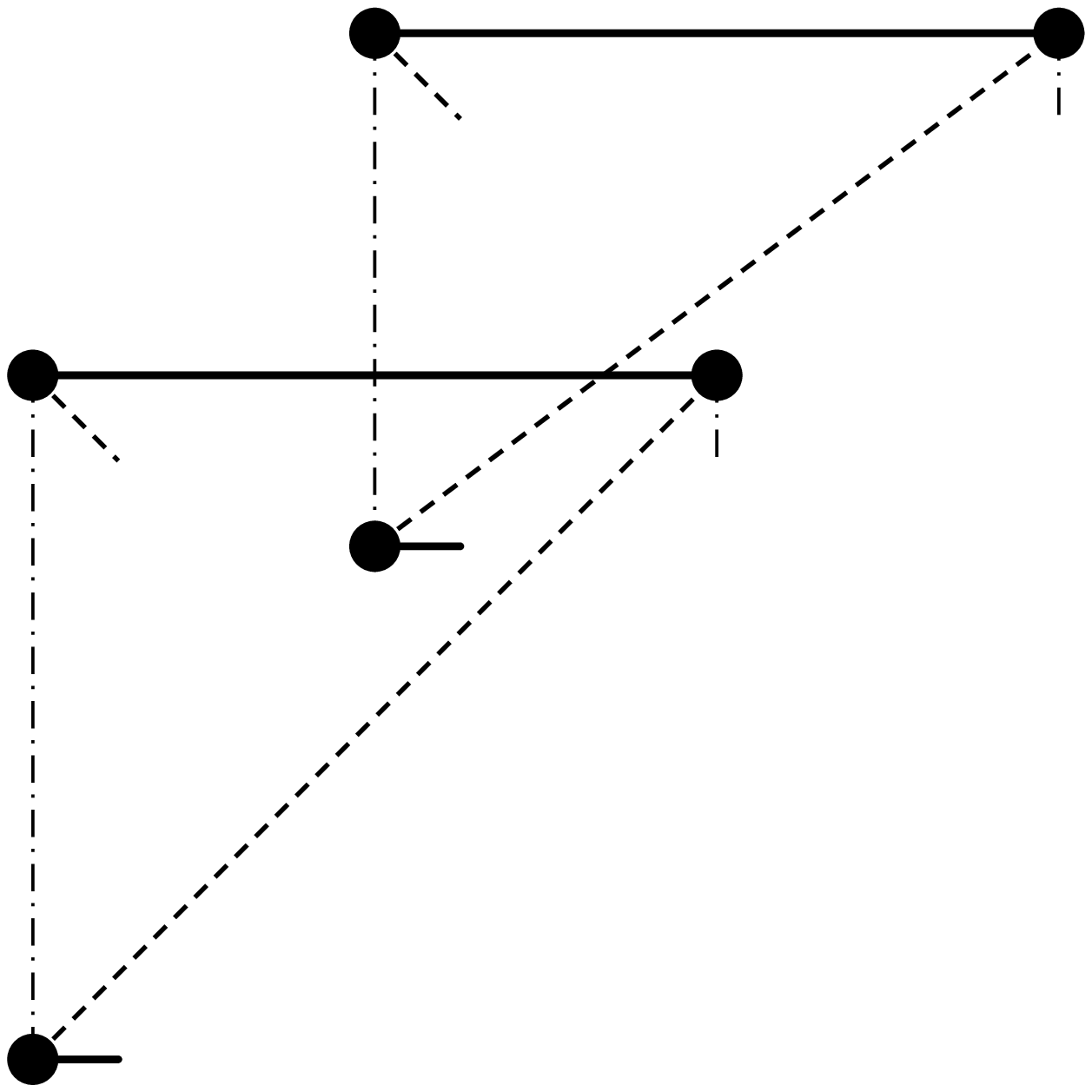}
    \end{tabular}
    \\
    $\mathcal{S}_{6,50}$ & $\mathcal{S}_{6,51}$ & $\mathcal{S}_{6,52}$ & $\mathcal{S}_{6,53}$ & $\mathcal{S}_{6,54}$ & $\mathcal{S}_{6,55}$
  \end{tabular}
  \caption{Classification of seminets with point rank up to six (II).}\label{Fig4}
\end{center}}
\end{figure}

Theorem \ref{thm_reg} has been implemented in the procedure {\em PLRCS} in {\em pls.lib} in order to determine in Table \ref{table6} the distribution of regular partial Latin squares of order up to $8$ according to their structures and main classes. This distribution is equivalent to that of seminets with point rank up to eight. A census of the main classes of seminets with point rank up to six is exposed in Figures \ref{Fig3} and \ref{Fig4}, where we can observe in particular the four configurations whose existence were already established by Havel \cite{Havel1985}: the {\em Fano configurations} $\mathcal{S}_{4,1}$ and $\mathcal{S}_{6,2}$, the {\em shattered Desargues configuration} $\mathcal{S}_{6,32}$ and the {\em Thomsen configuration} $\mathcal{S}_{6,33}$. Havel also determined the three configurations with point rank seven: the {\em hexagonal configuration} $\mathcal{H}$, the {\em first hybrid configuration} $\mathcal{C}_1$ and the {\em second hybrid configuration} $\mathcal{C}_2$. They correspond to the three main classes of partial Latin squares of type $(32^2,32^2,32^2)$ in Table \ref{table6}.
$${\scriptsize
\renewcommand{\tabcolsep}{2pt}
\begin{tabular}{ccccc}
\begin{tabular}{|c|c|c|}\hline
   1 & 2 & 3 \\ \hline
   2 & \  & 1 \\ \hline
   3 &  1 & \ \\ \hline
\end{tabular} & &
\begin{tabular}{|c|c|c|}\hline
   1 & 2 & 3 \\ \hline
   2 & 1  & \ \\ \hline
   3 & \ & 1 \\ \hline
\end{tabular} & & \begin{tabular}{|c|c|c|}\hline
   1 & 2 & 3 \\ \hline
   2 & 1  & \ \\ \hline
   3 & \ & 2 \\ \hline
\end{tabular}\\ $\mathcal{H}$ & & $\mathcal{C}_1$ & & $\mathcal{C}_2$\\
\end{tabular}}$$
Shortly after, Lyakh \cite{Lyakh1988} determined $21$ configurations with point rank $8$, which can be identified with the following partial Latin squares
$${\scriptsize
\renewcommand{\tabcolsep}{3pt}
\begin{tabular}{ccccccc}
\begin{tabular}{|c|c|c|c|}\hline
   1 & 2 & \ & \ \\ \hline
   3 & 4 & \ & \ \\ \hline
   \ & \ & 1 & 2 \\ \hline
   \  & \  & 3  & 4 \\ \hline
\end{tabular} &
\begin{tabular}{|c|c|c|c|}\hline
   1 & 2 & 3 & 4\\ \hline
   3 &  4 &  \ & \ \\ \hline
   \ & \ & 1 & 2 \\ \hline
   \  & \  & \  & \ \\ \hline
\end{tabular} &
\begin{tabular}{|c|c|c|c|}\hline
   1 & 2 & 3 & 4\\ \hline
   2 & 1 & 4 & 3 \\ \hline
   \ &  \ &  \ & \ \\ \hline
   \  & \  & \  & \ \\ \hline
\end{tabular} &
\begin{tabular}{|c|c|c|c|}\hline
   1 & 2 & 3 & 4\\ \hline
   \ & \ & 4 & 3 \\ \hline
   2 &  1 &  \ & \ \\ \hline
   \  & \  & \  & \ \\ \hline
\end{tabular} &
\begin{tabular}{|c|c|c|c|}\hline
   1 & 2 & 3 & \,\, \\ \hline
   2 & 1 & 4 & \ \\ \hline
   3 & 4 & \ & \  \\ \hline
   \  & \  & \  & \ \\ \hline
\end{tabular} &
\begin{tabular}{|c|c|c|c|}\hline
   1 & 2 & 3 & 4\\ \hline
   2 &  4 &  \ & \ \\ \hline
   \ & \ & 1 & 3 \\ \hline
   \  & \  & \  & \ \\ \hline
\end{tabular} &
\begin{tabular}{|c|c|c|c|}\hline
   1 & 2 & 3 & 4\\ \hline
   \ & 1 & \ & 3 \\ \hline
   4 &  \ &  2 & \ \\ \hline
   \  & \  & \  & \ \\ \hline
\end{tabular}\\
$\mathcal{F}_1$ & $\mathcal{F}_2$ & $\mathcal{F}_3$ & $\mathcal{F}_4$ & $\mathcal{F}_5$ & $\mathcal{F}_6$ & $\mathcal{F}_7$\\
\begin{tabular}{|c|c|c|c|}\hline
   \ & 2 & \ & 4\\ \hline
   4 & \ & 1 & \ \\ \hline
   \ & \ & 2 & 3 \\ \hline
   1  & 3  & \  & \ \\ \hline
\end{tabular} &
\begin{tabular}{|c|c|c|c|}\hline
   \ & 2 & \ & 4\\ \hline
   4 & \ & 1 & \ \\ \hline
   \ & \ & 2 & 3 \\ \hline
   3  & 1  & \  & \ \\ \hline
\end{tabular} &
\begin{tabular}{|c|c|c|c|}\hline
   \ & 2 & \ & 4\\ \hline
   1 & \ & 3 & \ \\ \hline
   \ & \ & 4 & 3 \\ \hline
   2  & 1  & \  & \ \\ \hline
\end{tabular} &
\begin{tabular}{|c|c|c|c|}\hline
   \ & 2 & \ & 4\\ \hline
   3 & \ & 1 & \ \\ \hline
   \ & \ & 4 & 3 \\ \hline
   2  & 1  & \  & \ \\ \hline
\end{tabular} & \begin{tabular}{|c|c|c|c|}\hline
   3 & 2 & 4 & \,\, \\ \hline
   1 & 3 & 2 & \ \\ \hline
   4 & 1 & \  & \ \\ \hline
   \  & \  & \  & \ \\ \hline
\end{tabular} & \begin{tabular}{|c|c|c|c|}\hline
   4 & 1 & 3 & 2\\ \hline
   2 & 3  & 4 & 1 \\ \hline
   \ &  \ &  \ & \ \\ \hline
   \  & \  & \  & \ \\ \hline
\end{tabular} &
\begin{tabular}{|c|c|c|c|}\hline
   3 & 4 & 2 & \,\, \\ \hline
   1 & 2 & 3 & \ \\ \hline
   4 & 1 & \  & \ \\ \hline
   \  & \  & \  & \ \\ \hline
\end{tabular}
\\ $\mathcal{F}_8$ & $\mathcal{F}_9$ & $\mathcal{F}_{10}$ & $\mathcal{F}_{11}$ & $\mathcal{F}_{12}$ & $\mathcal{F}_{13}$ & $\mathcal{F}_{14}$\\
\begin{tabular}{|c|c|c|}\hline
   1 & 3 & 2\\ \hline
   3 & 2 & 1\\ \hline
   2 & 1 & \ \\ \hline
\end{tabular} &
\begin{tabular}{|c|c|c|c|}\hline
   \ & 4 & 2 & 3\\ \hline
   3 & 2 & 1 & \\ \hline
   1 & \ & \ & 4\\ \hline
   \  & \  & \  & \ \\ \hline
\end{tabular} & \begin{tabular}{|c|c|c|c|}\hline
   \ & 2 & 4 & 3\\ \hline
   2 & 1 & 3 & \\ \hline
   1 & \ & \ & 4\\ \hline
   \  & \  & \  & \ \\ \hline
\end{tabular} & \begin{tabular}{|c|c|c|c|}\hline
   \ & 2 & 3 & 4\\ \hline
   1 & 3 & 2 & \\ \hline
   4 & \ & \ & 1\\ \hline
   \  & \  & \  & \ \\ \hline
\end{tabular} & \begin{tabular}{|c|c|c|c|}\hline
   \ & 3 & 4 & 2\\ \hline
   1 & 2 & 3 & \\ \hline
   4 & \ & \ & 1\\ \hline
   \  & \  & \  & \ \\ \hline
\end{tabular} & \begin{tabular}{|c|c|c|c|}\hline
   \ & 3 & 4 & 2\\ \hline
   2 & 1 & 3 & \\ \hline
   4 & \ & \ & 1\\ \hline
   \  & \  & \  & \ \\ \hline
\end{tabular} & \begin{tabular}{|c|c|c|c|}\hline
   \ & 4 & 3 & 2\\ \hline
   3 & 2 & 1 & \\ \hline
   4 & \ & \ & 1\\ \hline
   \  & \  & \  & \ \\ \hline
\end{tabular} \\ $\mathcal{F}_{15}$ & $\mathcal{F}_{16}$ & $\mathcal{F}_{17}$ & $\mathcal{F}_{18}$ & $\mathcal{F}_{19}$ & $\mathcal{F}_{20}$ & $\mathcal{F}_{21}$\\
\end{tabular}}$$
They correspond in Table \ref{table6} to
\begin{enumerate}[i.]
\item The two main classes of type $(4^2,2^4,2^4)$: $\mathcal{F}_3$ and $\mathcal{F}_{13}$.
\item The four main classes of type $(42^2,2^4,2^4)$: $\mathcal{F}_2$, $\mathcal{F}_4$, $\mathcal{F}_6$ and $\mathcal{F}_7$.
\item The main class of type $(3^22,3^22,3^22)$: $\mathcal{F}_{15}$.
\item The three main classes of type $(3^22,3^22,2^4)$: $\mathcal{F}_5$, $\mathcal{F}_{12}$ and $\mathcal{F}_{14}$.
\item The six main classes of type $(3^22,2^4,2^4)$: from $\mathcal{F}_{16}$ to $\mathcal{F}_{21}$.
\item Five of the eight main classes of type $(2^4,2^4,2^4)$: $\mathcal{F}_1$, $\mathcal{F}_8$, $\mathcal{F}_9$, $\mathcal{F}_{10}$ and $\mathcal{F}_{11}$.
\end{enumerate}
The next two main classes of type $(2^4,2^4,2^4)$ complete the list of Lyakh.
$${\scriptsize
\renewcommand{\tabcolsep}{3pt}
\begin{tabular}{cc}
\begin{tabular}{|c|c|c|c|}\hline
   1  & 2  & \  & \ \\ \hline
   \ & \ & 2 & 1 \\ \hline
    3 & \ & 4 & \ \\ \hline
   \ & 4 & \ & 3\\ \hline
\end{tabular} &
\begin{tabular}{|c|c|c|c|}\hline
   1  & 2  & \  & \ \\ \hline
   \ & \ & 3 & 4 \\ \hline
   4 & \ & 2 & \ \\ \hline
   \ & 3 & \ & 1\\ \hline
\end{tabular}\\
$\mathcal{F}_{22}$ & $\mathcal{F}_{23}$
\end{tabular}}$$
The eighth main class of type $(2^4,2^4,2^4)$ is not related to a configuration because there exist non-connected points in the corresponding seminet (see Figure \ref{FigSem8}).

\renewcommand{\tabcolsep}{3pt}

\renewcommand{\tabcolsep}{1pt}

\renewcommand{\tabcolsep}{0.5pt}
\begin{table}[htb]
{\tiny
\caption{Distribution into main classes of the set $\mathcal{R}^{\text{reg}}_{R,C,S}$.}
\label{table6}
\begin{center}
\resizebox{\textwidth}{!}{
\begin{tabular}{cccccccc}
\begin{tabular}{llllrr} \hline
$m$ & $z_R$ & $z_C$ & $z_S$ & $\rho^{\text{reg}}$ & MC\\ \hline
3 & $21$ & $21$ & $21$ & 1 & 1\\
4 & $2^2$ & $2^2$ & $2^2$ & 2 & 1\\
\ & \ & \ & $21^2$ & 4 & 1\\
\ & \ & \ & $1^4$ & 24 & 1\\
\ & \ & $21^2$ & $21^2$ & 4 & 1\\
5 & $32$ & $2^21$ & $2^21$ & 4 & 1\\
\ & \ & \ & $21^3$ & 12 & 1\\
\ & $31^2$ & $31^2$ & $2^21$ & 4 & 1\\
\ & \ & $2^21$ & $2^21$ & 8 & 1\\
\ & $2^21$ & $2^21$ & $2^21$ & 32 & 2\\
\ & \ & \ & $21^3$ & 24 & 1\\
6 & $42$ & $2^21^2$ & $2^21^2$ & 8 & 1\\
\ & $3^2$ & $2^3$ & $2^3$ & 12& 1\\
\ & \ & \ & $2^21^2$ & 36& 2\\
\ & \ & \ & $21^4$ & 144& 1\\
\ & \ & \ & $1^6$ & 720& 1\\
\ & \ & $2^21^2$ & $2^21^2$ & 48& 2\\
\ & \ & \ & $21^4$ & 48& 1\\
\ & $41^2$ & $2^21^2$ & $2^21^2$ & 16& 1\\
\ & $321$ & $321$ & $321$ & 1& 1\\
\ & \ &  & $31^3$ & 6& 1\\
\ & \ &  & $2^3$ & 12& 2\\
\ & \ &  & $2^21^2$ & 20& 4\\
\ & \ &  & $21^4$ & 24& 1\\
\ & \ & $2^3$ & $2^3$ & 36& 1\\
\ & \ & \ & $2^21^2$ & 120& 5\\
\ & \ & \ & $21^4$ & 288& 2\\
\ & \ & $2^21^2$ & $2^21^2$ & 160& 4\\
\ & $2^3$ & $2^3$ & $2^3$ & 144& 2\\
\ & \ & \ & $31^3$ & 72& 1\\
\ & \ &  & $2^21^2$ & 432& 4\\
\ & \ &  & $21^4$ & 1,296& 2\\
\ & \ &  & $1^6$ & 4,320& 1\\
\ & \ & $31^3$ & $31^3$ & 36& 1\\
\ & \ &  & $2^21^2$ & 144& 2\\
\ & \ & $2^21^2$ & $2^21^2$ & 624 & 7\\
\ & \ &  & $21^4$ & 288 & 1\\
\ & $2^21^2$ & $2^21^2$ & $2^21^2$ & 160& 3\\
7 & $43$ & $2^31$  & $2^31$ & 54 & 2\\
\ &   &   & $2^21^3$  & 144 & 2\\
\ &   &   & $21^5$  & 360 & 1\\
\ &   & $2^21^3$  & $2^21^3$  & 144 & 1\\
\ & $421$ & $321^2$  & $321^2$  & 4 & 1\\
\ &  &  & $2^31$  & 36 & 3\\
\ &   &   & $2^21^3$ & 48 & 2\\
\ &   & $31^4$  & $2^31$  & 144 & 1\\
\ &   & $2^31$  & $2^31$  & 162 & 4\\
\ &   &   & $2^21^3$  & 360 & 5\\
\ &   &   & $21^5$  & 360 & 1\\
\ &   & $2^21^3$  & $2^21^3$  & 144 & 1\\
\ & $3^21$ & $32^2$  & $32^2$  & 4 & 1\\
\ &  &  & $321^2$  & 12 & 2\\
\ &   &   &  $31^4$  & 48 & 1\\
\ &  &  & $2^31$  & 72 & 3\\
\ &   &   & $2^21^3$  & 192 & 4\\
\ &   &   & $21^5$  & 480& 1\\
\ &  & $321^2$  & $321^2$  & 24& 2\\
\ &   &   &  $31^4$  & 48& 1\\
\ &  &  & $2^31$  & 120& 5\\
\ &   &   & $2^21^3$  & 144& 3\\
\ &    & $2^31$  & $2^31$  & 612& 6\\
\ &   &   & $2^21^3$  &  1,008& 7\\
\ &   &   & $21^5$  & 720& 1\\
\ &   & $2^21^3$  & $2^21^3$  & 288& 1\\
\ & $32^2$ & $32^2$  & $32^2$  & 16& 3\\
\ &   &   & $321^2$  & 48& 5\\
\ &   &   & $31^4$  & 144& 2\\
\ &   &   & $2^31$  & 192& 7\\
\ &   &   & $2^21^3$  & 720& 12\\
\ &   &   & $21^5$  & 2,640& 5\\
\ &   &   & $1^7$  & 10,080& 1\\
\ &   & $321^2$  & $321^2$  & 112& 9\\
\ &   &   & $31^4$  & 192& 2\\
\ &   &   & $2^31$  & 456& 19\\
\ &  &   & $2^21^3$  & 816& 18\\
\\
\end{tabular}
&  &
\begin{tabular}{llllrr} \hline
$m$ & $z_R$ & $z_C$ & $z_S$ & $\rho^{\text{reg}}$ & MC\\ \hline
7 & $32^2$  & $321^2$  & $21^5$  & 480& 1\\
\ &   & $31^4$  & $2^31$  &  1,008& 4\\
\ &   &   & $2^21^3$  & 288& 1\\
\ &   & $2^31$  & $2^31$  & 1,692& 16\\
\ &   &   & $2^21^3$  & 3,744& 26\\
\ &   &   & $21^5$  & 6,480& 5\\
\ &   & $2^21^3$  & $2^21^3$  & 2,592 & 6\\
\ & $321^2$  & $321^2$  & $321^2$  & 144 & 5\\
\ &   &   & $2^31$  & 684 & 18\\
\ &   &   & $2^21^3$  & 264 & 5\\
\ &   & $31^4$  & $2^31$  & 432 & 2\\
\ &   & $2^31$  & $2^31$  & 2,556 & 21\\
\ &   &   & $2^21^3$  & 2,088 & 15\\
\ & $31^4$  & $2^31$  & $2^31$  & 3,456 & 3\\
\ & $2^31$ & $2^31$  & $2^31$  & 8,478 & 13\\
\ &   &   & $2^21^3$  & 10,152 & 16\\
\ &   & $2^21^3$  & $2^21^3$  & 2,160 & 3\\
8 & $53$ & $2^31^2$ & $2^31^2$ & 144 & 1\\
\ & \ & \ & $2^21^4$ & 288 & 1\\
\ & $4^2$ & $2^4$ & $2^4$ & 216 & 2\\
\ & \ & \ & $2^31^2$ & 528 & 3\\
\ & \ & \ & $2^21^4$ & 2,016 & 3\\
\ & \ & \ & $21^6$ & 8,640 & 1\\
\ & \ & \ & $1^8$ & 40,320 & 1\\
\ & \ & $2^31^2$ & $2^31^2$ & 792 & 4\\
\ & \ & \ & $2^21^4$ & 1,440 & 3\\
\ & \ & \ & $21^6$ & 1,440 & 1\\
\ & \ & $2^21^4$ & $2^21^4$ & 576 & 1\\
\ & $521$ & $321^3$ & $2^31^2$ & 72 & 1\\
\ & \ & $2^31^2$ & $2^31^2$ & 432 & 2\\
\ & \ & \ & $2^21^4$ & 576 & 1\\
\ & $431$ & $32^21$ & $32^21$ & 24 & 4\\
\ & \ & \ & $321^3$ & 72 & 6\\
\ & \ & \ & $31^5$ & 240 & 1\\
\ & \ & \ & $2^4$ & 192 & 4\\
\ & \ & \ & $2^31^2$ & 396 & 17\\
\ & \ & \ & $2^21^4$ & 768 & 8\\
\ & \ & \ & $21^6$ & 720 & 1\\
\ & \ & $321^3$ & $321^3$ & 108 & 2\\
\ & \ & \ & $2^4$ & 720 & 5\\
\ & \ & \ & $2^31^2$ & 720 & 10\\
\ & \ & \ & $2^21^4$ & 288  & 1\\
\ & \ & $31^5$ & $2^4$ & 2,880 & 1\\
\ & \ & \ & $2^31^2$ & 720 & 1\\
\ & \ & $2^4$ & $2^4$ & 864 & 2\\
\ & \ & \ & $2^31^2$ & 2,592 & 10\\
\ & \ & \ & $2^21^4$ & 7,488 & 7\\
\ & \ & \ & $21^6$ & 17,280  & 2\\
\ & \ & $2^31^2$ & $2^31^2$ & 3,744 & 15\\
\ & \ & \ & $2^21^4$ & 3,456 & 7\\
\ & $42^2$ & $3^21^2$ & $3^21^2$ & 8 & 1\\
\ & \ & \ & $32^21$ & 16 & 1\\
\ & \ & \ & $321^3$ & 48 & 1\\
\ & \ & \ & $2^4$ & 192 & 2\\
\ & \ & \ & $2^31^2$ & 336 & 4\\
\ & \ & \ & $2^21^4$ & 576 & 3\\
\ & \ & $32^21$ & $32^21$ & 72 & 8\\
\ & \ & \ & $321^3$ & 240 & 10\\
\ & \ & \ & $31^5$ & 720 & 2\\
\ & \ & \ & $2^4$ & 384 & 4\\
\ & \ & \ & $2^31^2$ & 1,104 & 23\\
\ & \ & \ & $2^21^4$ & 2,880 & 15\\
\ & \ & \ & $21^6$ & 5,760 & 2\\
\ & \ & $321^3$  & $321^3$ & 360 & 4\\
\ & \ & \ & $2^4$ & 1,728 & 6\\
\ & \ & \ & $2^31^2$ & 2,448 & 17\\
\ & \ & \ & $2^21^4$ & 1,728 & 3\\
\ & \ & $31^5$ & $2^4$ & 5,760 & 1\\
\ & \ & \ & $2^31^2$ & 2,880 & 1\\
\ & \ & $2^4$ & $2^4$ & 1,296 & 4\\
\ & \ & \ & $2^31^2$ & 5,184 & 11\\
\ & \ & \ & $2^21^4$ & 19,584 & 12\\
\ & \ & \ & $21^6$ & 69,120 & 3\\
\ & \ & \ & $1^8$ & 241,920 & 1\\
\ & \ & $2^31^2$ & $2^31^2$ & 10,368 & 24\\
\\
\end{tabular}
&  &
\begin{tabular}{llllrr} \hline
$m$ & $z_R$ & $z_C$ & $z_S$ & $\rho^{\text{reg}}$ & MC\\ \hline
8 & $42^2$ & $2^31^2$ & $2^21^4$ & 15,552 & 15\\
\ & \ & \ & $21^6$ & 8,640 & 1\\
\ & \ & $2^21^4$ & $2^21^4$ & 3,456 & 2\\
\ & $3^22$ & $3^22$ & $3^22$ & 4 & 1\\
\ & \ & \ & $3^21^2$ & 8 & 1\\
\ & \ & \ & $32^21$ & 48 & 4\\
\ & \ & \ & $321^3$ & 144 & 4\\
\ & \ & \ & $31^5$ & 480 & 1\\
\ & \ & \ & $2^4$ & 192 & 3\\
\ & \ & \ & $2^31^2$ & 720 & 11\\
\ & \ & \ & $2^21^4$ & 2,640 & 11\\
\ & \ & \ & $21^6$ & 10,080 & 3\\
\ & \ & \ & $1^8$ & 40,320 & 1\\
\ & \ & $3^21^2$ & $3^21^2$ & 16 & 1\\
\ & \ & \ & $32^21$ & 104 & 7\\
\ & \ & \ & $321^3$ & 240 & 5\\
\ & \ & \ & $31^5$ & 480 & 1\\
\ & \ & \ & $2^4$ & 480 & 4\\
\ & \ & \ & $2^31^2$ & 1,032 & 14\\
\ & \ & \ & $2^21^4$ & 1,920 & 7\\
\ & \ & \ & $21^6$ & 1,440 & 1\\
\ & \ & $32^21$ & $32^21$ &  396 & 29 \\
\ & \ & \ & $321^3$ & 1,020 & 43\\
\ & \ & \ & $31^5$ & 2,640 & 6\\
\ & \ & \ & $2^4$ & 1,440 & 15\\
\ & \ & \ & $2^31^2$ & 4,008 & 84\\
\ & \ & \ & $2^21^4$ & 9,792 & 51\\
\ & \ & \ & $21^6$ & 18,720 & 7\\
\ & \ & $321^3$  & $321^3$ & 1,440 & 12\\
\ & \ & \ & $31^5$ & 720 & 1\\
\ & \ & \ & $2^4$ & 4,032 & 14\\
\ & \ & \ & $2^31^2$ & 6,336 & 44\\
\ & \ & \ & $2^21^4$ & 5,184 & 9\\
\ & \ & $31^5$ & $2^4$ & 11,520 & 2\\
\ & \ & \ & $2^31^2$ & 7,200 & 3\\
\ & \ & $2^4$ & $2^4$ & 4,896 & 8\\
\ & \ & \ & $2^31^2$ & 14,832 & 31\\
\ & \ & \ & $2^21^4$ & 46,080 & 25\\
\ & \ & \ & $21^6$ & 146,880 & 6\\
\ & \ & \ & $1^8$ & 483,840 & 1\\
\ & \ & $2^31^2$ & $2^31^2$ & 26,208 & 53\\
\ & \ & \ & $2^21^4$ & 6,912 & 2\\
\ & \ & \ & $21^6$ & 17,280 & 2\\
\ & \ & $2^21^4$ & $2^21^4$ & 6,912 & 2\\
\ & $51^3$ & $321^3$ & $2^31^2$ & 216 & 1\\
\ & \ & $2^31^2$ & $2^31^2$ & 864 & 1\\
\ & $421^2$ & $421^2$ & $32^21$ & 16 & 2\\
\ & \ & \ & $2^4$ & 144 & 2\\
\ & \ & \ & $321^3$ & 24 & 1\\
\ & \ & \ & $2^31^2$ & 192 & 4\\
\ & \ & \ & $2^21^4$ & 96 & 1\\
\ & \ & $3^21^2$ & $3^21^2$ & 16 & 1\\
\ & \ & \ & $32^21$ & 48 & 3\\
\ & \ & \ & $2^4$ & 384 & 3\\
\ & \ & \ & $321^3$ & 96 & 2\\
\ & \ & \ & $2^31^2$ & 432 & 5\\
\ & \ & \ & $2^21^4$ & 192 & 1\\
\ & \ & $32^21$ & $32^21$ & 240 & 19\\
\ & \ & \ & $2^4$ & 960 & 10\\
\ & \ & \ & $41^4$ & 96 & 1\\
\ & \ & \ & $321^3$ & 528 & 22\\
\ & \ & \ & $2^31^2$ & 1,968 & 41\\
\ & \ & \ & $31^5$ & 480 & 1\\
\ & \ & \ & $2^21^4$ & 2,112 & 11\\
\ & \ & $2^4$ & $2^4$ & 2,592 & 4\\
\ & \ & \ & $41^4$ & 576 & 1\\
\ & \ & \ & $321^3$ & 3,168 & 12\\
\ & \ & \ & $2^31^2$ & 8,208 & 16\\
\ & \ & \ & $31^5$ & 5,760 & 2\\
\ & \ & \ & $2^21^4$ & 15,552 & 9\\
\ & \ & \ & $21^6$ & 8,640 & 1\\
\ & \ & $41^4$ & $2^31^2$ & 288 & 1\\
\ & \ & $321^3$ & $321^3$ & 288 & 3\\
\ & \ & \ & $2^31^2$ & 2,160 & 17\\
\ & \ & $2^31^2$ & $2^31^2$ & 9,648 & 21\\
\\
\end{tabular}
&  &
\begin{tabular}{llllrr} \hline
$m$ & $z_R$ & $z_C$ & $z_S$ & $\rho^{\text{reg}}$ & MC\\ \hline
8 & $421^2$ & $2^31^2$ & $2^21^4$ & 3,168 & 4\\
\ & $3^21^2$ & $3^21^2$ & $3^21^2$ & 32 & 1\\
\ & \ & \ & $32^21$ & 192 & 4\\
\ & \ & \ & $2^4$ & 1,248 & 5\\
\ & \ & \ & $41^4$ & 96 & 1\\
\ & \ & \ & $321^3$ & 288 & 2\\
\ & \ & \ & $2^31^2$ & 1,248 & 7\\
\ & \ & \ & $2^21^4$ & 576 & 2\\
\ & \ & $32^21$ & $32^21$ & 800 & 28\\
\ & \ & \ & $2^4$ & 3,648 & 19\\
\ & \ & \ & $41^4$ & 192 & 1\\
\ & \ & \ & $321^3$ & 1,344 & 240\\
\ & \ & \ & $2^31^2$ & 5,184 & 55\\
\ & \ & \ & $31^5$ & 960 & 1\\
\ & \ & \ & $2^21^4$ & 4,608 & 12\\
\ & \ & $2^4$ & $2^4$ & 13,248 & 8\\
\ & \ & \ & $41^4$ & 1,152 & 1\\
\ & \ & \ & $321^3$ & 8,064 & 14\\
\ & \ & \ & $2^31^2$ & 24,480 & 28 \\
\ & \ & \ & $31^5$ & 11,520 & 1\\
\ & \ & \ & $2^21^4$ & 38,016 & 14\\
\ & \ & \ & $21^6$ & 17,280 & 1\\
\ & \ & $41^4$ & $2^31^2$ & 576 & 1\\
\ & \ & $321^3$ & $321^3$ & 576 & 3\\
\ & \ & \ & $2^31^2$ & 4,176 & 15\\
\ & \ & $2^31^2$ & $2^31^2$ & 19,296 & 23\\
\ & \ & \ & $2^21^4$ & 5,184 & 5\\
\ & $32^21$ & $32^21$ & $32^21$ & 2,768 & 69\\
\ & \ & \ & $2^4$ & 9,504 & 59\\
\ & \ & \ & $41^4$ & 720 & 6\\
\ & \ & \ & $321^3$ & 5,328 & 117\\
\ & \ & \ & $2^31^2$ & 18,144 & 206\\
\ & \ & \ & $31^5$ & 8,640 & 11 \\
\ & \ & \ & $2^21^4$ & 26,016 & 77\\
\ & \ & \ & $21^6$ & 15,840 & 5\\
\ & \ & $2^4$ & $2^4$ & 27,072 & 16\\
\ & \ & \ & $41^4$ & 2,304 & 2\\
\ & \ & \ & $321^3$ & 22,176 & 77\\
\ & \ & \ & $2^31^2$ & 62,784 & 110\\
\ & \ & \ & $31^5$ & 48,960 & 9\\
\ & \ & \ & $2^21^4$ & 130,176 & 57\\
\ & \ & \ & $21^6$ & 207,360 & 7\\
\ & \ & $41^4$ & $321^3$ & 432 & 2\\
\ & \ & \ & $2^31^2$ & 2,880 & 5\\
\ & \ & $321^3$ & $321^3$ & 4,078 & 31\\
\ & \ & \ & $2^31^2$ & 19,512 & 137\\
\ & \ & \ & $2^21^4$ & 4,896 & 9\\
\ & \ & $2^31^2$ & $2^31^2$ & 72,576 & 133\\
\ & \ & \ & $31^5$ & 8,640 & 4\\
\ & \ & \ & $2^21^4$ & 47,232 & 42\\
\ & $2^4$ & $2^4$ & $2^4$ & 67,824 & 8\\
\ & \ & \ & $41^4$ & 5,184 & 2\\
\ & \ & \ & $321^3$ & 69,120 & 14\\
\ & \ & \ & $2^31^2$ & 177,120 & 25\\
\ & \ & \ & $31^5$ & 172,800 & 3\\
\ & \ & \ & $2^21^4$ & 475,200 & 20\\
\ & \ & \ & $21^6$ & 1,296,000 & 5\\
\ & \ & \ & $1^8$ & 3,628,800 & 2\\
\ & \ & $41^4$ & $41^4$ & 576 & 1\\
\ & \ & \ & $321^3$ & 3,456 & 2\\
\ & \ & \ & $2^31^2$ & 12,096 & 3\\
\ & \ & \ & $2^21^4$ & 3,456 & 1\\
\ & \ & $321^3$ & $321^3$ & 27,216 & 22\\
\ & \ & \ & $2^31^2$ & 90,720 & 54\\
\ & \ & \ & $31^5$ & 8,640 & 1\\
\ & \ & \ & $2^21^4$ & 58,752 & 10\\
\ & \ & $2^31^2$ & $2^31^2$ & 263,952 & 53\\
\ & \ & \ & $31^5$ & 86,400 & 3\\
\ & \ & \ & $2^21^4$ & 302,400 & 30\\
\ & \ & \ & $21^6$ & 129,600 & 2\\
\ & \ & $2^21^4$ & $2^21^4$ & 51,840 & 4\\
\ & $41^4$ & $2^31^2$ & $2^31^2$ & 4,320 & 2\\
\ & $321^3$ & $321^3$ & $2^31^2$ & 4,752 & 10\\
\ & \ & $2^31^2$ & $2^31^2$ & 36,288 & 24\\
\ & $2^31^2$ & $2^31^2$ & $2^31^2$ & 167,184 & 27\\
\ & \ & \ & $2^21^4$ & 33,696 & 7\\
\end{tabular}
&
\end{tabular}
}
\end{center}
}
\end{table}

\section{Conclusions and further work}

This paper has dealt with the enumeration and classification of partial Latin rectangles and seminets by means of computational algebraic geometry. Both combinatorial structures have been identified with the points of affine varieties defined by zero-dimensional radical ideals of polynomials. Their decompositions into finitely many disjoint subsets, each of them being the zeros of a triangular system of polynomial equations, have emerged as a useful technique to determine the distribution of $r\times s$ partial Latin rectangles based on $[n]$ into isotopic and main classes according to their size and types, for all $r,s,n\leq 6$, and that of non-compressible regular partial Latin squares of order $n\leq 8$. The latter is equivalent to that of seminets with point rank up to eight and has enabled us to complete a classification previously established by Lyakh. General formulas for the number of partial Latin squares of size up to six and a census of all the seminets with at most six points have also been established. A convenient generalization of the polynomial method exposed in this paper to the theory of $k$-seminets and that of non-compressible, regular and mutually regularly orthogonal partial Latin squares developed by U{\v{s}}an \cite{Uvsan1977} is established as further work.

\end{document}